\documentclass[11pt]{article}
\usepackage{amsfonts, amsmath, amssymb, amsgen, amsthm, amscd,
latexsym}
\usepackage[all]{xy}

\def\Diff{\mathop{\rm Diff}\nolimits}
\def\End{\mathop{\rm End}\nolimits}

\def\GL{\mathop{\rm GL}\nolimits}

\def\Sp{\mathop{\rm Sp}\nolimits}

\def\Id{\mathop{\rm Id}\nolimits}

\def\Ad{\mathop{\rm Ad}\nolimits}
\def\ad{\mathop{\rm ad}\nolimits}

\def\Ker{\mathop{\rm Ker}\nolimits}

\def\Tr{\mathop{\rm Tr}\nolimits}
\def\SL{\mathop{\rm SL}\nolimits}

\def\CSp{\mathop{\rm CSp}\nolimits}

\def\Tot{\mathop{\rm Tot}\nolimits}

\def\Cb{{\mathbb C}}

\def\Nb{{\mathbb N}}
\def\Rb{{\mathbb R}}

\def\Zb{{\mathbb Z}}

\def\Ac{{\cal A}}
\def\Bc{{\cal B}}
\def\Ec{{\cal E}}
\def\Fc{{\cal F}}

\def\Hc{{\cal H}}

\def\Ic{{\cal I}}
\def\Lc{{\cal L}}

\def\Pc{{\cal P}}
\def\Rc{{\cal R}}
\def\Sc{{\cal S}}
\def\Uc{{\cal U}}

\def\Kc{{\cal K}}
\def\Lc{{\cal L}}
\def\Cc{{\cal C}}

\def\Tc{{\cal T}}

\def\Dc{{\cal D}}

\def\a{\alpha}
\def\b{\beta}
\def\d{\delta}
\def\D{\Delta}
\def\g{\gamma}
\def\k{\kappa}
\def\G{\Gamma}
\def\lb{\lambda}
\def\Lba{\Lambda}
\def\om{\omega}
\def\Om{\Omega}
\def\s{\sigma}

\def\t{\theta}
\def\z{\zeta}
\def\ve{\varepsilon}
\def\vp{\varphi}

\def\x{\xi}

\def\ab{{\bf a}}
\def\boldb{{\bf b}}
\def\eb{{\bf e}}

\def\xb{{\bf x}}
\def\yb{{\bf y}}

\def\wb{{\bf w}}

\def\Gb{{\bf G}}

\def\Xb{{\bf X}}

\def\bash{\backslash}

\def\fl{\forall}
\def\ify{\infty}

\def\ot{\otimes}

\def\ra{\rightarrow}
\def\longra{\longrightarrow}

\def\sbs{\subset}

\def\rt{\triangleright}

\def\lt{\triangleleft}
\def\cl{\blacktriangleright\hspace{-4pt} < }
\def\al{>\hspace{-4pt}\vartriangleleft}

\def\acl{\blacktriangleright\hspace{-4pt}\vartriangleleft }

\def\hd{\overset{\ra}{\partial}}
\def \vd{\uparrow\hspace{-4pt}\partial}
\def\hs{\overset{\ra}{\sigma}}
\def \vs{\uparrow\hspace{-4pt}\sigma}
\def\hta{\overset{\ra}{\tau}}
\def \vta{\uparrow\hspace{-4pt}\tau}

\def \vB{\uparrow\hspace{-4pt}B}
\def\p{\partial}

\def\0D{\Delta^{(0)}}
\def\1D{\Delta^{(1)}}
\def\Db{\blacktriangledown}

\def\wg{\wedge}
\def\wt{\widetilde}
\def\td{\tilde}
\def\cop{{^{\rm cop}}}

\newcommand{\FA}{\mathfrak{A}}
\newcommand{\FB}{\mathfrak{B}}

\newcommand{\FH}{\mathfrak{H}}

\newcommand{\Fa}{\mathfrak{a}}
\newcommand{\Fg}{\mathfrak{g}}
\newcommand{\Fh}{\mathfrak{h}}

\newcommand{\Fl}{\mathfrak{l}}

\newcommand{\FX}{\mathfrak{X}}
\newcommand{\FY}{\mathfrak{Y}}
\newcommand{\Fn}{\mathfrak{n}}

\newcommand{\Fd}{\mathfrak{d}}
\newcommand{\Fp}{\mathfrak{p}}
\newcommand{\Fs}{\mathfrak{s}}
\newcommand{\FZ}{\mathfrak{Z}}

\def\Gbo{{\bf \Gamma}}
\def\Dbo{{\bf \Delta}}
\def\tpsi{\tilde{\psi}}

\newtheorem{theorem}{Theorem}[section]
\newtheorem{remark}[theorem]{Remark}
\newtheorem{proposition}[theorem]{Proposition}
\newtheorem{lemma}[theorem]{Lemma}
\newtheorem{corollary}[theorem]{Corollary}

\newtheorem{definition}[theorem]{Definition}

\def\build#1_#2^#3{\mathrel{
\mathop{\kern 0pt#1}\limits_{#2}^{#3}}}
\newcommand{\ps}[1]{~\hspace{-4pt}_{^{(#1)}}}
\newcommand{\ns}[1]{~\hspace{-4pt}_{_{{<#1>}}}}
\newcommand{\sns}[1]{~\hspace{-4pt}_{_{{<\overline{#1}>}}}}

\newcommand{\wid}[1]{\dot{\overbrace{#1}}}

\def\odots{\ot\dots\ot}
\def\wdots{\wedge\dots\wedge}
\newcommand{\nm}[1]{{\mid}#1{\mid}}
\def\one{{\bf 1}}

\numberwithin{equation}{section}

\def\a{\alpha}
\def\b{\beta}

\def\d{\delta}
\def\e{\epsilon}
\def\g{\gamma}
\def\k{\kappa}
\def\i{\iota}
\def\lb{\lambda}
\def\om{\omega}
\def\s{\sigma}
\def\t{\theta}
\def\ve{\varepsilon}

\def\vp{\varphi}

\def\z{\zeta}

\def\D{\Delta}
\def\G{\Gamma}

\def\Om{\Omega}

\def\Ph{\Phi}

\def\fl{\forall}
\def\ify{\infty}

\def\ot{\otimes}
\def\part{\partial}

\def\sbs{\subset}

\def\wdg{\wedge}

\newcommand{\ie}{{\it i.e.\/}\ }
\newcommand{\eg}{{\it e.g.\/}\ }
\newcommand{\cf}{{\it cf.\/}\ }

\def\ra{\rightarrow}
\def\longra{\longrightarrow}

\def\text{\hbox}

\def\bash{\backslash}

\def\fl{\forall}
\def\ify{\infty}

\def\ot{\otimes}

\def\ra{\rightarrow}

\def\sbs{\subset}

\def\wdg{\wedge}
\def\wt{\widetilde}

\def\Ad{\mathop{\rm Ad}\nolimits}

\def\Diff{\mathop{\rm Diff}\nolimits}

\def\End{\mathop{\rm End}\nolimits}

\def\GL{\mathop{\rm GL}\nolimits}

\def\Id{\mathop{\rm Id}\nolimits}
\def\exp{\mathop{\rm exp}\nolimits}
\def\mod{\mathop{\rm mod}\nolimits}

\def\Ker{\mathop{\rm Ker}\nolimits}

\def\KL{\Kc_\Lc}

\def\build#1_#2^#3{\mathrel{
\mathop{\kern 0pt#1}\limits_{#2}^{#3}}}

\numberwithin{equation}{section}
\begin{document}
\title{\bf Hopf algebras of primitive Lie pseudogroups and
Hopf cyclic cohomology}
\author{
\begin{tabular}{cc}
Henri Moscovici \thanks{Research
    supported by the National Science Foundation
    award no. DMS-0245481.}~~\thanks{Department of Mathematics,
    Columbus, OH 43210, USA }\quad and \quad  Bahram Rangipour
    \thanks {Department of Mathematics  and   Statistics,
     University of New Brunswick, Fredericton, NB, Canada}
      \end{tabular}}

\date{ \ }

\maketitle
\begin{abstract}
\noindent
We associate to each infinite primitive Lie pseudogroup
a Hopf algebra  of `transverse symmetries', by refining
a procedure due to Connes and the first author
in the case of the general pseudogroup.
The affiliated Hopf algebra can be viewed  as a `quantum group' counterpart of the
infinite-dimensional primitive Lie algebra of the pseudogroup.
It is first constructed via its action on the \'etale groupoid
associated to the pseudogroup, and then realized
as a bicrossed product of a
universal enveloping algebra by a Hopf algebra of regular functions on a formal group.
The bicrossed product structure allows to express its Hopf cyclic cohomology in terms of
a  bicocyclic bicomplex analogous to the Chevalley-Eilenberg complex.
As an application, we compute the relative Hopf cyclic cohomology modulo
the linear isotropy for the Hopf algebra of the general pseudogroup,
and find explicit  cocycle representatives for
the universal Chern classes in
Hopf cyclic cohomology. As another application, we determine
all Hopf cyclic cohomology groups for the Hopf algebra associated to the pseudogroup
of local diffeomorphisms of the line.

\end{abstract}

\section*{Introduction}

The transverse characteristic classes of foliations with holonomy in
a transitive Lie pseudogroup $\G$
of local diffeomorphisms of $\Rb^n$ are most effectively
described in the framework of the Gelfand-Fuks~\cite{GF} cohomology of the Lie
algebra of formal vector fields associated to $\G$
 (\cf \eg Bott-Haefliger~\cite{BottHaef}).
 In the dual, $K$-homological context,
the transverse characteristic  classes of general foliations
have been expressed by Connes and the first author (\cf~\cite{cm2, cm3, cm4})
in terms of the Hopf cyclic cohomology
of a Hopf algebra $\Hc_n$
canonically associated to the group $\Diff \Rb^n$.

In this paper we construct similar Hopf algebras for
all classical groups of diffeomorphisms, or equivalently  for the
infinite primitive Lie-Cartan pseudogroups~\cite{ECar} of
local $C^\infty$-diffeomorphisms.
The Hopf algebra $\Hc_\Pi$ associated to such a pseudogroup $\Pi$
 can be regarded as a `quantum group' analog of the infinite dimensional
primitive Lie algebra of $\Pi$ (\cf
Singer-Sternberg~\cite{SingStern}, Guillemin~\cite{Guill}). It is
initially constructed via its tautological action on the \'etale
groupoid associated to $\Pi$, and is then reconstructed, in a manner
reminiscent of a `quantum double', as the bicrossed product of a
universal enveloping algebra by a Hopf algebra of regular functions
on a formal group. In turn, the bicrossed product structure is
employed to reduce the computation of the Hopf cyclic cohomology of
$\Hc_\Pi$ to that a bicocyclic bicomplex analogous to the
Chevalley-Eilenberg complex. This apparatus is then applied to
compute the relative Hopf cyclic cohomology of $\Hc_n$ modulo $\Fg
\Fl_n$. We actually find explicit cocycles representing
 the Hopf cyclic analogues of the universal
Chern classes. In the case of $\Hc_1$, the improved technique allows
us to refine our previous computations~\cite{mr} and completely
determine the non-periodized Hopf cyclic cohomology of the Hopf
algebra $\Hc_1$ affiliated with the pseudogroup of local
diffeomorphisms of the line.
\medskip

We now give a brief outline of the main results.
Our construction of the Hopf algebra associated to a primitive Lie
pseudogroup $\Pi$ is modeled on
that of  $\Hc_n$ in~\cite{cm2}, which in turn was inspired
by a procedure due to  G. I. Kac~\cite{Kac} for producing
non-commutative and non-cocommutative quantum groups
out of  `matched pairs' of finite groups.
It relies on splitting the group $\Diff_\Pi$ of
globally defined diffeomorphisms of type $\Pi$
as a set-theoretical product of two subgroups,
\begin{equation*}
\Diff_\Pi = G_\Pi \cdot N_\Pi , \qquad G_\Pi \cap N_\Pi = \{e\}.
\end{equation*}
For a flat (\ie containing all the translations) primitive pseudogroup $\Pi$,
$G_\Pi$ is the subgroup consisting of the affine transformations of $\Rb^n$ that
are in $\Diff_\Pi$, while $N_\Pi$ is the subgroup consisting of those diffeomorphisms
in  $\Diff_\Pi$ that preserve the origin to order $1$.

The pseudogroup of contact transformations is the only infinite
primitive pseudogroup which is not flat.
In that case, we identify $\Rb^{2n+1}$ with the Heisenberg group $H_n$.
Instead of the vector translations, we let  $H_n$ act on itself by group
 left translations, and define the group of
`affine Heisenberg transformations' $G_\Pi$ as the
semidirect product of $H_n$ by the linear isotropy group  $G^0_\Pi$ consisting of
the linear contact transformations. As factor $N_\Pi$ we take the subgroup
of all contact diffeomorphisms preserving the origin
to order $1$ in the sense of Heisenberg calculus, \ie whose differential
at $0$ is the identity map of
the Heisenberg tangent bundle (\cf \eg~\cite{ponge}).
\smallskip

The factorization $\Diff_\Pi = G_\Pi \cdot N_\Pi$ allows to represent uniquely
any $\phi \in \Diff_\Pi$ as a product
$\phi  = \varphi \cdot \psi$, with $\varphi \in G_\Pi$ and  $\psi \in N_\Pi$.
Factorizing the product of any two elements $\varphi \in G_\Pi$ and  $\psi \in N_\Pi$
in the reverse order,
$\, \psi \cdot \varphi \, = \, (\psi\rt
\varphi)\cdot(\psi \lt \varphi)$,
one obtains a left action $\psi  \mapsto \tpsi (\varphi) :=
\psi \rt \varphi$ of $N_\Pi$ on $G_\Pi$, along with
a right action $\lt$
of $G_\Pi$ on $N_\Pi$. Equivalently, these actions are restrictions of
the natural actions of $ \Diff_\Pi$ on the coset spaces
$\Diff_\Pi \slash N_\Pi \cong G_\Pi$ and $G_\Pi \bash \Diff_\Pi \cong N_\Pi$.
\smallskip

The `dynamical' definition of the Hopf algebra $\Hc_\Pi$ associated to
the pseudogroup $\Pi$ is obtained by means of
 its action on the (discrete) crossed product algebra
 $ \Ac_\Pi = C^\infty (G_\Pi) \rtimes  \Diff_\Pi$, which arises as follows.
One starts with a fixed basis $\{X_i\}_{1\leq i \leq m}$
for the Lie algebra
$\Fg_\Pi$ of $G_\Pi$. Each $X \in \Fg_\Pi$ gives rise to
a left-invariant vector field on  $G_\Pi$,
which is then extended to a linear operator on $\Ac_\Pi$, in the
most obvious fashion:
$ \, X (f \, U_{\phi^{-1}} )\, = \, X (f) \, U_{\phi^{-1}}$, where $\, f \in C^\infty (G_\Pi)$
 and $\, \phi \in \Diff_\Pi$.
One has
\begin{equation*}
 U_{\phi^{-1}} \, X_i \, U_{\phi} \,  =  \,  \sum_{j=1}^m  \G_i^j (\phi)\, X_j\, , \qquad i=1, \ldots , m ,
\end{equation*}
with $ \G_i^j (\phi)  \in C^\infty (G_{\rm cn})$, and we define corresponding
multiplication operators on $\Ac_\Pi$ by taking
\begin{equation*}
  \D_i^j (f \, U_{\phi^{-1}}) =   (\Gbo(\phi)^{-1})_i^j  \, f \, U^\ast_\phi ,
   \quad  \text{where} \quad  \Gbo (\phi) = \big(\G_i^j (\phi) \big)_{1 \leq i,j \leq m} .
\end{equation*}
 As an algebra, $\Hc_\Pi$ is generated by
 the operators $X_k$'s and $\D_i^j$'s. In particular, $\Hc_\Pi$
 contains all iterated commutators
 \begin{equation*}
 \D^j_{i, k_1 \ldots k_r} := [X_{k_r}, \ldots , [X_{k_1}, \D_i^j] \ldots] ,
 \end{equation*}
which are multiplication
operators by the functions
 \begin{equation*}
 \G^j_{i, k_1 \ldots k_r} (\phi) := X_{k_r} \ldots X_{k_1} ( \G_i^j (\phi)) , \qquad \phi \in  \Diff_\Pi .
 \end{equation*}
 For any $a, b \in \Ac (\Pi_{\rm cn})$,  one has
\begin{align*}
&  X_k (a b) \, = \, X_k (a) \, b +  \sum_j \D^j_k (a) \, X_j (b) , \\ \label{leibG}
&  \D_i^j (ab)  \, = \, \sum_k  \D_i^k (a) \, \D_k^j (b) ,
 \end{align*}
and by multiplicativity
every $h \in \Hc (\Pi_{\rm cn})$ satisfies a `Leibniz rule'
 of the form
 \begin{equation*}
  h (a b) \, = \,  \sum h_{(1)} (a) \,
h_{(2)} b)  \, , \qquad \fl \, a , b \in \Ac_\Pi .
\end{equation*}
The  operators $\, \D_{\bullet \cdots \bullet}^{\bullet} \,$
satisfy the following Bianchi-type identities:
\begin{equation*}
 \D_{i, j}^k \, - \,  \D_{j, i }^k \, = \,
\sum_{r, s} c^k_{rs}\, \D_{i}^r \, \D_{j}^s  \ - \, \sum_\ell c_{i j}^\ell \, \D_{\ell}^k ,
\end{equation*}
where $c^i_{j k}$ are the structure constants of the Lie algebra $\Fg_\Pi$.
 \medskip

\begin{theorem} \label{Theorem A}
Let  $\FH_\Pi$ be the abstract Lie algebra generated by
the operators $\{ X_k , \, \D^j_{i, k_1 \ldots k_r} \}$ and their commutation
relations.
\begin{enumerate}
\item The algebra $\Hc_\Pi$ is isomorphic to the quotient of the universal
enveloping algebra
 $\Uc (\FH_\Pi)$ by the ideal $\Bc_\Pi$ generated by the Bianchi
identities.
\item The Leibniz rule determines uniquely
 a coproduct, with respect to which
 $\Hc_\Pi $ a Hopf algebra and $\Ac_\Pi$ an $\Hc_\Pi $-module algebra.
    \end{enumerate}
 \end{theorem}

\medskip

We next describe the \textit{bicrossed product} realization of $\Hc_\Pi$.
Let $\Fc_\Pi$ denote the algebra of functions on $N$ generated by the
jet `coordinates'
\begin{equation*}
\eta^j_{i, k_1 \ldots k_r} (\psi) \, : = \,  \G^j_{i, k_1 \ldots
k_r} (\psi)(e) , \qquad \psi \in N ;
\end{equation*}
the definition is obviously independent of the choice
of basis for $\Fg_\Pi$. Furthermore,
 $\Fc_\Pi$  is a Hopf algebra with coproduct uniquely and well-defined
by the rule
\begin{equation*}
 \D f (\psi_1, \psi_2) \, : = \ f (\psi_1 \circ \psi_2) , \qquad \fl \psi_1 , \psi_2 \in N ,
\end{equation*}
and with antipode
\begin{equation*}
 Sf (\psi): =  f (\psi^{-1}) , \qquad \psi \in N .
\end{equation*}
Now the universal enveloping algebra $ \, \Uc_\Pi = \Uc (\Fg_\Pi)$
can be equipped with a right  $\Fc_\Pi$-comodule coalgebra structure
$\Db:\Uc_\Pi \ra\Uc_\Pi \ot \Fc_\Pi $ as follows.
Let $\{X_I = X_1^{i_1} \cdots X_m^{i_m} \, ; \, i_1, \ldots , i_m \in \Zb^+   \}$
be the PBW basis  of $\Uc_\Pi$ induced by the chosen basis of $\Fg_\Pi$.
Then
\begin{equation*}
 U_{\psi^{-1}} \, X_I \, U_\psi \, = \, \sum_J  \b_I^J (\psi)\, X_J,
\end{equation*}
with $\b_I^J  (\psi)$  in
  the algebra of functions on $G_\Pi$ generated by $ \G^j_{i, K} (\psi)$.
 One obtains
 a coaction $\Db:\Uc_\Pi \ra\Uc_\Pi \ot \Fc_\Pi $ by defining
\begin{equation*}
\Db (X_I)\,  = \, \sum_J  X_J \ot \b_I^J (\cdot)(e) ;
\end{equation*}
again, the definition is independent
of the choice of basis.

 The right action $\lt$ of $G_\Pi$ on $N_\Pi$ induces an action of $G_\Pi$ on $\Fc_\Pi$
  and hence a left action of $\Uc_\Pi$ on $\Fc_\Pi$, that makes  $\Fc_\Pi$
 a  left $\Uc_\Pi$-module algebra.
 \medskip

\begin{theorem} \label{Theorem B}
With the above operations,
 $\Uc_\Pi$ and $\Fc_\Pi$ form a matched pair of Hopf algebras, and
 their bicrossed product $\Fc_\Pi\acl \Uc_\Pi$
 is canonically isomorphic to the Hopf algebra  $\Hc_\Pi^{\cop}$.
\end{theorem}
\medskip

The Hopf algebra $\Hc_\Pi$ also comes equipped with a modular character
$\d = \d_{\Pi}$, extending the
infinitesimal modular character $\,  \d (X) = \Tr (\ad X)$, $\, X \in \Fg_{\Pi}$.
The corresponding module $\Cb_\d$, viewed also as a trivial comodule,
defines a `modular pair in involution', \cf~\cite{cm3}, or a particular case
of an `SAYD module-comodule', \cf~\cite{hkrs1}. Such a datum allows to
specialize Connes'  $Ext_\Lba$-definition~\cite{Cext} of cyclic cohomology
 to the context of Hopf algebras
(\cf~\cite{cm2, cm3}). The resulting Hopf
cyclic cohomology, introduced in~\cite{cm2} and extended to SAYD coefficients
in~\cite{hkrs2},
incorporates both Lie algebra and group
cohomology and provides the appropriate cohomological tool
for the treatment of symmetry in noncommutative geometry.
However, its computation
for general, \ie non-commutative and non-cocommutative,
Hopf algebras poses quite a challenge. In the case of $\Hc_n$,
it has been shown in~\cite{cm2} that $HP^\ast (\Hc_n ; \Cb_\d)$ is canonically isomorphic
to the Gelfand-Fuks cohomology of the Lie algebra $\Fa_n$ of formal vector fields
on $\Rb^n$, result which allowed to transfer the transverse characteristic
classes of foliations from $K$-theory classes into $K$-homology characteristic
classes. There are very few instances of direct calculations so far
(see \eg~\cite{mr}, where the periodic Hopf cyclic cohomology
of several variants of $\Hc_1$ has been directly computed), and a general
machinery for performing such computations is only beginning to emerge.

We rely on the bicrossed product structure of
the Hopf algebra $\Hc_\Pi$ to reduce the computation
of its Hopf cyclic cohomology to that of a simpler, bicocyclic
bicomplex. The latter combines the Chevalley-Eilenberg
 complex of the Lie algebra
 $\Fg = \Fg_\Pi$ with coefficients in $\Cb_\d\ot \Fc^{\ot \bullet}$
 and the coalgebra cohomology complex of $\Fc = \Fc_\Pi$ with
coefficients in $\wg^\bullet \Fg$ :
\begin{align*}
\begin{xy} \xymatrix{  \vdots\ar[d]^{\p_\Fg} & \vdots\ar[d]^{\p_\Fg}
 &\vdots \ar[d]^{\p_\Fg} & &\\
\Cb_\d \ot \wg^2\Fg \ar[r]^{\b_\Fc} \ar[d]^{\p_\Fg}&
 \Cb_\d\ot\Fc\ot \wg^2\Fg  \ar[r]^{\b_\Fc} \ar[d]^{\p_\Fg}
 &\Cb_\d\ot\Fc^{\ot 2} \ot \wg^2\Fg \ar[r]^{~~~~~~~~~\b_\Fc}   \ar[d]^{\p_\Fg}& \hdots&\\
\Cb_\d \ot \Fg \ar[r]^{\b_\Fc} \ar[d]^{\p_\Fg}&  \Cb_\d\ot\Fc\ot \Fg
\ar[r]^{\b_\Fc} \ar[d]^{\p_\Fg}&\Cb_\d\ot\Fc^{\ot 2}
 \ot \Fg  \ar[d]^{\p_\Fg} \ar[r]^{~~~~~~~~~\b_\Fc}& \hdots&\\
\Cb_\d \ot \Cb \ar[r]^{\b_\Fc}&  \Cb_\d\ot\Fc\ot \Cb
\ar[r]^{\b_\Fc}&\Cb_\d\ot\Fc^{\ot 2} \ot \Cb \ar[r]^{~~~~~~~\b_\Fc}
& \hdots& .  }
\end{xy}
\end{align*}

\begin{theorem} \label{Theorem C}
\begin{enumerate}
\item The above bicomplex computes the periodic Hopf cyclic cohomology
$HP^\ast (\Hc_\Pi ; \Cb_{\d})$.

\item There is a relative version of the above bicomplex that computes
the relative periodic Hopf cyclic cohomology
$HP^\ast (\Hc_\Pi , \Uc(\Fh); \Cb_{\d})$,
for any reductive subalgebra $\Fh$ of the linear isotropy Lie algebra
$\Fg^0_\Pi$.
\end{enumerate}
\end{theorem}
\medskip

As the main application in this paper, we compute the  periodic  Hopf cyclic cohomology
of $\Hc_n$ relative to $\Fg \Fl_n$ and find explicit cocycle representatives
for its basis, as described below.

For each partition $\lb =  (\lb_1 \geq  \ldots \geq  \lb_k)$
 of the set $ \{ 1, \dots , p \}$, where $\, 1 \leq p \leq n$,
 we let $\lb \in S_p$ also denote a permutation whose
 cycles have lengths $\lb_1 \geq  \ldots \geq  \lb_k$,  \ie
 representing
 the corresponding conjugacy class $[\lb] \in [S_p] $.
 We then define
\begin{equation*}
C_{p, \lb} := \sum  (-1)^\mu \one \ot \eta^{j_1}_{ \mu(1),j_{{\lb
(1)}}} \wdots \eta^{j_{p}}_{\mu(p),j_{{\lb (p)}}}\ot
X_{\mu(p+1)}\wdots
 X_{\mu(n)} ,
\end{equation*}
where the summation is over  all $\, \mu \in S_n$ and all $ \, 1\le j_1,j_2,\dots , j_p\le n$.

\begin{theorem} \label{Theorem D}
 The cochains $\{ C_{p, \lb} \, ; \, 1 \leq p \leq n , \, \, [\lb] \in [S_p] \}$
 are cocycles and their classes form a basis of
the  group  $HP^\e (\Hc_n , \Uc(\Fg \Fl_n); \Cb_{\d})$,
where $\e \equiv n$ {\em mod} $2$, while
 $HP^{1-\e} (\Hc_n , \Uc(\Fg \Fl_n); \Cb_{\d}) = 0$.
\end{theorem}

The correspondence with the universal Chern classes is obvious.
Let
$$\Pc_n[c_1, \dots c_n]=\Cb[c_1, \dots c_n]/\Ic_n
$$
denote the truncated polynomial ring,  where $\deg(c_j)=2j$, and
 $\Ic_n$ is the ideal generated by the monomials of degree
 $> 2n$.   To each  partition $\lb$ as above, one associates the
degree $2p$ monomial
 \begin{equation*}
c_{p, \lb} :=  c_{\lb_1} \cdots c_{ \lb_k} , \qquad \lb_1 +  \ldots
+ \lb_k = p ;
\end{equation*}
the corresponding classes $\{ c_{p, \lb} \, ; \, 1 \leq p \leq n ,
\, \, \lb \in [S_p] \}$ form a
 basis of the vector space $\Pc_n[c_1, \dots c_n]$.
\medskip

 A second application is the complete determination of the Hochschild cohomology
 of $\Hc_1$ and of the non-periodized
  Hopf cyclic cohomology groups $\, HC^q (\Hc_1 ; \Cb_\d)$, where $\d$ is the
  modular character of the Hopf algebra $\Hc_1$.
 \tableofcontents


\section{Construction via Hopf actions}\label{secHflat}

This section is devoted to the `dynamical' construction of the Hopf algebras
associated to primitive Lie pseudogroups of infinite type. For the sake of
clarity, we start with the case of the general pseudogroup (\cf also~\cite{cm6}),
where the technical details can be handled in the most transparent fashion.
The other cases of flat pseudogroups can be treated in a similar manner.
In order to illustrate the slight adjustments needed to cover them,
we describe in some detail  the Hopf algebras affiliated to the
volume preserving and the symplectic pseudogroups.
On the other hand, the case of
the contact pseudogroup, which is the only non-flat one, requires a certain
change of geometric viewpoint,
namely replacing the natural motions of $\Rb^{2n+1}$ with the
natural motions of the Heisenberg group $H_n$.

 \subsection{Hopf algebra of the general pseudogroup}\label{subsecHn}

Let  $F{\Rb}^n \ra {\Rb}^n$ be the frame bundle on ${\Rb}^n$, which
we identify to ${\Rb}^n \times \GL(n, {\Rb})$ in the obvious way:
the 1-jet at $0 \in {\Rb}^n$ of a germ of a local diffeomorphism
$\phi$ on ${\Rb}^n$ is viewed as
 the pair
\begin{equation}\label{frameid}
 \big(x := \phi (0), \, {\bf y}:= \phi_0^\prime (0)\big) \in {\Rb}^n \times \GL(n, {\Rb}) ,
  \qquad \phi_0 (x) := \phi (x) - \phi (0).
 \end{equation}
 The flat connection on $F{\Rb}^n \ra {\Rb}^n$ is given by the matrix-valued
 $1$-form $\, \om = ( \om^i_j ) \,$ where,
with the usual summation convention,
\begin{equation} \label{trivcon}
 \om^i_j \, :=  \, ({\bf y}^{-1})^i_{\mu} \, d{\bf y}^{\mu}_j \, = \, ({\bf y}^{-1} \, d{\bf y})^i_j
 \, , \qquad i, j =1, \ldots , n \, ,
 \end{equation}
and the canonical form is the vector-valued $1$-form $\, \t = ( \t^k
) \,$,
\begin{equation}  \label{can}
 \t^k \, :=  \, ({\bf y}^{-1})^k_{\mu} \, dx^{\mu} \, = \, ({\bf y}^{-1} \, dx)^k
 \, , \qquad k =1, \ldots , n \,
 \end{equation}
 The basic horizontal vector fields for the above connection are
\begin{equation}\label{horiz}
 X_k \, = \, y_k^{\mu} \, \part_{\mu} \, , \quad   k =1, \ldots , n \, ,
 \quad \text{where} \quad \part_{\mu} = \frac{\part} {\part \, x^{\mu}}  \, ,
 \end{equation}
and the fundamental vertical vector fields associated to the
standard basis of ${\Fg \Fl} (n, \Rb )$, formed by the elementary
matrices $\{ E_i^j \, ; \, 1 \leq i, j \leq n \}$, have the
expression
\begin{equation}\label{vert}
 Y_i^j \, = \,  y_i^{\mu} \, \part_{\mu}^j \, , \quad  i, j =1, \ldots , n \, ,
 \quad \text{where}
\quad \part_{\mu}^j  := \frac{\part}{\part \, y_j^{\mu}} .
 \end{equation}
Let $G := {\Rb}^n \rtimes \GL(n, {\Rb})$ denote the group
 of affine motions of ${\Rb}^n. $
 \smallskip

\begin{proposition}\label{leftinv}
The vector fields
 $\{ X_k, Y_i^j \, ;  \,   i, j, k =1, \ldots , n \}$ form a basis  of left-invariant
 vector fields on the group  $G$.
 \end{proposition}

\proof Represent $G $ as the subgroup of $\GL(n+1, \Rb)$ consisting
of the matrices $\displaystyle {\bf a} = \begin{pmatrix}{\bf y} &
x\\ 0 & 1 \end{pmatrix}$, with ${\bf y} \in  \GL(n, {\Rb})$, and $x
\in {\Rb}^n $. Let $\{ e_k, E_i^j \, ;  \,   i, j, k =1, \ldots , n
\}$ be the standard basis of the Lie algebra  ${\Fg}:= {\Rb}^n
\rtimes {\Fg \Fl} (n, \Rb )$, and denote by $\{ \tilde{e}_k,
\tilde{E}_i^j \, ;  \,   i, j, k =1, \ldots , n \}$  the
corresponding left-invariant vector fields. By definition, at the
point ${\bf a}$,
 $\tilde{e}_k$ is tangent
to the curve
$$  t \mapsto \begin{pmatrix}{\bf y} & x\\ 0 & 1 \end{pmatrix}
\begin{pmatrix}{\bf 1} & t e_k\\ 0 & 1 \end{pmatrix} =
\begin{pmatrix}{\bf y} & t {\bf y} e_k + x\\ 0 & 1 \end{pmatrix}
$$
and therefore coincides with $ X_ k \, = \, \sum_\mu y_k^{\mu} \,
\part_{\mu} $, while $\tilde{E}_i^j $ is tangent to
$$  t \mapsto \begin{pmatrix}{\bf y} & x\\ 0 & 1 \end{pmatrix}
\begin{pmatrix} e^{t E_i^j } & 0\\ 0 & 1 \end{pmatrix} =
\begin{pmatrix}{\bf y} e^{t E_i^j } & 0\\ 0 & 1 \ \end{pmatrix}
$$
which is precisely $ Y_i^j \, = \, \sum_\mu  y_i^{\mu} \,
\part_{\mu}^j $.
\endproof

The group of diffeomorphisms
 ${\Gb} := \Diff {\Rb}^n$
acts on $F{\Rb}^n$, by the natural lift of the tautological action
to the frame level:
\begin{equation} \label{frameact}
\wt{\vp} (x, {\bf y}) := \left( \vp (x), {\vp}^{\prime} (x) \cdot
{\bf y} \right) \, , \quad \text{where} \quad {\vp}^{\prime}
(x)^{i}_{j} = \part_{j} \, {\bf \vp}^{i}  (x) \, .
\end{equation}
Viewing here ${\Gb}$ as a discrete group, we  form the crossed
product algebra
$$
   {\Ac} \, : = \, C_c^{\ify} (F{\Rb}^n ) \rtimes {\Gb}  \, .
$$
As a vector space, it is spanned by monomials of the form $\,f \,
U_{\vp}^* \,$, where $\, f \in C_c^{\ify} (F{\Rb}^n) \, $ and  $\,
U_{\vp}^* \,$ stands for $\, \wt{\vp}^{-1} $, while the product is
given by the multiplication rule
\begin{equation}
f_1 \, U_{\vp_1}^* \cdot f_2 \, U_{\vp_2}^* = f_1 (f_2 \circ
\wt{\vp}_1) \, U_{\vp_2 \vp_1}^* \, .
\end{equation}
Alternatively, $\Ac$ can be regarded as the subalgebra of the
endomorphism algebra $\Lc \left(C_c^{\ify} (F{\Rb}^n ) \right) =
\End_{\Cb}\big(C_c^{\ify} (F{\Rb}^n )\big) $, generated by the
multiplication and the translation operators
\begin{eqnarray}
M_f (\xi)   &=& \,  f \, \xi \, , \quad f \in C_c^{\ify} (F{\Rb}^n) \, , \, \xi \in C_c^{\ify} (F{\Rb}^n ) \\
U_{\vp}^* (\xi)  &=& \,  \xi \circ \wt{\vp}  \, , \qquad \vp \in
{\Gb} \, , \, \xi \in C_c^{\ify} (F{\Rb}^n ) \, .
\end{eqnarray}

Since the right action of $\GL(n,\Rb)$ on  $F{\Rb}^n$ commutes with
the action of ${\Gb}$, at the Lie algebra level one has
\begin{equation} \label{Ydisp}
U_{\vp} \, Y_i^j  \, U_{\vp}^* \, = \, Y_i^j  \, , \qquad  \vp  \in
{\Gb} \, .
\end{equation}
This allows to promote the vertical vector fields to derivations of
$\Ac$. Indeed, setting
\begin{equation}
Y_i^j (f \, U_{\vp}^*) \, = \, Y_i^j ( f) \, U_{\vp}^*  \, , \quad f
\, U_{\vp}^* \in {\Ac} \, ,
\end{equation}
the extended operators satisfy the derivation rule
\begin{equation} \label{Yrule}
Y_i^j (a \, b) \, = \, Y_i^j (a) \, b \, + \, a \,Y_i^j (b) \,   ,
\quad a, b \in {\Ac} \, ,
\end{equation}

We also prolong the horizontal vector fields to linear
transformations $X_k  \in \Lc \, ({\Ac}) $, in a similar fashion:
\begin{equation}
X_k (f \, U_{\vp}^* ) = X_k (f) \, U_{\vp}^* \, , \quad f \,
U_{\vp}^* \in {\Ac} \, .
\end{equation}
The resulting operators are no longer ${\Gb}$-invariant. Instead of
(\ref{Ydisp}), they satisfy
\begin{equation} \label{Xphi}
 U_{\vp}^* \, X_k \, U_{\vp} \,  =  \,  X_{k} \,  - \, \g_{jk}^i (\vp) \, Y_i^j  \, ,
\end{equation}
where
\begin{equation} \label{gijk}
\g_{jk}^i (\vp) (x, {\bf y}) \, =\, \left( {\bf y}^{-1} \cdot
{\vp}^{\prime} (x)^{-1} \cdot \part_{\mu} {\vp}^{\prime} (x) \cdot
{\bf y}\right)^i_j \, {\bf y}^{\mu}_k \, .
\end{equation}
Using the left-invariance of the vector fields $X_k$ and
\eqref{Xphi}, or just the explicit formula \eqref{gijk}, one sees
that $ \, \vp \mapsto \g_{jk}^i (\vp) \,$ is a group $1$-cocycle on
${\Gb}$ with values in $C^{\ify} (F{\Rb}^n ) $; specifically,
\begin{equation} \label{gcocy}
\g_{jk}^i (\vp \circ \psi)  \, =\,  \g_{jk}^i (\vp) \circ \tpsi \, +
\, \g_{jk}^i (\psi) , \qquad \fl \, \vp, \psi \in \Gb.
\end{equation}

As a consequence of \eqref{Xphi}, the operators $X_k  \in \Lc \,
({\Ac}) $ are no longer derivations of $\Ac$, but satisfy instead a
non-symmetric Leibniz rule:
\begin{equation} \label{Xrule}
X_k(a \, b) \, = \, X_k (a) \, b \, + \, a \,X_k (b) \, + \,
\d_{jk}^i (a) \, Y_{i}^{j} (b)
 \,   , \quad
a, b \in {\Ac} \, ,
\end{equation}
where the linear operators $\, \d_{jk}^i   \in \Lc \, ({\Ac}) $ are
defined by
\begin{equation} \label{d}
\d_{jk}^i (f \, U_{\vp}^*) \, =\, \g_{jk}^i  (\vp) \, f \, U_{\vp}^*
\, .
\end{equation}
Indeed, on taking $\, a=f_1 \, U_{\vp_1}^*$, $\, b = f_2 \,
U_{\vp_2}^*$,  one has
\begin{eqnarray*}
X_k (a \cdot b) &=&  X_k (f_1 \, U_{\vp_1}^* \cdot f_2 \,
U_{\vp_2}^*) \,
 = \, X_k (f_1 \cdot U_{\vp_1}^* \, f_2 \, U_{\vp_1}) \, U_{\vp_2 \vp_1}^* \\ \nonumber
&=& X_k (f_1) \, U_{\vp_1}^* \cdot f_2 \, U_{\vp_2}^* \, + \, f_1 \,
U_{\vp_1}^* \cdot X_k (f_2 \, U_{\vp_2}^*)  \\ \nonumber &+& f_1 \,
U_{\vp_1}^* \cdot (U_{\vp_1} \, X_k \,  \, U_{\vp_1}^* \,  \, - \,
X_k) ( f_2 \, U_{\vp_2}^*) \, ,
\end{eqnarray*}
which together with (\ref{Xphi}) and the cocycle property
\eqref{gcocy} imply (\ref{Xrule}).

The same cocycle property shows that the operators $\, \d_{jk}^i $
are derivations:
\begin{equation} \label{drule}
\d_{jk}^i (a \, b) \, = \, \d_{jk}^i  (a) \, b \, + \, a \,
\d_{jk}^i  (b) \,   , \quad a, b \in {\Ac} \, ,
\end{equation}

The operators $\, \{X_k , \, Y^i_j \} \,$ satisfy the commutation
relations of the group of affine transformations of ${\Rb}^n$:
\begin{eqnarray} \label{aff}
[Y_i^j , Y_k^{\ell}] &=& \d_k^j Y_i^{\ell} - \d_i^{\ell} Y_k^j \, ,
\\ \nonumber [Y_i^j , X_k] &=& \d_k^j X_i \, , \qquad \qquad [X_k ,
X_{\ell}]\, =\, 0 \, .  \nonumber
\end{eqnarray}
The successive commutators of the operators $\, \d_{jk}^i $'s with
the $X_{\ell}$'s yield  new generations of
\begin{equation} \label{highd}
\d_{jk \,  \ell_1 \ldots \ell_r}^i \, := \, [X_{\ell_r} , \ldots
[X_{\ell_1} , \d_{jk}^i] \ldots ] \, ,
\end{equation}
which involve multiplication by higher order jets of diffeomorphisms
\begin{eqnarray} \label{d'}
\d_{jk \,  \ell_1 \ldots \ell_r}^i \, ( f \, U_{\vp}^*) \, &:=& \,
 \g_{jk \,  \ell_1 \ldots \ell_r}^i (\vp)\, f \, U_{\vp}^* \, , \qquad \text{where}  \\
  \g_{jk \,  \ell_1 \ldots \ell_r}^i (\vp)\, &:=&
\, X_{\ell_r} \cdots X_{\ell_1} \big(\g_{jk}^i (\vp)\big) \, .
\nonumber
\end{eqnarray}
Evidently, they commute among themselves:
\begin{equation} \label{abel}
[\d_{jk \,  \ell_1 \ldots \ell_r}^i , \, \d_{j'k' \,  {\ell}'_1
\ldots  {\ell}_r}^{i'}] \, = \, 0 \, .
\end{equation}

The commutators between the $Y_\nu^\lb$'s and $\d_{j k}^i$'s,
which can be easily
obtained from the explicit expression (\ref{gijk}) of the cocycle
$\g$, are as follows:
\begin{equation*} \label{Yd}
 [Y_\nu^\lb , \d_{j k}^i ] \, =\,  \d_j^\lb \ \d_{\nu k}^i  \, + \, \d_k^\lb \ \d_{j \nu}^i
 \, - \, \d_\nu^i \ \d_{j k}^\lb  \, .
\end{equation*}
More generally, one checks by induction the relations
\begin{equation}  \label{Yhighd}
 [Y_\nu^\lb , \d_{j_1 j_2 \,  j_3 \ldots j_r}^i ] \, =\, \sum_{s=1}^r \, \d^\lb_{j_s} \,
 \d^i_{j_1  j_2 \,  j_3  \ldots j_{s-1} \nu j_{s+1} \ldots j_r}
 \, - \, \d_\nu^i \, \d_{j_1 j_2 \,  j_3  \ldots j_r}^\lb  \, .
\end{equation}

The  commutator relations \eqref{aff}, \eqref{highd}, \eqref{abel},
\eqref{Yhighd} show that the subspace $\Fh_n$ of $\Lc (\Ac)$
generated by the operators
\begin{equation} \label{gens}
\{X_k , \, Y^i_j , \, \d_{jk \,  \ell_1  \ldots \ell_r}^i \, ;  \,
i, j, k, \ell_1  \ldots \ell_r  =1, \ldots , n , \, r \in \Nb\}
\end{equation}
forms a {\it Lie algebra}.
 \medskip

 We let
$\Hc_n$ denote the  {\it unital subalgebra} of  $\, \Lc (\Ac)$
generated by
 $\Fh_n$.
 Unlike the codimension $1$ case (cf.~\cite{cm2}),
$\Hc_n$ does not coincide with the universal enveloping algebra $\FA
({\Fh}_n)$ when $n \geq 2$.
 Indeed, first of all the operators $\, \d_{jk \,  \ell_1 \ldots \ell_r}^i \,$
are not all distinct; the order of the first two lower indices or of
the last $r$ indices is immaterial. Indeed, the expression of the cocycle $\g$,
\begin{equation} \label{sgijk}
\g_{jk}^i (\vp) (x, {\bf y}) \, =\, ( {\bf y}^{-1} )^i_\lb \,
({\vp}^{\prime} (x)^{-1})^\lb_\rho
 \, \part_{\mu} \part_\nu {\vp}^{\rho} (x) \, {\bf y}^\nu_j \, {\bf y}^{\mu}_k \, ,
\end{equation}
is clearly symmetric in the indices $j$ and $k$. The symmetry
in the last $r$ indices follows from the definition (\ref{highd})
and the fact that, the connection being flat, the horizontal vector
fields commute. It can also be directly seen from the explicit
formula
\begin{eqnarray} \label{highg}
&\g_{jk \,  \ell_1 \ldots \ell_r}^i (\vp) (x, {\bf y}) \, = \\
\nonumber
 &= \, ({\bf y}^{-1})^i_\lb \part_{\b_r} \ldots \part_{\b_1}
\left(({\vp}^{\prime} (x)^{-1})^\lb_\rho
 \part_{\mu} \part_\nu {\vp}^{\rho} (x) \right) {\bf y}^\nu_j  {\bf y}^{\mu}_k  {\bf y}^{\b_1}_{\ell_1}
 \ldots {\bf y}^{\b_r}_{\ell_r}\, .
\end{eqnarray}
\smallskip

\begin{proposition} \label{Bianchi}
The  operators $\, \d_{\bullet \cdots \bullet}^{\bullet} \,$
satisfy the identities
\begin{equation}  \label{bianchi}
 \d_{j \ell \,  k}^i \, - \,  \d_{j k \,  \ell}^i \, = \,
 \d_{j k}^s \, \d_{s \ell}^i  \ - \, \d_{j \ell}^s \, \d_{s k}^i .
\end{equation}
 \end{proposition}

\proof These Bianchi-type identities are an expression of the fact
 that the underlying connection is flat. Indeed, if
 we let $a, b \in {\Ac}$ and
apply  \eqref{Xrule} we obtain
\begin{eqnarray*}
&X_\ell X_k(a \, b) \, \, =\, X_\ell X_k (a) \, b \, + \, X_k (a)
\,X_\ell (b) \, +
\, \d_{j \ell}^i (X_k (a)) \, Y_{i}^{j}  (b) \\
&+\, X_\ell (a) \, X_k (b) \, + \, a \, X_\ell X_k ( b) \, +
\, \d_{j \ell}^i (a) \, Y_{i}^{j}  (X_k (b)) \\
&+\,X_\ell (\d_{j k}^i (a)) \, Y_{i}^{j}  (b) \, + \, \d_{j k}^i
(a)\, X_\ell (Y_{i}^{j}  (b) ) \, + \, \d_{s \ell}^r (\d_{j k}^i
(a)) \, Y_{r}^{s}  (Y_{i}^{j}  (b)) .
\end{eqnarray*}
Since $\, [X_k, X_\ell]=0$, by antisymmetrizing in $k, \ell$ it
follows that
\begin{eqnarray*}
& \d_{j \ell}^i (X_k (a)) \, Y_{i}^{j}  (b) \, + \, \d_{j \ell}^i
(a) \, Y_{i}^{j}  (X_k (b)) \, + \,
X_\ell (\d_{j k}^i (a)) \, Y_{i}^{j}  (b) \\
&+ \, \d_{j k}^i (a)\, X_\ell (Y_{i}^{j}  (b) ) \,
+ \, \d_{s \ell}^r (\d_{j k}^i (a)) \, Y_{r}^{s}  (Y_{i}^{j}  (b)) \\
&= \d_{j k}^i (X_\ell (a)) \, Y_{i}^{j}  (b) \, + \, \d_{j k}^i (a)
\, Y_{i}^{j}  (X_\ell (b)) \, + \,
X_k (\d_{j \ell}^i (a)) \, Y_{i}^{j}  (b) \\
&+ \, \d_{j \ell}^i (a)\, X_k (Y_{i}^{j}  (b) ) \, + \, \d_{s k}^r
(\d_{j \ell}^i (a)) \, Y_{r}^{s}  (Y_{i}^{j}  (b)).
\end{eqnarray*}
Using the `affine' relations \eqref{aff} and the symmetry of $\,
\d_{j k}^i$ in the lower indices one readily obtains the equation
\eqref{bianchi}.
\endproof
 \medskip

In view of this result, the algebra $\Hc_n$ admits a basis similar
to the Poincar\'e-Birkhoff-Witt basis of a universal enveloping
algebra. The notation needed to specify such a basis involves two
kinds of multi-indices. The first kind are of the form
\begin{equation} \label{1kind}
 I \, = \left\{ i_1 \leq \ldots \leq i_p \, ; \left( ^{j_1}_{k_1} \right) \leq \ldots \leq
 \left( ^{j_q}_{k_q} \right) \right\} \, ,
\end{equation}
while the second kind are of the form $\, \displaystyle K \, = \, \{
\k_1 \leq \ldots \leq \k_r \} $, where
\begin{equation} \label{2kind}
\k_s \, = \, \left( \begin{matrix} &i_s &\qquad \qquad\qquad \hfill
\cr
 j_s &\leq & k_s \leq \, \,
\ell_1^s \leq \ldots \leq \ell_{p_s}^s \hfill \cr\end{matrix}
\right)  \, , \quad s=1, \ldots , r \, ;
\end{equation}
in both cases the inner multi-indices are ordered lexicographically.
We then denote
\begin{equation}\label{pbw}
Z_I = X_{i_1} \ldots X_{i_p} \, Y_{k_1}^{j_1} \ldots Y_{k_q}^{j_q}
\quad \hbox{and} \quad \d_K = \d_{j_1 \, k_1 \,   \ell_1^1 \ldots
\ell_{p_1}^1}^{i_1} \ldots \d_{j_r \, k_r \,   \ell_1^r \ldots
\ell_{p_r}^r}^{i_r} \, .
\end{equation}
\smallskip

\begin{proposition} \label{free}
 The monomials $\, \d_K \, Z_I \,$,
ordered lexicographically, form a linear basis of $\Hc_n $.
  \end{proposition}
\medskip

\begin{proof}
We need to prove that if $\, c_{I, \k} \in \Cb$ are such that
\begin{equation} \label{0}
 \sum_{I, K} \, c_{I, K} \,  \d_K \, Z_I   \, (a) \, =  \, 0 \, ,
\quad \fl \,  a \in \Ac \, ,
\end{equation}
then $\, c_{I, K} = 0$, for any $(I, K)$.

To this end, we evaluate (\ref{0}) on all elements of the form $\, a
\, = \, f \, U_{\vp}^* $ at the point
$$ e \, = \, (x=0, {\bf y} = {\bf I}) \in F {\Rb}^n = {\Rb}^n \times \GL(n,  {\Rb}^n) \, .
$$
In particular, for any fixed but arbitrary $\vp \in {\Gb}$, one
obtains
\begin{equation}  \label{00}
\sum_{I} \left(\sum_{K} \, c_{I, K}  \, \g_{\k} (\vp) ( e)\right)
   (Z_I \, f) ( e) = 0 \,  , \quad \fl \, f \in C_c^{\ify} (F{\Rb}^n ) \, .
\end{equation}
Since the $Z_I$'s form a PBW basis of $\FA ({\Rb}^n \rtimes {\Fg
\Fl} (n, \Rb)) $, which can be viewed as the algebra of
left-invariant differential operators on $F{\Rb}^n $, the validity
of (\ref{00}) for any $f \in C_c^{\ify} (F{\Rb}^n )$ implies the
vanishing for each $I$ of the corresponding coefficient. One
therefore obtains, for any fixed $I$,
\begin{equation}  \label{01}
\sum_{K} \, c_{I, K}  \, \g_{K}  (\vp) ( e) \,
 \, = 0 \,  , \quad \fl \, \vp \in {\Gb} \, .
\end{equation}
To prove the vanishing of all the coefficients, we shall use
induction on the \textit{height} of $\, \displaystyle K \, = \, \{
\k_1 \leq \ldots \leq \k_r \} $; the latter is defined by counting
the total number of horizontal derivatives of its largest
components:
$$
\vert K \vert \, = \, \ell^r_1 + \cdots + \ell^r_{p_r} \, .
$$
\medskip

We start with the case of height $0$, when the identity (\ref{01})
reads
\begin{equation*}
\sum_{K} \, c_{I, K}  \, \g_{j_1 k_1}^{i_1} (\vp)( e)  \cdots
 \g_{j_r k_r}^{i_r} (\vp) ( e ) \,
 \, = 0 \,  , \quad \fl \, \vp \in {\Gb} \, .
\end{equation*}
Let $ {\Gb}_0$ be the subgroup of all $\vp \in {\Gb}$ such that $\vp
(0) = 0$.
 Choosing $\vp$  in the subgroup ${\Gb}^{(2)} (0) \subset {\Gb_0}$
consisting of the diffeomorphisms  whose 2-jet at $0$ is of the form
\begin{equation*}
J_0^2 ( \vp)^i (x) \, = \, x^i \, + \, \frac{1}{2}  \sum_{j, k =1}^n
\xi^{i}_{j k} x^j x^k, \quad  \xi \in \Rb^{n^{3}}, \quad \xi^{i}_{j
k} =   \xi^{i}_{k j} \, ,
\end{equation*}
and using (\ref{sgijk}), one obtains:
\begin{equation*}
\sum_{K} \, c_{I, K}  \, \xi_{j_1 k_1}^{i_1}  \cdots
 \xi_{j_r k_r}^{i_r} \, = 0 , \quad
\xi^{i}_{j k}  \in \Rb^{n^{3}}, \quad  \xi^{i}_{j k} =   \xi^{i}_{k
j}  \, .
\end{equation*}
It follows that all coefficients $\, c_{I, K} \, = \, 0$.
\medskip

Let now $N \in \Nb$ be the largest height of  occurring in
(\ref{01}). By varying $\vp $ in the subgroup ${\Gb}^{(N+2)} (0)
\subset {\Gb}_0$ of all diffeomorphisms whose $(N+2)$-jet at $0$ has
the form
\begin{eqnarray*}
J_0^{N+2}  (\vp)^i (x) &=  \, x^i + \frac{1}{(N+2)!}  \sum_{j, k,
\a_1, \ldots , \a_{N+2} } \xi^i_{j k \a_1 \ldots \a_{N}} x^j x^k
x^{\a_1} \cdots x^{\a_{N}} ,  \\ \nonumber
  \xi^i_{j k \a_1 \ldots \a_{N}} & \in \Cb^{n^{N+3}}  ,  \, \,
 \xi^i_{j k \a_1 \ldots \a_{N}} = \xi^i_{k j \a_{\s(1)} \ldots \a_{\s(N)}} ,
 \, \, \fl \, \, \text{permutation} \, \, \s \, , \nonumber
\end{eqnarray*}
and using (\ref{highg}) instead of (\ref{sgijk}), one derives as
above
 the vanishing
of all coefficients $c_{I, \k}  $ with $\vert \k \vert = N$. This
lowers the height  in (\ref{01})  and thus completes the induction.
 \end{proof}
\bigskip

Let $\, \FB$ denote the ideal of  $\FA ({\Fh}_n)$ generated by the
combinations of the form
\begin{equation}  \label{idealbianchi}
 \d_{j \ell \,  k}^i \, - \,  \d_{j k \,  \ell}^i \, - \,
 \d_{j k}^s \, \d_{s \ell}^i  \ + \, \d_{j \ell}^s \, \d_{s k}^i ,
\end{equation}
which according to \eqref{bianchi} vanish when viewed in $\Hc_n$.

\begin{corollary} \label{quotient}
The algebra $\Hc_n$ is isomorphic to the quotient of the universal enveloping
algebra
 $\FA ({\Fh}_n)$ by the ideal $\FB$.
\end{corollary}

\proof Extending the notation for indices introduced above, we form
multi-indices of the third kind, $\, \displaystyle K' \, = \, \{
\k'_1 \leq \ldots \leq \k'_r \} $, which are similar to those of the
second kind \eqref{2kind} except that we drop the requirement  $k_s
\leq \, \ell_1^s$, {\it i.e.}
\begin{equation} \label{3kind}
\k'_s \, = \, \left( \begin{matrix} &i_s &\qquad \qquad\qquad \hfill
\cr
 j_s &\leq & k_s  \, , \,
\ell_1^s \leq \ldots \leq \ell_{p_s}^s \hfill \cr\end{matrix}
\right)  \, , \qquad s=1, \ldots , r \, .
\end{equation}
We then form a Poincar\'e-Birkhoff-Witt basis of  $\FA ({\Fh}_n)$
out of the monomials  $\, \d_{K'} \, Z_I \,$, ordered
lexicographically. Let $\pi : \FA ({\Fh}_n) \ra \Hc_n$ be the
tautological algebra homomorphism, which sends the generators of
${\Fh}_n$ to the same symbols in $\Hc_n$. In particular, it sends
the basis elements of  $\FA ({\Fh}_n)$ to the corresponding elements
in  $ \Hc_n$
\begin{equation*}
\pi ( \d_{K'} \, Z_I) \, = \,  \d_{K'} \, Z_I ,
\end{equation*}
but the monomial in the right hand side belongs to the basis of  $
\Hc_n$  only when the components of all $\k'_s \in K'$ are in the
increasing order.

Evidently, $\, {\FB} \sbs \Ker \pi$. To prove the converse, assume
 that
 \begin{equation*}
 u =  \sum_{I, K'} \, c_{I, K'} \,  \d_{K'}\, Z_I \in \FA ({\Fh}_n)
  \end{equation*}
satisfies
\begin{equation}  \label{ker}
\pi(u)   \, \equiv \,
 \sum_{I, K'} \, c_{I, K'} \,  \d_{K'}\, Z_I   \, =  \,0 .
 \end{equation}
In order to show that $\, u $ belongs to the ideal $\FB$, we shall
again use induction, on the height of $u$. The height $0$ case is
obvious, because the $0$-height monomials remain linearly
independent in $ \Hc_n$, so $\pi (u)=0$ implies $c_I =0$ for each
$I$, and therefore $u=0$.

Let now $N\geq 1$ be the largest height of  occurring in $u$. For
each  $K'$ of height $N$, denote by $K$ the multi-index with the
corresponding components $\k'_s \in K'$ rearranged in the increasing
order. In view of \eqref{bianchi}, one can replace each $\, \d_{K'}$
in the equation \eqref{ker} by $\d_K + $ {\em lower height}, because
the difference belongs to $\FB$. Thus, the top height part of $\, u$
becomes $\displaystyle \sum_{I, |K| = N} c_{I, K'} \,  \d_K\, Z_I $,
and so
\begin{equation*}
u \, = \, v   \,+ \,
  \sum_{I, |K| = N} c_{I, K'} \,  \d_K\, Z_I   \quad (\text{mod}\, \FB ),
 \end{equation*}
where $v$ has height at most $N-1$. Then \eqref{ker} takes the form
\begin{equation*}
\pi(v)  \,+ \,
  \sum_{I, |K| = N} c_{I, K'} \,  \d_K\, Z_I \, =  \,0 ,
 \end{equation*}
and from Proposition \ref{free} it follows that the coefficient of
each $ \,  \d_K\, Z_I $ vanishes. One concludes that
\begin{equation*}
u \, = \, v
 \quad (\text{mod}\, \FB ) .
 \end{equation*}
On the other hand, by the induction hypothesis, $\pi (v) =0$ implies
$v \in \FB$.
\endproof
\bigskip

In order to state the next result, we associate to any element $h^1
\ot \ldots \ot  h^p \in \Hc_n^{ \ot ^p} $ a
\textit{multi-differential} operator, acting on $\Ac$, by the
following formula
\begin{eqnarray} \label{T}
T(h^1 \ot \ldots \ot  h^p) \, (a^1\ot \ldots \ot  a^p)
\, = \, h^1(a^1) \cdots h^p (a^p) \, , \\
\text{where}  \quad  h^1, \ldots , h^p \in \Hc_n \, \quad \text{and}
\quad a^1, \ldots , a^p \in\Ac \, ;\nonumber
\end{eqnarray}
the linear extension of this assignment will be denoted by the same
letter.
\smallskip

\begin{proposition} \label{multid}
For each $p \in \Nb$, the linear transformation $ T\, : \, \Hc_n^{
\ot ^p}  \longra \, \Lc (\Ac^{ \ot^ p} , \Ac) $
  is injective.
\end{proposition}
\smallskip

\begin{proof}  For $p=1$, $T$ gives the standard action of
$\Hc_n$ on $\Ac$, which was just shown to be faithful. To prove that
$ \Ker T \, = \, 0 \, $ for an arbitrary $p \in \Nb$, assume that
\begin{equation*}
H \, =\, \sum_{\rho} \, h^{1}_{\rho} \ot \cdots \ot h^{p}_{\rho} \,
\in
 \Ker T  \, .
\end{equation*}
After fixing a Poincar\'e-Birkhoff-Witt basis  as above, we may
uniquely express each $h^{j}_{\rho}$ in the form
$$
h^j_{\rho} = \, \sum_{I_{j}, K_{j}} \ C_{\rho , \, I_{j}, K_{j}}
 \,  \d_{K_{j}} \, Z_{I_{j} } \, , \quad
\text{with} \quad C_{\rho , \, I_{j}, K_{j}} \in \Cb \, .
$$
Evaluating $T(H)$ on elementary tensors of the form $ \, f_1
U_{\vp_1}^* \ot \cdots  \ot f_p U_{\vp_p}^*  \, $, one obtains
\begin{equation*}
\sum_{\rho, I , K} \, C_{\rho, \, I_{1}, K_{1}} \cdots C_{\rho, \,
I_{p}, K_{p}}
 \,  \d_{K_{1}} \left(Z_{I_{1} } (f_1) U_{\vp_1}^* \right)  \cdots \,
\d_{K_{p}} \left( Z_{I_{p} } (f_p) U_{\vp_p}^* \right) \, = \, 0 \,
.
\end{equation*}
Evaluating further at a point $\, u_1 = (x_1, {\bf y}_1) \in F
{\Rb}^n$, and denoting
$$
u_2 = \wt{\vp}_1 (u_1) \, ,  \ldots , \, u_{p} = \wt{\vp}_{p-1}
(u_{p-1}) \, ,
$$
the above identity gives
\begin{eqnarray*}
\sum_{\rho, I , K} \, C_{\rho, \, I_{1}, K_{1}} \cdots C_{\rho, \,
I_{p}, K_{p}} &\cdot&  \g_{K_{1}} ({\vp}_1) (u_1)  \cdots \g_{K_{p}}
({\vp}_p) (u_p) \\ \nonumber &\cdot&  Z_{I_{1} } (f_1) (u_1) \cdots
Z_{I_{p} } (f_p) (u_p)\, = \, 0 \, . \nonumber
\end{eqnarray*}

Let us fix points $\, u_1, \, \ldots , \, u_p \in  F {\Rb}^n$ and
then diffeomorphisms $\psi_0, \, \psi_1, \, \ldots , \, \psi_p$,
such that
$$  u_2 = \wt{\psi}_1 (u_1) \, ,  \ldots , \, u_{p} =
\wt{\psi}_{p-1} (u_{p-1}) \, .
$$
Following a line of reasoning similar to that of the preceding
proof, and iterated with respect to the points $\, u_1, \ldots , u_p
\,$, we can infer that for  each $p$-tuple of indices of the first
kind $ \, ( I_1 , \ldots , I_p) \, $ one has
\begin{equation*}
\sum_{\rho, K} \, C_{\rho, \, I_{1}, K_{1}} \cdots C_{\rho, \,
I_{p}, K_{p}} \cdot \g_{K_{1}} ({\vp}_1) (u_1)  \cdots \g_{K_{p}}
({\vp}_p) (u_p)
 \, = \, 0 \, .
\end{equation*}
Similarly, making repeated use of diffeomorphisms of the form
$$
\psi_k  \circ \vp \quad \text{with} \quad \vp \in  {\Gb}^{(N)} (u_k)
\, , \quad k=\, 1, \ldots, p \, ,
$$
for sufficiently many values of $N$, we can eventually conclude that
for any $ \, ( K_1 , \ldots , K_p) \,$
\begin{equation*} \label{nul*}
\sum_{\rho} \, C_{\rho, \, I_{1}, K_{1}} \cdots C_{\rho, \, I_{p},
K_{p}}
 \, = \, 0 \, .
\end{equation*}
This proves that $\, H \, = \, 0$.
\end{proof}
\bigskip

The crossed product algebra $ {\Ac}  = C_c^{\ify} (F{\Rb}^n )
\rtimes {\Gb}$ carries a canonical trace, uniquely determined up to
a scaling factor. It is defined as the linear functional $\tau:\Ac
\ra \Cb$,
\begin{equation} \label{tr}
\tau \, (f \, U_{\vp}^* )  \, = \, \left\{ \begin{matrix}
\displaystyle
  \int_{F{\Rb}^n} \, f  \, \varpi \, , \quad \text{if}
\quad \vp = Id \, , \cr\cr \quad 0 \, , \qquad \qquad
\text{otherwise} \, .
\end{matrix} \right.
\end{equation}
Here $\varpi$  is the volume form attached to the canonical framing
given by the flat connection
\begin{equation*} \label{vol}
\varpi \, = \, \bigwedge_{k=1}^n \t^k   \wedge \bigwedge_{(i, j)}
\om^i_j  \qquad \text{(ordered lexicographically)} \, .
\end{equation*}
 The tracial property
\begin{equation*}
\tau (a\, b) \, = \,  \tau(b \, a) \, , \qquad \fl \, a, b \in \Ac ,
 \end{equation*}
 is a consequence of the
${\Gb}$-invariance of the volume form $\varpi$. In turn, the latter
follows from the fact that
\begin{equation*}
\wt{\vp}^*(\t) = \t \quad \text{and} \quad  \wt{\vp}^* (\om) \, = \,
\om \, + \, \g \cdot \t  \, ;
 \end{equation*}
indeed,
\begin{equation*}
\wt{\vp}^*(\varpi) \, =\,  \bigwedge_{k=1}^n \t^k
 \wedge \bigwedge_{(i, j)} \left(\om^i_j  + \g^i_{j\ell}(\vp) \t^\ell \right)
\,  = \, \bigwedge_{k=1}^n \t_k   \wedge \bigwedge_{(i, j)} \om^i_j
\, .
 \end{equation*}
This trace satisfies an invariance property relative to the
\textit{modular character} of $\Hc_n $. The latter,
 $\d :\Hc_n \ra \Cb$, extends the \textit{trace} character of
 ${\Fg \Fl} (n, \Rb)$, and is
 defined on the algebra generators as follows:
 \begin{equation} \label{chard}
\d(Y_i^j) = \d_i^j , \quad  \d(X_k) = 0, \quad \d(\d_{jk}^i ) =0 ,
 \qquad i, j, k =1, \ldots , n \, .
 \end{equation}
 Clearly, this definition is compatible with the relations \eqref{bianchi}
 and therefore extends to a character of the algebra $\Hc_n$.
 \medskip

\begin{proposition} \label{ibp}
For any $a, b \in \Ac$ and $\, h \in \Hc_n$  one has
\begin{equation} \label{it}
\tau (h(a)) \, = \,  \d(h)\, \tau(a) \, .
\end{equation}
\end{proposition}
\medskip

\begin{proof}
It suffices to verify the stated identity on the algebra generators
of $\Hc_n$. Evidently, both sides vanish if $h = \d^i_{jk}$. On the
other hand, its restriction to the Lie algebra ${\Fg}= \Rb^n \rtimes
{\Fg \Fl} (n, \Rb)$ is just the restatement, at the level of the Lie
algebra, of  the invariance property of the left Haar measure on $G
= \Rb^n \rtimes \GL (n, \Rb)$ with respect to right translations.
\end{proof}
\medskip

\begin{proposition}
There exists a unique anti-automorphism $ \wt{S} : \Hc_n \ra \Hc_n$
such that
\begin{equation} \label{sit}
\tau (h(a) \, b) \, = \, \tau \, (a \, \wt{S} (h) (b)) \, ,
\end{equation}
for any $h \in \Hc_n$ and $\, a, b  \in \Ac$. Moreover, $\wt{S}$ is
involutive:
 \begin{equation} \label{ainv}
 \wt{S}^2 \, = \,\Id .
 \end{equation}
\end{proposition}
\medskip

\begin{proof}
Using  the `Leibnitz rule'  (\ref{Yrule}) for vertical vector
fields, and the invariance property (\ref{it}) applied to the
product $a b$, $a, b \in \Ac$, one obtains
\begin{equation}  \label{Ytap}
\displaystyle \tau \, (Y_i^j (a) \, b ) \, = \, - \, \tau \, (a \,
Y_i^j (b)) \,
 + \, \d_{i}^{j} \, \tau \, (a \, b ) \, , \qquad \fl \, a, b \in \Ac \, . \hfill
\end{equation}
On the other hand, for the basic horizontal vector fields,
(\ref{Xrule}) and  (\ref{it}) give
\begin{equation}  \label{Xtap}
\begin{matrix}
 \displaystyle \tau \, (X_k (a) \, b ) \,  &=& \, - \,
\tau \, (a \, X_k (b) ) \, - \, \tau \, ( \d_{jk}^i (a) \, Y_i^j (b)
) \, \hfill \cr \displaystyle &=& \, - \, \tau \, (a \, X_k (b) ) \,
+ \, \tau \, (a \, \d_{jk}^i (Y_i^j (b) ) \, ; \hfill  \cr
\end{matrix}
\end{equation}
the second equality uses the $1$-cocycle nature of $\g^i_{jk}$. The
same property implies
\begin{equation}   \label{dtap}
\displaystyle \tau \, (\d^i_{jk} (a) \, b) \, = \, - \, \tau \, (a
\, \d^i_{jk}  (b)) \, , \qquad \fl \, a, b \in \Ac
   \, . \hfill
\end{equation}
Thus, the generators of $\Hc_n$ satisfy \textit{integration by
parts} identities of the form \eqref{sit}, with
\begin{eqnarray} \label{sitgen}
\wt{S} (Y_i^j ) &=&\,  - Y_i^j  \, + \, \d_{i}^{j} \\
\wt{S} (X_k ) &=&\,  - X_k  \, + \, \d^i_{jk} \, Y_i^j \\
\wt{S} (\d^i_{jk} ) &=&\, - \d^i_{jk}
\end{eqnarray}
Since the pairing $(a, b) \mapsto \tau(a\, b)$ is non-degenerate,
the above operators are uniquely determined.

 Being obviously multiplicative, the `integration by parts'
rule extends from generators to all elements $h \in \Hc_n$, and
uniquely defines a map
 ${\wt S} : \Hc_n \ra \Hc_n$
satisfying \eqref{sit}. In turn, this very identity implies that
${\wt S}$ is  a homomorphism from $\Hc_n$ to $\Hc_n^{\rm op}$, as
well as the fact that ${\wt S}$ is involutive.
 \end{proof}
\bigskip

Relying on the above results, we are now in a position to equip  $\,
\Hc_n \,$ with a canonical Hopf structure.
 \medskip

\begin{theorem} \label{hopfHn}
There exists a unique Hopf algebra structure on $ \Hc_n $ with
respect to which $ \Ac $ is a left $\, \Hc_n $-module algebra.
\end{theorem}
\medskip

\begin{proof}
The formulae (\ref{Yrule}), (\ref{Xrule}) and (\ref{drule}) extend
by multiplicatively to a general `Leibnitz rule' satisfied by any
element $h \in \Hc_n$, of the form
\begin{equation} \label{HA}
 h (a  b) \ = \ \sum_{(h)} \, h_{(1)} (a) \, h_{(2)} (b) \, ,
\qquad  h_{(1)}, h_{(2)} \in \Hc_n \, , \quad a, b \in \Ac
\end{equation}
By Proposition \ref{multid}, this property uniquely determines the
{\it coproduct} map $\D : \Hc \ra \Hc \ot \Hc$,
\begin{equation} \label{D}
\D ( h) \ = \ \sum_{(h)} \, h_{(1)} \ot h_{(2)} \, ,
\end{equation}
that satisfies
\begin{equation} \label{SA}
 T(\D h) (a\ot b)\, = \,  h (a  b)   \, .
 \end{equation}
Furthermore, the coassociativity of $\D$ becomes a consequence of
the associativity of $\Ac$, because after applying $T$ it amounts to
the identity
\begin{equation*}
h ( (a b) c) \, = \,  h ( a (b c)) \, , \quad \fl \, h \in \Hc_n \,
, \,  a, b \in \Ac \, .
 \end{equation*}
 Similarly, the property that $\D$ is an algebra homomorphism follows from
 the fact that $ \Ac $ is a left $\, \Hc_n $-module. By the very definition of
 the coproduct, $ \Ac $ is actually a left $\, \Hc_n $-{\it module algebra}.
 \smallskip

 The {\it counit} is defined by
\begin{equation} \label{coun}
 \ve( h)  \, = \,  h (1)   \, ;
 \end{equation}
when transported via $T$, its required properties amount to the
obvious identities
\begin{equation*}
h (a \, 1) \, = \,  h ( 1 \, a) \, = \,h(a) \, , \quad   \fl \, h
\in \Hc_n \, , \,  a \in \Ac \, .
 \end{equation*}
 \smallskip

It remains to show the existence of {\it antipode}. We first check
that the anti-automorphism $\wt{S}$ is a {\it twisted antipode},
i.e. satisfies for any $h \in \Hc$,
\begin{eqnarray} \label{twan1}
( \Id \ast \wt{S}) (h) := \, \sum_{(h)} \, h_{(1)}\, \wt{S}
(h_{(2)})\, &=&\,  \d (h) \,1 , \\ \label{twan2}
 ( \wt{S} \ast \Id) (h) :=\, \sum_{(h)} \, \wt{S} (h_{(1)}) \, h_{(2)}\, &=&\,  \d (h) \,1 .
\end{eqnarray}
Indeed, with $a, b \in \Ac$ arbitrary, one has
\begin{equation*}
\tau (a \, \d (h) b))= \tau (h (ab))  =  \sum_{(h)}  \tau (h_{(1)}
(a) \, h_{(2)} (b))
 =  \sum_{(h)} \tau (a \,  ( \wt{S} (h_{(1)}) h_{(2)} )(b))  ,
\end{equation*}
which proves \eqref{twan1}. Similarly, but also using the tracial
property,
\begin{equation*}
\tau (a \, \d (h) b))= \tau (h (ba))  =  \sum_{(h)}  \tau (h_{(2)}
(a) \, h_{(1)} (b))
 =  \sum_{(h)} \, \tau (a \,  ( \wt{S} (h_{(2)}) h_{(1)} )(b))  ,
\end{equation*}
or equivalently
\begin{equation*}
 \sum_{(h)} \, \wt{S} (h_{(2)}) h_{(1)} \,=\, \d (h) \,1 ;
\end{equation*}
applying $\wt{S}$ to both sides yields \eqref{twan2}.

Now let $\check{\d} \in \Hc_n^*$ denote the convolution inverse of
the character $\d \in \Hc_n^*$, which on generators is given by
 \begin{equation} \label{chard'}
\check{\d}(Y_i^j) = - \d_i^j , \quad \check{\d}(X_k) = 0, \quad
\check{\d}(\d_{jk}^i ) =0 ,
 \qquad i, j, k =1, \ldots , n \, .
 \end{equation}
 Then  $S := \check{\d} \ast \wt{S}$ is an algebra anti-homorphism
 which satisfies the antipode requirement
 \begin{equation*}
 \sum_{(h)} S (h_{(1)}) h_{(2)} \,=\, \ve (h) \,1  \,=\,
  \sum_{(h)} h_{(1)} S(h_{(2)} )
\end{equation*}
on the generators, and hence for any $h \in \Hc_n$.
 \end{proof}
\bigskip

\subsection{The general case of a flat primitive pseudogroup}  \label{subsecHG}

Let $\Pi$ be a flat primitive Lie pseudogroup of local $C^\infty$-diffeomorphisms
of $\Rb^m$ . Denote by $F_{\Pi}\Rb^m$ the
 sub-bundle of  $F\Rb^m$ consisting of the
 $\Pi$-frames on $\Rb^m$. It consists of
 the 1-jets at $0 \in {\Rb}^n$ of the germs of local diffeomorphisms
$\phi \in \Pi$. Since  $\Pi$ contains the translations,
$F_{\Pi}\Rb^m$ can be identified, by the restriction of the map
\eqref{frameid},
 to ${\Rb}^m \times G_0  (\Pi)$,
 where $G_0 (\Pi) \sbs \GL(m, \Rb)$ is the linear isotropy group,
 formed of the Jacobians at $0$ of the local diffeomorphisms
$\phi \in \Pi$ preserving the origin.

 The flat connection on $F\Rb^m$ restricts to a connection form on
 $F_{\Pi}\Rb^m$ with values in the Lie algebra
 ${\Fg}_0  (\Pi)$ of $G_0  (\Pi)$,
\begin{equation*}
 \om_{\Pi}\, :=\,   {\bf y}^{-1}  \, d{\bf y} \in {\Fg}_0 (\Pi) ,  \qquad {\bf y} \in G_0 (\Pi) .
  \end{equation*}
 The basic horizontal
 vector fields on $F_{\Pi}\Rb^m$ are restrictions of those on $F\Rb^{m}$,
\begin{equation} \label{XG}
 X_k \, = \,  y_k^{\mu} \, \frac{\part} {\part x^{\mu}}\, , \quad   k =1, \ldots , 2n \, ,
 \qquad {\bf y} = (y_i^j) \in G_0(\Pi),
 \end{equation}
and the fundamental vertical vector fields are
\begin{equation}  \label{YG}
 Y_i^j  \,=\,  y_i^{\mu} \,  \frac{\part}{\part y_j^{\mu}}  , \quad  i,  j =1, \ldots , 2n ,
  \qquad {\bf y} = (y_i^j)  \in G_0 (\Pi) ;
  \end{equation}
when assembled into a matrix-valued vector field, ${\bf Y} =
\big(Y_i^j \big)$ takes values in the Lie subalgebra ${\Fg}_0  (\Pi)
\sbs {\Fg \Fl}(m, \Rb)$.
\smallskip

By virtue of Proposition \ref{leftinv}, \eqref{XG}, \eqref{YG} are
also left-invariant vector fields, that give a framing of the group
of affine $ \Pi$-motions $G (\Pi) := {\Rb}^{m} \rtimes G_0(\Pi)$.
\smallskip

The group ${\Gb} (\Pi) := \Diff (\Rb^m) \cap \Pi$ of  global
$\Pi$-diffeomorphisms acts on  $F_{\Pi}\Rb^m$  by prolongation, and
the corresponding crossed product algebra
$$
   {\Ac (\Pi)}\, : = \, C_c^{\ify} ( F_{\Pi}\Rb^m)  \rtimes {\Gb} (\Pi)
$$
is a subalgebra of $\Ac$. After promoting the above vector fields
$X_k$  and $Y_i^j$ to linear transformations in $ \Lc \,
({\Ac}(\Pi)) $, one automatically obtains, as in \S \ref{subsecHn},
the affiliated multiplication operators
 $\, \d_{jk}^i   \in \Lc \, ({\Ac}(\Pi))  $,
\begin{equation*}
\d_{jk}^i (f \, U_{\vp}^*) \, =\, \g_{jk}^i  (\vp) \, f \, U_{\vp}^*
\, , \quad \vp \in {\Gb} (\Pi) ,
\end{equation*}
and then their higher `derivatives' $\, \d_{jk \,  \ell_1 \ldots
\ell_r}^i  \in \Lc \, ({\Ac (\Pi)})$,
\begin{equation*}
\d_{jk \,  \ell_1 \ldots \ell_r}^i (f \, U_{\vp}^*) \, =\, \g_{jk \,
\ell_1 \ldots \ell_r}^i  (\vp) \, f \, U_{\vp}^*  \, , \quad \fl \,
\vp \in {\Gb} (\Pi) .
\end{equation*}
These are precisely the restrictions of the corresponding operators
in $\Lc ({\Ac})$, characterized by the property that the $m \times
m$-matrix defined by their `isotropy part' is $
{\Fg}_0(\Pi)$-valued:
\begin{equation} \label{isotype}
\big(\g_{j \bullet \cdots \bullet}^i  (\vp)\big)_{1\leq i, j \leq m}
 \in {\Fg}_0(\Pi) \sbs {\Fg \Fl}(m, \Rb)\, , \quad \vp \in {\Gb} (\Pi).
\end{equation}
They form a Lie subalgebra ${\Fh}(\Pi)$ of ${\Fh}_m$ and satisfy the
Bianchi-type identities of Proposition \ref{Bianchi}.

We let ${\Hc(\Pi)}$ denote the generated subalgebra of $ \Lc \,
({\Ac (\Pi)})$ generated by the above operators, while ${\FB}(\Pi)$
stands for the ideal generated by the identities \eqref{bianchi}.
\smallskip

\begin{theorem} \label{hopfGHn}
The algebra ${\Hc(\Pi)} $ is isomorphic to the quotient of the
universal enveloping algebra $\FA \big({\Fh}(\Pi)\big)$ by the ideal
${\FB}(\Pi)$, and can be equipped with a unique Hopf algebra
structure with respect to which $ {\Ac (\Pi)} $ is a left
$\,{\Hc(\Pi)} $-module algebra.
\end{theorem}

\proof The proof amounts to a mere repetition of the steps followed in the
previous subsection to establish Theorem \ref{hopfHn}. It suffices
to notice that all the arguments remain valid when the isotropy-type
of the generators is restricted to the linear isotropy Lie algebra $
{\Fg}_0(\Pi)$.
\endproof

To illustrate the construction of the Hopf algebra ${\Hc(\Pi)}$ in a
concrete fashion, we close this section with a more detailed
discussion of the two main subclasses of flat primitive Lie
pseudogroups: volume preserving and symplectic.


\subsubsection{Hopf algebra of the volume preserving pseudogroup}

The sub-bundle $F_{\rm s}\Rb^n$ of $F\Rb^n$ consists in this case of
all {\em special} (unimodular) frames on $\Rb^n$, defined by taking
the $1$-jet at $0 \in \Rb^n$ of germs of local diffeomorphisms
$\phi$ on ${\Rb}^n$ that preserve the volume form, {\it i.e.}
\begin{equation*}
\phi^\ast(dx^1\wdots dx^n)= dx^1\wdots dx^n .
\end{equation*}
By means of the identification \eqref{frameid}, $F_{\rm s}\Rb^n
\simeq {\Rb}^n \times \SL(n, {\Rb})$.

  The flat
connection is given by the ${\Fs \Fl} (n, \Rb)$-valued $1$-form $\,
_{\rm s}\om = ( _{\rm s}\om^i_j ) \,$
\begin{eqnarray} \label{stcon} \nonumber
 _{\rm s}\om^i_j &:=&  \om^i_j \,= \, ({\bf y}^{-1})^i_{\mu} \, d{\bf y}^{\mu}_j  ,
 \qquad i\neq j =1, \ldots , n ; \\
  _{\rm s}\om_i^i &:=& \,\om_i^i -  \om_n^n , \qquad i =1, \ldots , n-1.
 \end{eqnarray}
 The basic horizontal
 vector fields on $F_{\rm s}{\Rb}^n$ are restrictions of those on $F\Rb^n$,
\begin{equation*}
 _{\rm s}X_k \, = \,  y_k^{\mu} \, \part_{\mu} \, , \quad   k =1, \ldots , n \, ,
 \qquad {\bf y} \in \SL(n, {\Rb})
 \end{equation*}
while the fundamental vertical vector fields are
\begin{eqnarray*}
 _{\rm s}Y_i^j &:=& Y_i^j  \,=\,  y_i^{\mu} \, \part_{\mu}^j , \qquad  i \neq j =1, \ldots , n ,\\
 _{\rm s}Y_i^i &:=& \, \frac{1}{2} \left(Y_i^i -  Y_n^n\right) , \qquad i =1, \ldots , n-1.
 \end{eqnarray*}
By Proposition \ref{leftinv}, these are also the left-invariant
vector fields on the group $G_{\rm s} := {\Rb}^n \rtimes \SL(n,
{\Rb})$ associated to the standard basis of ${\Fg}_{\rm s} := \Rb^n
\rtimes {\Fs \Fl} (n, \Rb)$.

The group ${\Gb}_{\rm s} := \Diff (\Rb^n , {\rm vol})$ of volume
preserving  diffeomorphisms acts on $F_{\rm s} {\Rb}^n$ as in
\eqref{frameact}, and the corresponding crossed product algebra
$$
   \Ac_{\rm s} \, : = \, C_c^{\ify} (F_{\rm s}{\Rb}^n ) \rtimes {\Gb}_{\rm s}
$$
is a subalgebra of $\Ac$.

The analogue of \eqref{Xphi} is
\begin{equation}
 U_{\vp}^* \, _{\rm s}X_k \, U_{\vp} \,  =
  \,  _{\rm s}X_{k} \,  - \, _{\rm s}\g_{jk}^i (\vp) \, _{\rm s}Y_i^j  \, ,
\end{equation}
where
\begin{eqnarray*}
_{\rm s}\g_{jk}^i &=& \g_{jk}^i ,   \qquad  i \neq j =1, \ldots , n
, \qquad
\text{and} \\
 _{\rm s}\g_{ik}^i &=& \g_{ik}^i - \g_{nk}^n , \qquad   i =1, \ldots , n-1.
\end{eqnarray*}
One thus obtains the multiplication operators $\, _{\rm s}\d_{jk \,
\ell_1 \ldots \ell_r}^i  \in \Lc \, ({\Ac_{\rm s}})$,
\begin{equation*}
_{\rm s}\d_{jk\,  \ell_1 \ldots \ell_r}^i (f \, U_{\vp}^*) \, =\,
_{\rm s}\g_{jk\,  \ell_1 \ldots \ell_r}^i  (\vp) \, f \, U_{\vp}^*
\, ,
\end{equation*}
which continue to satisfy the Bianchi-type identities of Proposition
\ref{Bianchi}.
\smallskip

Denoting by $S\Hc_n$ the subalgebra of  $ \Lc \, ({\Ac_{\rm s}})$
generated by the above operators, one equips it with the canonical
Hopf structure
 with respect to which
$ \Ac_{\rm s} $ is a left $\, S\Hc_n $-module algebra. We remark
that the Hopf algebra $\, S\Hc_n $ is  {\em unimodular},  in the
sense that its antipode is involutive.

The {\em conformal volume preserving} case is similar, except that the linear
isotropy subgroup is ${\rm CSL}(n, \Rb) = \Rb^+ \times \SL(n, \Rb)$.

\subsubsection{Hopf algebra of the symplectic pseudogroup}

Let ${\Gb}_{\rm sp} \subset {\Gb}$ be the subgroup of
diffeomorphisms of $\Rb^{2n}$ preserving the symplecting form,
\begin{equation}
\Om= dx^1\wg dx^{n+1}+\dots +dx^n\wg dx^{2n}.
\end{equation}
Denote by $F_{\rm sp}\Rb^{2n}$ the sub-bundle  of $F\Rb^{2n}$ formed
of symplectic frames on $\Rb^{2n}$, {\em i.e.} those defined by
taking the $1$-jet at $0 \in \Rb^{2n}$ of germs of local
diffeomorphisms $\phi$ on ${\Rb}^{2n}$ preserving the form $\Om$.
Via \eqref{frameid} it can be identified  to ${\Rb}^{2n} \times
\Sp(n, {\Rb})$. In turn, $\Sp(n, {\Rb})$ is identified to the
subgroup of matrices $A\in \GL(2n, {\Rb})$ satisfying
\begin{equation}\label{spg}
^tA\, J\, A\, =\, J , \quad \text{where}
 \quad J=\left(\begin{matrix} 0 & \Id_n \\ -\Id_n & 0 \end{matrix}\right) ,
\end{equation}
while its Lie algebra ${\Fs \Fp} (n, {\Rb})$ is formed of matrices
${\bf a} \in {\Fg \Fl}(2n, {\Rb})$ such that
\begin{equation}\label{splie}
^t{\bf a}\, J\, + \, J\,{\bf a} \, = \, 0,
\end{equation}
 Thus the flat
connection is given by the ${\Fs \Fp} (n, \Rb)$-valued $1$-form
\begin{equation*}
 _{\rm sp}\om\, :=\,   {\bf y}^{-1}  \, d{\bf y} \in {\Fs \Fp} (n, \Rb) ,  \qquad {\bf y} \in \Sp(n, {\Rb}) ,
  \end{equation*}
 the basic horizontal
 vector fields on $F_{\rm sp}{\Rb}^{2n}$ are restrictions of those on $F\Rb^{2n}$,
\begin{equation*}
 _{\rm sp}X_k \, = \,  y_k^{\mu} \, \part_{\mu} \, , \quad   k =1, \ldots , 2n \, ,
 \qquad {\bf y} \in \Sp(n, {\Rb}) ,
 \end{equation*}
and the fundamental vertical vector fields are given by
 the ${\Fs \Fp} (n, \Rb)$-valued vector field $ (_{\rm sp}Y_i^j ) \in {\Fs \Fp} (n, \Rb)$
\begin{equation*}
 _{\rm sp}Y_i^j  \,=\,  y_i^{\mu} \, \part_{\mu}^j , \quad  i,  j =1, \ldots , 2n ,
  \qquad {\bf y} \in \Sp(n, {\Rb}) .
  \end{equation*}
They also form a basis of left-invariant vector fields on the group
$G_{\rm sp} := {\Rb}^{2n} \times \Sp(n, {\Rb})$ associated to the
standard basis of ${\Fg}_{\rm sp} := \Rb^{2n} \rtimes {\Fs \Fp} (n,
\Rb)$.

The group ${\Gb}_{\rm sp} := \Diff (\Rb^n ,\Om)$ of symplectic
diffeomorphisms acts on  $F_{\rm sp}{\Rb}^{2n}$ by prolongation, and
the corresponding crossed product algebra
$$
   {\Ac}_{\rm sp}\, : = \, C_c^{\ify} (F_{\rm sp}{\Rb}^{2n}) \rtimes {\Gb}_{\rm sp}
$$
is a subalgebra of $\Ac$. The vector fields $_{\rm sp}X_k$  and
$_{\rm sp}Y_i^j$ extend to linear transformations in $ \Lc \,
({\Ac_{\rm sp}}) $, and their action brings in multiplication
operators
 $\, _{\rm sp}\d_{jk \,  \ell_1 \ldots \ell_r}^i   \in \Lc \, ({\Ac_{\rm sp}}) $,
\begin{equation*}
_{\rm sp}\d_{jk\,  \ell_1 \ldots \ell_r}^i (f \, U_{\vp}^*) \, =\,
 _{\rm sp}\g_{jk\,  \ell_1 \ldots \ell_r}^i  (\vp) \, f \, U_{\vp}^*  \, .
\end{equation*}
 We let $Sp\Hc_n$ denote the generated subalgebra of
$ \Lc \, ({\Ac_{\rm sp}})$ generated by the above operators. It
acquires a unique Hopf algebra structure  such that $ \Ac_{\rm sp} $
is a left $\, Sp\Hc_n $-module algebra. Like  $\, S\Hc_n $, the Hopf
algebra $\, Sp\Hc_n $ is unimodular.

The {\em conformal symplectic} case is again similar, except that the linear
isotropy subgroup is ${\rm CSp}(n, \Rb) = \Rb^+ \times \Sp(n, \Rb)$.
\subsection{Hopf algebra of the contact pseudogroup}

 We denote by ${\Gb}_{\rm cn} \subset {\Gb}$ the subgroup of orientation preserving
diffeomorphisms of $\Rb^{2n+1}$ which leave invariant the contact form
\begin{equation}\label{contform}
\a:= \,  - dx^0 + \frac{1}{2}\sum_{i=1}^n(x^idx^{n+i} - x^{n+i}dx^i).
\end{equation}
 The vector field $E_0:= \p_0 \equiv \frac{\p}{\p x_0}$
 satisfies $\i_{E_0} (\a) =1$ and $\i_{E_0} (d\a) = 0$, {\it i.e.} represents
 the {\em Reeb vector field} of the contact structure. The {\em
contact distribution} $\,\Dc := \Ker \a$ is spanned by the vector fields $\, \{E_1, \ldots , E_{2n}\}$,
 \begin{equation}  \label{Ei}
 E_i = \p_i - \frac{1}{2} x^{n+i} E_0, \qquad
E_{n+i} =\p_{n+i}+ \frac{1}{2}x^{i} E_0 , \qquad1 \leq i \leq n .
\end{equation}
The Lie brackets between these vector fields are precisely those
of the Heisenberg Lie algebra $\Fh_n$:
\begin{equation*}
[E_i, E_{j+n}] = \d^i_j E_0, \quad [E_0, E_i] = 0 , \quad [E_0, E_{j+n}] = 0 ,
\quad i,j=1,\dots, n .
\end{equation*}
This gives an identification of  $\Rb^{2n+1}$, whose standard basis we
denote by   $\, \{ {\eb}_0 , {\eb}_1 , \ldots,  {\eb}_{2n}\}$,
with the Heisenberg group
$H_n$, whose group law is given by
\begin{align} \label{Hlaw}
&\qquad  \qquad {\xb} \ast {\yb} = {\xb} + {\yb} + \b ({\xb}^\prime, {\yb}^\prime)
\, {\eb}_0 \, \equiv \, \big(x^0 + y^0 + \b (\xb^\prime, \yb^\prime), \, \xb^\prime, \, \yb^\prime\big)\\   \notag
&\text{where} \quad {\xb} = (x^0, x^1, \ldots, x^{2n}), \quad
{\yb} = (y^0, y^1, \ldots, y^{2n}) \in H_n \equiv \Rb^{2n+1} ,\\ \notag
&\qquad  \qquad {\xb}^\prime = ( x^1, \ldots, x^{2n}) , \quad
{\yb}^\prime =   (y^1, \ldots, y^{2n}) \in \Rb^{2n} \\  \label{sympb}
&\text{and} \qquad  \b ({\xb}^\prime, {\yb}^\prime) =
\frac{1}{2}\sum_{i=1}^n(x^i\, y^{n+i} - x^{n+i}\, y^i) .
\end{align}

\begin{lemma}
The vector fields $\nobreak{\{E_0, E_1, \ldots, E_{2n}\}}$ are
respectively the left-invariant vector fields on $H_n$ determined by
the basis $\{ {\eb}_0 , {\eb}_1 , \dots, {\eb}_{2n}\}$ of its Lie
algebra $\Fh_n$.
\end{lemma}

\proof Straightforward, given that in the above
realization the exponential map $\exp: \Fh_n \ra H_n$
coincides with the identity map $\Id: \Rb^{2n+1} \ra \Rb^{2n+1}$.
\endproof

We recall that a {\em contact
diffeomorphism} is a diffeomorphism   $\phi: \Rb^{2n+1} \ra {\Rb}^{2n+1}$ such that
$\, \phi_\ast (\Ker \a)=\Ker \a, \, $  or equivalently $\,
\phi^\ast(\a)=f\a$ for a nowhere vanishing function $f : \Rb^{2n+1} \ra \Rb$;
$\phi$ is orientation preserving iff  $f > 0$ and
is called a  {\em strict contact diffeomorphism} if $\, f \equiv 1$.

In particular, the group left translations $L_{\ab}:H_n \ra H_n$, $\, {\ab}  \in H_n$,
\begin{align*}
L_{\ab} ({\xb}) \, = \, {\ab} \ast  {\xb} \, = \,\big(a^0 + x^0 +
\b ({\ab}^\prime,  {\xb}^\prime) , \, {\ab}^\prime + {\xb}^\prime\big)
\qquad \fl \,  {\xb} \in H_n ,
\end{align*}
are easily seen to be strict contact diffeomorphism.
They replace the
usual translations in $\Rb^{2n+1} $. Together with the linear transformations
preserving the symplectic form $\b : \Rb^{2n} \times \Rb^{2n} \ra   \Rb^{2n} $
of \eqref{sympb},
\begin{align*}
S_A ({\xb}) \, = \, (x^0 , A {\xb}^\prime) ,
\qquad A \in \Sp (n, \Rb) , \qquad {\xb}^\prime = (x^1, \ldots , x^{2n}) ,
\end{align*}
they form a group of strict contact transformations
which is isomorphic to $H_n \rtimes \Sp (n, \Rb)$.
Indeed, one has
\begin{align*}
&(S_A \circ  L_{\ab})({\xb}) = S_A \big(x^0 + a^0 +
\b ({\ab}^\prime,  {\xb}^\prime) , \, {\xb}^\prime + {\ab}^\prime\big) \\
&= \big(x^0 + a^0 + \b ({\ab}^\prime,  {\xb}^\prime),
A ( {\xb}^\prime + {\ab}^\prime) \big) =
\big(x^0 + a^0 + \b (A {\ab}^\prime,  A {\xb}^\prime),
A( {\xb}^\prime + {\ab}^\prime) \big) \\
& \qquad = L_{S_A ({\ab})} \big(S_A ({\xb})\big) =  (L_{S_A ({\ab})} \circ S_A)({\xb}) ,
\end{align*}
therefore $\quad S_A \circ  L_{\ab} \circ S_A^{-1} = L_{S_A ({\ab})} $;  also
\begin{align*}
&(L_{\ab} \circ S_A) \circ (L_{\boldb} \circ S_B) = L_{\ab} \circ (S_A \circ L_{\boldb})\circ S_B
=  L_{\ab} \circ (L_{S_A ({\boldb})} \circ S_A)  \circ S_B \\
&\qquad = L_{{\ab} \ast S_A ({\boldb})} \circ S_{A B} .
\end{align*}
We shall enlarge this group the $1$-parameter group of {\it contact
homotheties} $\, \{ \mu_t ; \, t \in \Rb_+^\ast\} $, defined by
\begin{align*}
\mu_t ({\xb}) \, = \, (t^2 x^0 , t {\xb}^\prime) , \qquad  {\xb} \in H_n ;
\end{align*}
one has
\begin{align*}
\mu_t \circ S_A = S_A \circ \mu_t \qquad \text{and} \qquad \mu_t \circ L_{\ab} =
L_{\mu_t({\ab})} \circ  \mu_t .
\end{align*}
We denote by $G_{\rm cnt}$ the group
generated by the contact transformations $\, \{L_{\ab} , S_A , \mu_t ; \,
{\ab} \in H_n,  A \in \Sp (n, \Rb) , t \in \Rb_+^\ast\} $. It can be identified with
the semidirect product
\begin{align*}
 G_{\rm cn} \, \cong \, H_n \rtimes \CSp (n, \Rb)  ,
\end{align*}
where $\, \CSp (n, \Rb) = \Sp (n, \Rb) \times \Rb_+^\ast$ is
the {\it conformal symplectic group}.
Via the identification
\begin{align}  \label{param}
 L_{\ab} \circ S_A \circ \mu_t  \, \cong \, ({\ab}, A, t) ,
 \qquad {\ab} \in  H_n , \, A \in \Sp (n, \Rb) , \, t \in \Rb_+^\ast ,
\end{align}
 the multiplication law is
\begin{align} \label{contaff}
  (\ab, A, t) \cdot  (\boldb, B, s)  \, = \,  (\ab \ast  S_A (\mu_t (\boldb)), AB,  ts) .
\end{align}
\smallskip

Besides the usual tangent bundle $T \Rb^{2n+1}$, it will be convenient to introduce
a version of it that arises naturally when the contact
structure is treated as special case of a Heisenberg manifold (\cf \eg~\cite{ponge}).
Denoting by $\, \Rc : = T \Rb^{2n+1} \slash \Dc$
the line bundle determined by the class of the Reeb vector field
modulo the contact distribution, the {\it H-tangent bundle} is the direct sum
$\, T^{\rm H} \Rb^{2n+1} :=\Rc \oplus \Dc$. If $\, \phi \in \Gb_{\rm cn}$
is a contact diffeomorphism, then its tangent map at
  $\xb \in \Rb^{2n+1}$,
 $\,  \phi_{\ast \xb} : T_\xb \Rb^{2n+1} \ra T_{\phi(\xb)}\Rb^{2n+1}$,
 leaves the contact distribution $\Dc$ invariant, and hence induces a
corresponding {\it H-tangent map}, which will be denoted
$\,  \phi^{\rm H} _{\ast \xb} : T_\xb^{\rm H} \Rb^{2n+1} \ra T_{\phi(\xb)}^{\rm H} \Rb^{2n+1}$.

\begin{lemma} \label{HKac}
 Given any contact diffeomorphism $\, \phi \in \Gb_{\rm cn}$, there is a unique
 $\, \vp \in G_{\rm cn}$, such that $\, \psi = \vp^{-1} \circ \phi$ has the properties
\begin{align} \label{Hunip}
 \psi (0) = 0 \qquad \text{and} \qquad
  \psi^{\rm H} _{\ast 0} = \Id : T_0^{\rm H} \Rb^{2n+1} \ra T_0^{\rm H} \Rb^{2n+1}.
  \end{align}
\end{lemma}

\proof  Consider the H-tangent map $\nobreak{ \phi^{\rm H} _{\ast 0}
: T_0^{\rm H} \Rb^{2n+1} \ra T_{\ab}^{\rm H} \Rb^{2n+1}}$,  where
$\ab = \phi (0)$.   Then relative to the moving frame $\{E_1,
\ldots, E_{2n}\}$ for $\Dc$, $\phi^{\rm H} _{\ast 0} \mid \Dc :
\Dc_0  \ra \Dc_{\ab} $ is given by a conformal symplectic matrix
$A_\phi (0)\in \CSp (n, \Rb) $ with conformal factor $t_\phi(0)^2 >
0$. Furthermore, one can easily see that
  \begin{align*}
 &\phi_{\ast 0} (\eb_0) \, = t_\phi(0)^2 E_0 \mid_{\ab} + \wb_\phi ,
\qquad \text{with } \quad \wb_\phi \in  \Dc_{\ab} ,
   \end{align*}
and so the restriction  $\phi^{\rm H} _{\ast 0} \mid \Rc : \Rc_0
\ra \Rc_{\ab} $ can be identified with
  the scalar $\, t_\phi (0)^2$ relative to the moving frame  $E_0$ of $\Rc$.

 Since $\,   A = t_\phi^{-1} (0)A_\phi (0) \in \Sp(n, \Rb) $, we can form
 \begin{align} \label{contaffdiff}
 \vp : = L_\ab \circ S_A \circ \mu_t  \in G_{\rm cn}, \qquad
 \ab =  \phi (0) , \qquad t = t_\phi (0) ,
   \end{align}
and  we claim that  this  $\, \vp \in G_{\rm cn}$
 fulfills the required property. Indeed, $L^{-1}_\ab \circ \vp$
 has the expression
  \begin{align*}
( L^{-1}_\ab \circ \vp) ({\xb}) = \big( t_\phi(0)^2 x^0 , \, t_\phi(0) A {\xb}^\prime \big)
   \end{align*}
 and therefore the matrix of $\, \vp_{\ast 0} : T_0 \Rb^{2n+1} \ra
T_{\ab} \Rb^{2n+1} $ is the same as that of
$ \phi^{\rm H} _{\ast 0} : T_0^{\rm H} \Rb^{2n+1} \ra T_{\ab}^{\rm H} \Rb^{2n+1}$, \ie
\begin{equation} \label{Htang}
(L^{-1}_\ab \circ \vp)_{\ast 0}  \, =\,  \left(\begin{array}{c|c}   t_\phi(0)^2&   0\,  \ldots \,0\\
\hline\\
0& \\
\vdots& t_\phi(0) A\\
0&
\end{array}\right)
 \, =\,  \left(\begin{array}{c|c}   t_\phi(0)^2&   0\,  \ldots \,0\\
\hline\\
0& \\
\vdots& A_\phi(0)\\
0&
\end{array}\right).
\end{equation}
Thus, $\, \psi = \vp^{-1} \circ \phi$ satisfies \eqref{Hunip}.

To prove uniqueness, observe that if
$\vp =  L_{\ab} \circ S_A \circ \mu_t \in G_{\rm cn} $ fulfills \eqref{Hunip}, then
$\, \ab =0$ and, since $\,  \vp_{\ast 0} : T_0 \Rb^{2n+1} \ra
T_0 \Rb^{2n+1} $ is of the form \eqref{Htang}, $A = \Id$ and $t=1$.
\endproof

The subgroup of $\, \Gb_{\rm cn}$ defined by the two conditions in \eqref{Hunip}
will be denoted $\, N_{\rm cn} $.  Lemma \ref{HKac}  gives
a {\it Kac decomposition} for the group of orientation preserving contact
diffeomorphisms:
\begin{align} \label{HKdec}
\Gb_{\rm cn} \, = \, G_{\rm cn} \cdot N_{\rm cn} .
\end{align}
Note though that the group  $\, N_{\rm cn}$ is no longer pro-unipotent.
However, it has a pro-unipotent normal subgroup of finite codimension,
namely
\begin{align*}
U_{\rm cn} := \{\psi \in N_{\rm cn} \, \mid \,  \psi (0) = 0 , \quad
  \psi_{\ast 0} = \Id : T_0 \Rb^{2n+1} \ra T_0 \Rb^{2n+1} \} ;
\end{align*}
this gives rise to a group extension
 \begin{align} \label{Hext}
\Id \ra  U_{\rm cn}  \ra  N_{\rm cn}  \ra \Rb^{2n} \ra 0 ,
\end{align}
with the last arrow given by the tangent map at $0$.
\medskip

The decomposition \eqref{HKdec} gives rise to a pair of actions
of $\, \Gb_{\rm cn}$: a left action $\rt$ on $\, G_{\rm cn} \cong
\Gb_{\rm cn} \slash  N_{\rm cn}$, and a right action  $\lt$
on $\, N_{\rm cn} \cong
\Gb_{\rm cn} \bash  G_{\rm cn}$.
To understand the left action,
 let $\, \phi \in \Gb_{\rm cn}$ and
$\, \vp = L_{\xb} \circ S_A \circ \mu_t \in G_{\rm cn}$. By Lemma  \ref{HKac},
 \begin{align*}
 \phi \circ \vp \, = \,  (\phi \rt \vp) \circ  \psi , \qquad  \text{with} \quad
 \phi \rt \vp \in G_{\rm cn}
  \quad \text{and} \quad \psi \in N_{\rm cn} .
 \end{align*}
Write $ \,  \phi \rt \vp = L_{\boldb} \circ S_B \circ \mu_s$.
If in  \eqref{contaffdiff} one replaces
 $\phi$ by  $\phi\circ \vp$ then,
\begin{align*}
A_{\phi\circ \vp} (0) = A_\phi (\xb) A_\vp (0) = t A_\phi (\xb) A ,
\qquad t_{\phi\circ \vp} (0) = t_\phi (\xb) t ,
\end{align*}
hence
\begin{align*}
\phi \rt \vp \, = \,  L_{\phi (\xb)} \circ S_{t_\phi(\xb)^{-1} A_\phi(\xb) A} \circ \mu_{t_\phi t}
  \end{align*}

  Using the parametrization \eqref{param}, the explicit description of the action
  can be recorded as follows.

   \begin{lemma}
   Let $\, \phi \in \Gb_{\rm cn}$ and $(\ab, A, t) \in G_{\rm cn}$. Then
   \begin{equation}  \label{rtact}
  \tilde\phi  (\xb, A, t) :=  \phi \rt  (\xb, A, t) \, = \,
  (\phi (\xb) , \,  t^{-1}_\phi(\xb) A_\phi(\xb)  A, \, t_\phi(\xb) t ) .
   \end{equation}
   \end{lemma}
 \medskip

 With the above notational convention, let
 $R_{(\boldb, B, s)}$ denote the right translation by an element
 $(\boldb, B, s) \in \Sp (n, \Rb)$, and let $\phi \in \Gb_{\rm cn}$.
We want to understand the commutation relationship between these two
transformations.
 By \eqref{contaff} and \eqref{rtact}, one has
 \begin{eqnarray} \nonumber
&(\tilde\phi \circ R_{(\boldb, B, s)}) (\xb, A, t)=
 \tilde\phi \big( (\xb, A, t) \cdot (\boldb, B, s)\big)  =
 \tilde\phi (\xb \ast S_A (\mu_t(\boldb)), AB, ts) \\ \nonumber
&=\big(\phi (\xb \ast S_A (\mu_t(\boldb))), t^{-1}_\phi(\xb)  A_\phi(\xb)  AB, t_\phi(\xb)  t s\big) .
 \end{eqnarray}
On the other hand,
 \begin{eqnarray} \nonumber
&(R_{(\boldb, B, s)} \circ \tilde\phi)(\xb, A, t)=
(\phi (\xb), \, t^{-1}_\phi(\xb)   A_\phi(\xb)  A, \, t_\phi(\xb)  t ) \cdot (\boldb, B, s) \\  \nonumber
&= (\phi (\xb) \ast  t^{-1}_\phi(\xb)  S_{A_\phi(\xb)  A} (\mu_{t_\phi(\xb)  t }(\boldb)),
\, t^{-1}_\phi(\xb)   A_\phi (\xb) A B , \, t_\phi(\xb)  t s) .
 \end{eqnarray}
Although these two answers are in general different, when $\, \boldb = 0$
 they do coincide, and we record this fact in the following statement.

 \begin{lemma} \label{commact}
 The left action of  $\, \Gb_{\rm cn}$ on  $\, G_{\rm cn}$ commutes
 with the right translations by the elements of the subgroup
 $\, \CSp (n, \Rb)$.
 \end{lemma}
 \medskip

 As in the flat case, we proceed to associate to the pseudogroup
 $\Gb_{\rm cn}$ of orientation preserving
diffeomorphisms of $\Rb^{2n+1}$ a Hopf algebra
 $\Hc (\Pi_{\rm cn})$, realized via
 its Hopf action on the crossed product algebra
 $ \Ac (\Pi_{\rm cn}) = C^\infty (G_{\rm cn}) \rtimes \Gb_{\rm cn}$.
 This type of construction actually applies whenever
 one has a Kac decomposition of the form \eqref{HKac}.

One starts with a fixed basis $\{X_i\}_{1\leq i \leq m}$
for the Lie algebra
$\Fg_{\rm cn}$ of $G_{\rm cn}$. Each $X \in \Fg_{\rm cn}$ gives rise to
a left-invariant vector field $X$ on  $G_{\rm cn}$,
which is then extended to a linear operator on $\Ac (\Pi_{\rm cn})$,
\begin{equation*}
X (f \, U_{\phi}^* )\, = \, X (f) \, U_{\phi}^* \, , \qquad \phi \in  \Gb_{\rm cn} \, .
\end{equation*}
One has
\begin{equation} \label{UXU}
 U_{\phi}^* \, X_i \, U_{\phi} \,  =  \,  \sum_{j=1}^m  \G_i^j (\phi)\, X_j\, , \qquad i=1, \ldots , m ,
\end{equation}
with $ \G_i^j (\phi)  \in C^\infty (G_{\rm cn})$. The  matrix
of functions $\, \Gbo (\phi) = \big(\G_i^j (\phi) \big)_{1 \leq i,j \leq m} $ automatically
satisfies the cocycle identity
\begin{equation}  \label{cocid}
  \Gbo (\phi \circ \psi)\, = \, ( \Gbo (\phi)  \circ \psi) \cdot  \Gbo (\psi)  \, ,
   \qquad \phi, \psi \in  \Gb_{\rm cn} \, .
\end{equation}
We next denote by $\, \D_i^j (\phi) $ the following multiplication operator
on $\Ac (\Pi_{\rm cn})$:
\begin{equation*}
  \D_i^j (f \, U^\ast_\phi) =   (\Gbo(\phi)^{-1})_i^j  \, f \, U^\ast_\phi
  \, , \qquad i, j=1, \ldots , m \, .
\end{equation*}

With this notation, we define $\Hc (\Pi_{\rm cn})$ as the
 subalgebra of  linear operators on $\Ac (\Pi_{\rm cn})$ generated by
 the operators $X_k$'s and $\D_i^j$'s, $ i, j, k = 1, \ldots, m$. In particular, $\Hc_\Pi$
 contains all iterated commutators
 \begin{equation*}
 \D^j_{i, k_1 \ldots k_r} := [X_{k_r}, \ldots , [X_{k_1}, \D_i^j] \ldots] ,
 \end{equation*}
 \ie the multiplication
 operators by the functions on $G$,
 \begin{equation*}
 \G^j_{i, k_1 \ldots k_r} (\phi) := X_{k_r} \ldots X_{k_1} ( \G_i^j (\phi)) , \qquad \phi \in \Gb .
 \end{equation*}

 \begin{lemma}
 For any $a, b \in \Ac (\Pi_{\rm cn})$,  one has
\begin{align} \label{leibX}
&  X_k (a b) \, = \, X_k (a) \, b +  \sum_j \D^j_k (a) \, X_j (b) , \\ \label{leibG}
&  \D_i^j (ab)  \, = \, \sum_k  \D_i^k (a) \, \D_k^j (b) .
 \end{align}
  \end{lemma}

  \proof
With $\, a=f_1 \, U_{\phi_1}^*$, $\, b = f_2 \,
U_{\phi_2}^*$, and assembling the $X_k$'s into a column
vector $\Xb$ and the $\D_i^j$'s into a matrix $ \Dbo$, one has
\begin{eqnarray*}
\Xb (a \cdot b) &=&  \Xb (f_1 \, U_{\phi_1}^* \, f_2 \,
U_{\phi_2}^*) \,
 = \, \Xb (f_1 \, U_{\phi_1}^* \, f_2 \, U_{\phi_1}) \, U_{\phi_2 \phi_1}^* =\\
&=& \Xb (f_1) \, U_{\phi_1}^* \, f_2 \, U_{\phi_2}^* \, + \, f_1 \,
\Xb (U_{\phi_1}^* \, f_2 \, U_{\phi_1}) \,
U_{\phi_1}^* U_{\phi_2}^* = \\
 &=& \Xb (a) \, b \, +\, f_1 \, U_{\phi_1}^* (U_{\phi_1} \, \Xb \, U_{\phi_1}^*) (f_2) \,U_{\phi_2}^\ast
 = \qquad \quad \text{[using \eqref{UXU}]} \\
 &=&  \Xb (a) \, b \, +\,  f_1 U_{\phi_1}^* \Gbo (\phi_1^{-1}) \Xb (f_2) \,U_{\phi_2}^\ast = \\
  &=&  \Xb (a) \, b \, +\,  f_1
 (\Gbo (\phi_1^{-1})\circ \phi_1)\, U_{\phi_1}^*\,  \Xb (f_2) U_{\phi_2}^\ast =
 \quad \text{[using \eqref{cocid}]} \\
  &=&  \Xb (a) \, b \, +\,  f_1 \, \Gbo (\phi_1)^{-1}\, U_{\phi_1}^* \, \Xb (b)\, = \,
  \Xb (a) \, b \, +\,  \Dbo (a) \, \Xb (b) ,
\end{eqnarray*}
 which proves \eqref{leibX}.

 The identity \eqref{leibG} is merely a reformulation of the cocycle identity \eqref{cocid}.
 \endproof

As a consequence, by multiplicativity
every $h \in \Hc (\Pi_{\rm cn})$ satisfies a \textit{Leibniz rule}
 of the form
 \begin{equation} \label{leibh}
  h (a b) \, = \,  \sum h_{(1)} (a) \,
h_{(2)} b)  \, , \qquad \fl \, a , b \in \Ac (\Pi) .
\end{equation}

\medskip

\begin{proposition} \label{BP}
The  operators $\, \D_{\bullet \cdots \bullet}^{\bullet} \,$
satisfy the (Bianchi) identities
\begin{equation}  \label{BI}
 \D_{i, j}^k \, - \,  \D_{j, i }^k \, = \,
\sum_{r, s} c^k_{rs}\, \D_{i}^r \, \D_{j}^s  \ - \, \sum_\ell c_{i j}^\ell \, \D_{\ell}^k ,
\end{equation}
where $c^i_{j k}$ are the structure constants of $\Fg_{\rm cn}$,
\begin{equation}  \label{SC}
 [X_j , X_k] \, = \, \sum_i c^i_{j k} X_i .
\end{equation}
 \end{proposition}

\proof  Applying  \eqref{leibX} one has, for any $\, a, b \in \Ac (\Pi_{\rm cn})$,
\begin{eqnarray*}
 X_i X_j(a \, b)&=& X_i \big(X_j (a) \, b \, +\, \sum_s \D_j^s (a) X_s (b)\big) =\\
&=&X_i (X_j (a)) \, b \, + \, \sum_r \D^r_i (X_j(a))\,  X_r( b) \, +
\, \sum_{s} X_i(\D_j^s (a)) \,  X_s (b) \\
&+&\sum_{r, s} (\D_i^r (\D_j^s (a)) X_r ( X_s (b)) ,
\end{eqnarray*}
and thus the commutators can be expressed as follows:
\begin{eqnarray*}
[X_i , X_j](a \, b)&=& [X_i , X_j](a)\, b \, +\sum_r \D^r_i (X_j(a))\,  X_r( b) \, -
\sum_s \D^s_j (X_i(a))\,  X_s( b) \\
&+& \, \sum_{s} X_i(\D_j^s (a)) \,  X_s (b)\, - \,  \sum_{r} X_j(\D_i^r (a)) \,  X_r (b) +\\
&+& \, \sum_{r, s} (\D_i^r (\D_j^s (a)) X_r ( X_s (b)) \, - \, \sum_{r, s} (\D_j^r (\D_i^s (a)) X_r ( X_s (b)) \\
&=& [X_i , X_j](a)\, b \, -\, \sum_r \big(\D^r_{i, j} (a) -  \D_{j,i}^r (a)\big)\,  X_r( b) \, +\\
&+& \, \sum_{r, s} (\D_i^r \D_j^s) (a) \sum_k c^k_{rs} X_k (b)  .
 \end{eqnarray*}
 On the other hand, by \eqref{SC}, the left hand side equals
\begin{eqnarray*}
\sum_\ell c^\ell_{i j} \, X_\ell (a \, b) &=&
\sum_\ell c^\ell_{i j}\, X_\ell (a) \, b \, + \, \sum_\ell c^\ell_{i j}\, \D^k_\ell (a) \, X_k (b)\\
&=& [X_i , X_j](a)\, b \, + \, \sum_\ell c^\ell_{i j} \, \D^k_\ell (a) \, X_k (b) .
  \end{eqnarray*}
Equating the two expressions one obtains after cancelation
\begin{eqnarray*}
&\sum_\ell c^\ell_{i j} \, \D^k_\ell (a) \, X_k (b) \, =\,
-\, \sum_k \big(\D^k_{i, j} (a) -  \D_{j,i}^k (a)\big)\,  X_k( b) \, +\\
&+\, \sum_k \sum_{r, s} c^k_{rs} \, (\D_i^r \D_j^s) (a)\, X_k (b)  .
  \end{eqnarray*}
Since $\, a, b \in \Ac (\Pi)$ are arbitrary and the $X_k$'s are linearly independent, this
gives the claimed identity.
\endproof
 \medskip

 Let $\FH_{\rm cn}$ be the Lie algebra generated by
 the operators $X_k$ and $ \D^j_{i, k_1 \ldots k_r} $,
 $ i, j, k_1 \ldots k_r = 1, \ldots, m$, $r \in \Nb$. Following the same line
of arguments as in the proof of Corollary \ref{quotient}, one can
establish its exact analog.

\begin{proposition} \label{qcn}
The algebra $\Hc(\Pi_{\rm cn})$ is isomorphic to the quotient of the universal
enveloping algebra
 $\Uc(\FH_{\rm cn})$ by the ideal $\Bc_{\rm cn}$ generated by the Bianchi
 identities \eqref{BI}.
\end{proposition}

Actually, one can be quite a bit more specific about the above
cocycles as well as about the corresponding
Bianchi identities, if one uses an appropriate basis of the Lie algebra $\Fg_{\rm cn}$.
Recalling that $\Fg_{\rm cn}$ is a semidirect product
of the Heisenberg Lie algebra $\Fh_n$ by the Lie algebra $\Fg_{\rm csp}$
of the conformal symplectic group $\CSp(n, \Rb)$, one can choose the
basis $\{X_i\}_{1\leq i \leq m}$ such that the first $2n+1$ vectors are
the basis $\{E_i\}_{0\leq i \leq 2n}$ of  $\Fh_n$, while the rest form
 the canonical basis $\{Y_i^j , Z\}$ of $\Fg_{\rm csp}$, with $Z$ central.
 By Lemma \ref{commact}, for any $\phi \in \Gb_{\rm cn}$,
 \begin{equation}
 U^\ast_\phi \, Y \, U_\phi \, = \, Y , \qquad  Y \in \Fg_{\rm csp} .
\end{equation}
Thus,  the elements of $\Fg_{\rm csp}$ act as derivations on $\Ac
(\Pi)$, and therefore give rise to `tensorial identities'. The only
genuine `Bianchi identities' among \eqref{BI} are those generated by
the lifts of the canonical framing $\{E_0, E_1, \ldots , E_{2n} \}$
of  $TH_n$ to left-invariant vector fields  $\{X_0, X_1, \ldots ,
X_{2n} \}$
 on $G_{\rm cn}$.
\medskip

\begin{proposition} \label{lifts}
The left-invariant vector fields on $G_{\rm cn}$ corresponding to
the canonical basis of the Heisenberg Lie algebra are as follows:
 \begin{align} \label{X0}
& X_0  \mid_{(\xb, \, A, \,  s)} \, = \, s^2\, \frac{\p}{\p x^0}   \, = \, s^2\, E_0 \, , \\ \label{Xj}
&X_j  \mid_{(\xb, \, A, \,  s)} \, = \, s \sum_{i=1}^{2n} a^i_j\, E_i \,  , \quad
\qquad 1 \leq j \leq 2n.
\end{align}
\end{proposition}

 \proof We start with the lift of $E_0$. Since
$\, \exp (t \eb_0) = t \eb_0 $, one has for any $F \in C^\infty (G_{\rm cn})$,
\begin{align*}
&X_0 F (\xb, A, s) =
\frac{d}{dt} \mid_{t=0} F ( (\xb, A, s)  \cdot (t \eb_0 , \Id, 1)) =\\
&=\frac{d}{dt} \mid_{t=0} F
 (\xb \ast t  S_A (\mu_s (\eb_0)) ,\, A,\,  s) ;
 \end{align*}
as $\, S_A (\mu_s (\eb_0)) = s^{2} \eb_0$, we can continue as follows:
\begin{align*}
&= \frac{d}{dt} \mid_{t=0} F (\xb \ast t  s^2 \, \eb_0 , \,  A, \, s)
=\frac{d}{dt} \mid_{t=0} F
( x^0+ t  s^2, \, x^1, \ldots , x^{2n}) ,\, A, \,  s) = \\
& \quad = s^2  \, \frac{\p F}{\p x^0} (\xb, \, A, \,  s) .
\end{align*}
This proves \eqref{X0}.

Next, for $1\leq j \leq 2n$, let $X_j $ denote the lift of $E_j$ to $G_{\rm cn}$.
Again,  using that $\, \exp (t \eb_j) = t \eb_j $ in $H_n$, one has
\begin{align*}
&X_j F (\xb, A, s) =
\frac{d}{dt} \mid_{t=0} F ( (\xb, A, s)  \cdot (t \eb_j , \Id, 1)) =\\
&=\frac{d}{dt} \mid_{t=0} F
 (\xb \ast t  S_A (\mu_s (\eb_j)) ,\, A,\,  s) ;
  \end{align*}
because  $\, S_A (\mu_s (\eb_j)) = s \ab_j$, with $ \ab_j$ denoting the $j$th
column in the matrix $A$, the above is equal to
\begin{align*}
&  = \frac{d}{dt} \mid_{t=0} F (\xb \ast t  s \,\ab_j, \,  A, \, s)
\frac{d}{dt} \mid_{t=0} F( (x^0 + t  s \b (\xb' ,\ab_j), \xb' + t s \ab_j ), \,  A, \, s)  \\
&=\frac{d}{dt} \mid_{t=0} F
( x^0+ \frac{t  s}{2}\sum_{i=1}^n (x^i a^{n+i}_j - a^i_j x^{n+i}),
x^1+ ts a^1_j , \ldots ,  x^{2n} + ts a^1_{2n}) , A, s) \\
&= \frac{s}{2}  \sum_{i=1}^n \left(x^i a^{n+i}_j - a^i_j x^{n+i}\right) \,
\frac{\p F}{\p x^0} (\xb, \, A, \,  s) +
s  \sum_{k=1}^{2n} a_{k, j} \frac{\p F}{\p x^k} (\xb, \, A, \,  s).
\end{align*}
Thus,
\begin{align*}
&X_j  \mid_{(\xb, \, A, \,  s)} \, = \,
\frac{s}{2}  \sum_{i=1}^n \left(x^i a^{n+i}_j - a^i_j x^{n+i}\right) \, \frac{\p}{\p x^0} +
s  \sum_{k=1}^{2n} a^k_j \frac{\p }{\p x^k} =\\
&= \frac{s}{2} \sum_{i=1}^n \left(x^i a^{n+i}_j - a^i_j x^{n+i}\right) \, \frac{\p}{\p x^0} +
s  \sum_{i=1}^{n} \left( a^i_j \frac{\p }{\p x^i}
+ a^{n+i}_j \frac{\p }{\p x^{n+i}} \right) \\
&= s \sum_{i=1}^n a^i_j \left(\frac{\p }{\p x^i} - \frac{1}{2} x^{n+i} \frac{\p}{\p x^0} \right)
+ s \sum_{i=1}^n a^{n+i}_j \left(  \frac{\p }{\p x^{n+i}} +\frac{1}{2} x^i  \frac{\p}{\p x^0}  \right) \\
&= \qquad  s \sum_{i=1}^n \big( a^i_j\, E_i \, + \, a^{n+i}_j \, E_{n+i} \big) ,
  \end{align*}
which is the expression in \eqref{Xj}.
 \endproof

\begin{remark} \label{tautology}
The formulae \eqref{X0}, \eqref{Xj}, which taken together are the exact analogue
of the formula \eqref{horiz},
 simply express the fact that the transition
matrix from the basis $\{E_0, E_1, \ldots , E_{2n}\}$ to the basis
$\{X_0, X_1, \ldots , X_{2n} \}$ of the horizontal subspace of
$T_{(\xb , A , s)} G_{\rm cn} \cong T_\xb H_n$ is precisely the matrix
\end{remark}
\begin{equation*}
 \left(\begin{array}{c|c}   s^2&   0\,  \ldots \,0\\
\hline\\
0& \\
\vdots& s A\\
0&
\end{array}\right).
\end{equation*}
\medskip

We now give a few examples the cocycles $ \G_i^j (\phi)  \in C^\infty (G_{\rm cn})$,
$\phi \in \Gb_{\rm cn}$,
corresponding to the horizontal vector fields $\{X_0, X_1, \ldots , X_{2n} \}$.
 \medskip

Starting with $X_0 $, one has
\begin{align*}
&(U^\ast_\phi  X_0 U_\phi )\, F (\xb, A, s) =  (U^\ast_\phi  s^2 U_\phi )(U^\ast_\phi
\frac{\p}{\p x^0}  U_\phi ) \, F (\xb, A, s)=\\
&=  t_\phi(\xb)^2 s^2\\ &
\frac{\p}{\p y^0}
 F \big(\phi^{-1}(\yb) ,  t_{\phi^{-1}} (\yb)^{-1}  A_{\phi^{-1}} (\yb) B,
  t_{\phi^{-1}} (\yb) z \big) \mid_{(\phi (\xb) , \, t_\phi(\xb)^{-1} A_\phi(\xb) A,\,  t_\phi(\xb) s)}
  =\\
&=  t_\phi(\xb)^2 s^2 \,\Big(\p_{\xb} F \cdot \, \frac{\p (\phi^{-1})}{\p y^0}  \, + \,
\p_A F \cdot  \frac{\p}{\p y^0} \big(t_{\phi^{-1}} (\yb)^{-1}  A_{\phi^{-1}} (\yb)\big) B
\, +   \\
& +\p_s F  \frac{\p}{\p y^0} \big( t_{\phi^{-1}} (\yb) z\big) \Big)
\mid_{(\phi (\xb) , \, t_\phi(\xb)^{-1} A_\phi(\xb) A,\,  t_\phi(\xb) s)}
\,  =\\
& = t_\phi(\xb)^2 s^2 \,\Big(\p_{\xb} F \cdot \, \frac{\p
(\phi^{-1})}{\p y^0} - \, t_{\phi^{-1}} (\yb)^{-2} \frac{\p}{\p y^0}
 \big(t_{\phi^{-1}} (\yb) \big) \p_A F \cdot \big(A_{\phi^{-1}} (\yb) B \big) + \\
&+ t_{\phi^{-1}} (\yb)^{-1}  \p_A F \cdot
 \frac{\p}{\p y^0} \big(A_{\phi^{-1}} (\yb)\big) B+\\
&+\p_s F  \frac{\p}{\p y^0} \big( t_{\phi^{-1}} (\yb) \big) z \Big)
\mid_{(\phi (\xb) , \, t_\phi(\xb)^{-1} A_\phi(\xb) A,\,
t_\phi(\xb) s)}.
 \end{align*}
Taking into account that
\begin{align} \label{aux1}
 &  t_{\phi^{-1}} (\phi (\xb)) \, t_\phi(\xb) = 1 , \qquad
 A_{\phi^{-1}} (\phi (\xb)) \, A_\phi(\xb) = \Id , \\  \label{aux2}
 &\text{and} \qquad \frac{\p (\phi^{-1})^0}{\p x^0} (\phi (\xb)) = t_{\phi^{-1}} (\phi (\xb))^2 =
 t_\phi(\xb)^{-2} ,
\end{align}
one obtains after evaluation at
$(\phi (\xb) , \, t_\phi(\xb)^{-1} A_\phi(\xb) A,\,  t_\phi(\xb) s) \in G_{\rm cn}$
\begin{align} \notag
&U^\ast_\phi  X_0 U_\phi \, =\\
& = \, X_0  + t_\phi(\xb)^2 s^2
\,\sum_{i=1}^{2n} \frac{\p (\phi^{-1})^i}{\p x^0} (\phi (\xb)) \,
\frac{\p }{\p x^i} - t_\phi(\xb)^3 s^2 \,\frac{\p t_{\phi^{-1}}}{\p
x^0} (\phi (\xb))\, \p_A  \cdot A \\ \label{G0} & + \, t_\phi(\xb)^2
s^2 \, \p_A  \cdot \frac{\p A_{\phi^{-1}}}{\p x^0} (\phi (\xb))
A_\phi(\xb) A \, +\,
 \frac{\p  t_{\phi^{-1}} }{\p x^0}  (\phi (\xb))\,t_\phi (\xb)^3 s^3  \frac{\p}{\p s}.
\end{align}
To find the cocycles
of the form $\G_i^0$'s, with $i=1, \ldots , 2n$,
we use \eqref{Ei} to replace the partial derivatives by the horizontal vector fields,
 \begin{equation}  \label{invEi}
 \p_i = E_i + \frac{1}{2} x^{n+i} E_0, \qquad
\p_{n+i} = E_{n+i} - \frac{1}{2}x^{i} E_0 , \qquad1 \leq i \leq n ,
\end{equation}
and rewrite
the second term in the right hand side of \eqref{G0} as follows
 \begin{align*}
&II_{\rm term} =t_\phi (\xb)^2 \, s^2 \left( \sum_{i=1}^{n}
\frac{\p (\phi^{-1})^i}{\p x^0} (\phi (\xb)) \,  \frac{\p }{\p x^i} +
 \sum_{i=n+1}^{2n} \frac{\p (\phi^{-1})^i}{\p x^0} (\phi (\xb)) \,  \frac{\p }{\p x^i}\right) =\\
 &t_\phi (\xb)^2 \, s^2 \sum_{i=1}^{n}
\Big( \frac{\p (\phi^{-1})^i}{\p x^0} (\phi (\xb)) \, (E_i +
\\
& +\frac{1}{2} x^{n+i} E_0) +
 \frac{\p (\phi^{-1})^i}{\p x^0} (\phi (\xb)) \, (E_{n+i} - \frac{1}{2}x^{i} E_0 )\Big) \\
&= t_\phi(\xb)^2 \,
 \b \big(\frac{\p (\phi^{-1})^\prime}{\p x^0} (\phi (\xb)) , \, \xb^\prime \big) \, X_0  \, + t_\phi(\xb)^2  \,  s^2\, \sum_{i=1}^{2n}  \frac{\p
(\phi^{-1})^i}{\p x^0} (\phi (\xb)) \, E_i .
\end{align*}
 We next invert the formula \eqref{Xj}, \cf Remark \ref{tautology},
 \begin{align} \label{Ej}
E_i \mid_{(\xb, \, A, \,  s)} \, = \, s^{-1} \sum_{j=1}^{2n} \check{a}^j_i\, X_j \,  , \quad
\qquad 1 \leq i \leq 2n ,
\end{align}
where $( \check{a}^i_j) = A^{-1}$, to obtain
\begin{align*}
&II_{\rm term}
= t_\phi(\xb)^2  \,
 \b \big(\frac{\p (\phi^{-1})^\prime}{\p x^0} (\phi (\xb)) , \, \xb^\prime \big) \, X_0  \, + \,
  t_\phi(\xb)^2  \,  s \,
\sum_{i, j=1}^{2n}  \frac{\p (\phi^{-1})^i}{\p x^0} (\phi (\xb)) \,  \check{a}^j_i\, X_j .
\end{align*}

We have thus shown that
\begin{align} \label{UX00U}
&\G_0^0 (\phi)(\xb, \, A, \,  s)\, = \, \Id \, +     \, t_\phi(\xb)^2  \,
 \b \big(\frac{\p (\phi^{-1})^\prime}{\p x^0} (\phi (\xb)) , \, \xb^\prime \big) , \\ \label{UX0jU}
&\G_0^i (\phi)(\xb, \, A, \,  s)\,=\,  t_\phi(\xb)^2  \,  s\,
\sum_{j=1}^{2n}  \frac{\p (\phi^{-1})^j}{\p x^0} (\phi (\xb)) \,  \check{a}^i_j , \quad i=1, \ldots , 2n .
\end{align}
In particular, when restricted to $\psi \in N_{\rm cn}$
and evaluated  at the neutral element $e= (1, \, \Id, \,  1) \in G_{\rm cn}$,
these cocycles take the simple form
\begin{align} \label{a00psi}
&\G_0^0 (\psi)(0, \, \Id, \,  1)\,=\,  \Id  ,  \\ \label{a0jpsi}
&\G_0^j (\psi)(0, \, \Id, \,  1)\,=\,  \frac{\p (\psi^{-1})^j}{\p x^0} (0) , \qquad j=1, \ldots , 2n .
\end{align}

By comparison with the flat case, these cocycles and their derivatives
 give the only new type of coordinate
functions on the group $\psi \in N_{\rm cn}$, all the rest being completely
analogous to the $\eta^\bullet_{\bullet \ldots \bullet}$ coordinates of \eqref{fcoord}.
\medskip

One last ingredient needed for the construction of the Hopf algebra, is
provided by the following lemma.

\begin{lemma}
The left Haar volume form of the group $G_{\rm cn}$ is invariant
under the action $\rt$ of $\Gb_{\rm cn}$.
\end{lemma}

\proof Up to a constant factor, the left-invariant volume form of $G_{\rm cn}$
is given, in the coordinates \eqref{param}, by
\begin{align} \label{volfor}
\varpi_{\rm cn} := \, \a \wg d \a^n \wg \varpi_{\rm Sp} \wg s^{-2(n+1)} \frac{ds}{s} ,
\end{align}
where $\varpi_{\rm Sp}$ is the  left-invariant volume form of $Sp(N, \Rb)$.
Using the formula \eqref{rtact}
expressing the action of
$\, \phi \in \Gb_{\rm cn}$ on $G_{\rm cn}$, in conjunction  with the left invariance of
$\varpi_{\rm Sp}$  and the fact that
$$
\phi^*(\a \wg d \a^n ) \, = \, t_\phi^{2(n+1)} \a \wg d \a^n ,
$$
one immediately sees that $\, \tilde\phi^* (\varpi_{\rm cn}) = \varpi_{\rm cn}$.
\endproof

As a consequence, we can define an invariant
trace $\tau = \tau_{\rm cn}$ on the crossed product algebra
 $ \Ac (\Pi_{\rm cn}) = C^\infty (G_{\rm cn}) \rtimes \Gb_{\rm cn}$
 by precisely the same formula \eqref{tr}. Furthermore,
 the following counterpart of Proposition \ref{ibp} holds.

\begin{proposition} \label{cnibp}
The infinitesimal modular character $\,  \d (X) = \Tr (\ad X)$, $\, X \in \Fg_{\rm cn}$, extends
uniquely to a character  $\d = \d_{\rm cn}$ of $\Hc(\Pi_{\rm cn})$, and the
trace $\tau = \tau_{\rm cn}$ is $\Hc(\Pi_{\rm cn})$-invariant relative to this
character, \ie
\begin{equation} \label{cnit}
\tau (h(a)) \, = \,  \d(h)\, \tau(a) \, , \qquad \fl \, a, b \in \Ac(\Pi_{\rm cn}).
\end{equation}
\end{proposition}

\proof On the the canonical basis of $\Fg_{\rm cn}$, the character $\d$ takes
the values
\begin{align*}
\d (E_i) = 0 , \quad 0 \leq i \leq 2n , \qquad
\d(Y_i^j) = 0 , \qquad \text{and} \qquad \d(Z) = 2n+2 ;
\end{align*}
indeed, $\ad (E_i)$'s are nilpotent, $\Ad(Y_i^j)$'s are unimodular, and
\begin{align*}
[Z, E_0] = 2 E_0, \quad [Z, E_i] =  E_i , \, \fl \,  1 \leq i \leq 2n ,
\qquad  [Z, Y] = 0 , \, \fl \, Y \in \Fg_{\rm cn}.
\end{align*}
The rest of the proof is virtually identical to that of Prop. \ref{ibp}.
\endproof

Finally, following the same line of arguments which led to Theorem \ref{hopfHn},
one obtains the corresponding analog.

\begin{theorem}  \label{cnHopf}
There exists a unique Hopf algebra structure on $\Hc(\Pi_{\rm cn})$, such
that its tautological action makes $ \Ac (\Pi_{\rm cn}) $  a left module algebra.
\end{theorem}
\bigskip


\section{Bicrossed product realization}

In this section we reconstruct (or rather deconstruct)
the Hopf algebra affiliated to a primitive
Lie pseudogroup  as a bicrossed product of a matched pair of
Hopf algebras.
In the particular case of $\Hc_1$, this has been proved in~\cite{hm},
by direct algebraic calculations that rely on the detailed knowledge of
its presentation.
By contrast, our method is completely geometric
and for this reason applicable to the entire class of Lie pseudogroups
admitting a Kac-type decomposition.

We recall below the most basic notions concerning the bicrossed product construction,
referring the reader to Majid's monograph~\cite{maj} for
a detailed exposition.

\medskip

Let  $\Uc$ and $\Fc$ be two Hopf algebras.  A linear map
$$ \Db:\Uc\ra\Uc\ot \Fc , \qquad \Db u \, = \, u\ns{0} \ot u\ns{1} \, ,
$$
defines a {\em right coaction}, and thus equips
 $\Uc$ with a  {\em right $\Fc$-comodule coalgebra} structure,  if the
following conditions are satisfied for any $u\in \Uc$:
\begin{align}\label{cc1}
&u\ns{0}\ps{1}\ot u\ns{0}\ps{2}\ot u\ns{1}= u\ps{1}\ns{0}\ot
u\ps{2}\ns{0}\ot u\ps{1}\ns{1}u\ps{2}\ns{1}\\\label{cc2}
&\epsilon(u\ns{0})u\ns{1}=\epsilon(u)1.
\end{align}
One  can then  form a  cocrossed product coalgebra $\Fc\cl\Uc$, that has
$\Fc\ot \Uc$ as underlying vector space and the following coalgebra
structure:
\begin{align}\label{cocross}
&\Delta(f\cl u)= f\ps{1}\cl u\ps{1}\ns{0}\ot  f\ps{2}u\ps{1}\ns{1}\cl u\ps{2}, \\
&\epsilon(h\cl k)=\epsilon(h)\epsilon(k).
\end{align}
In a dual fashion,   $\Fc$ is called a {\em left $\Uc$ module algebra}, if
$\Uc$ acts from the left on $\Fc$ via a left action
$$
\rt : \Fc\ot \Uc \ra \Fc
$$
which satisfies the following condition for any $u\in \Uc$, and
$f,g\in \Fc$ :
\begin{align}\label{ma1}
&u\rt 1=\epsilon(u)1\\\label{ma2}
 &u\rt(fg)=(u\ps{1}\rt
f)(u\ps{2}\rt g).
\end{align}
This time  we can endow the underlying vector space $\Fc\ot \Uc$ with
 an algebra structure, to be denoted by $\Fc\al \Uc$, with  $1\al 1$
 as its unit and the product given by
\begin{equation}
(f\al u)(g\al v)=f \;u\ps{1}\rt g\al u\ps{2}v
\end{equation}

$\Uc$ and $\Fc$ are said to form a {\em matched pair}
of Hopf algebras if they are equipped, as above,  with an action and a coaction
which satisfy the following compatibility conditions:
following conditions for any $u\in\Uc$, and any $f\in \Fc$.
\begin{align}\label{mp1}
&\epsilon(u\rt f)=\epsilon(u)\epsilon(f), \\  \label{mp2}
&\Delta(u\rt f)=u\ps{1}\ns{0} \rt f\ps{1}\ot
u\ps{1}\ns{1}(u\ps{2}\rt f\ps{2}), \\  \label{mp3}
 &\Db(1)=1\ot
1, \\ \label{mp4} &\Db(uv)=u\ps{1}\ns{0} v\ns{0}\ot
u\ps{1}\ns{1}(u\ps{2}\rt v\ns{1}),\\  \label{mp5} &u\ps{2}\ns{0}\ot
(u\ps{1}\rt f)u\ps{2}\ns{1}=u\ps{1}\ns{0}\ot
u\ps{1}\ns{1}(u\ps{2}\rt f).
\end{align}
One can then form  a new Hopf algebra $\Fc\acl \Uc$, called the {\em bicrossed product}
of the matched pair  $(\Fc , \Uc)$ ; it has $\Fc\cl \Uc$ as underlying coalgebra,
$\Fc\al \Uc$ as underlying algebra and the antipode is defined by
\begin{equation}\label{anti}
S(f\acl u)=(1\acl S(u\ns{0}))(S(fu\ns{1})\acl 1) , \qquad f \in \Fc , \, u \in \Uc.
\end{equation}

\bigskip

\subsection{The flat case}
 As mentioned in the introduction, the matched pair of Hopf algebras
 arises from a matched pair of groups, via a splitting \`a la G.I. Kac~\cite{Kac}.

 \begin{proposition}
 Let $\Pi$ be a flat primitive Lie pseudogroup of infinite type, $F_{\Pi}\Rb^m$
the principal bundle of  $\Pi$-frames on $\Rb^m$. There is a
canonical splitting of the group ${\Gb} = \Diff(\Rb^m) \cap  \Pi$,
as a cartesian product $\, {\Gb} = G \cdot N$, with $\, G \simeq
F_{\Pi}\Rb^m$ the group of affine  $\Pi$-motions of $\Rb^m$, and $\,
N = \{ \phi \in {\Gb} \, ; \, \phi (0) = 0, \,  \phi' (0) = \Id \}$.
 \end{proposition}

 \proof Let $\phi \in {\Gb} $. Since ${\Gb}$ contains the translations,
 then $\, \phi_0 := \phi - \phi (0) \in {\Gb} $, and $\phi_0 (0)= 0$.
 Moreover, the affine diffeomorphism
 \begin{equation} \label{Kac1}
 \vp (x)\, := \, \phi'_0 (0) \cdot x \, + \, \phi (0) , \qquad \fl \, x \in \Rb^m
 \end{equation}
also belongs to ${\Gb} $, and has the same $1$-jet at $0$ as $\phi$.
Therefore, the diffeomorphism
\begin{equation}  \label{Kac2}
  \psi (x) \, := \,  \vp^{-1} (\phi (x)) = \phi'_0 (0)^{-1} \left(\phi (x) - \phi (0)\right),
  \quad \fl \, x \in \Rb^m
 \end{equation}
 belongs  to $N$, and the canonical decomposition is
 \begin{equation} \label{Kac}
  \phi\, = \,  \vp \circ \psi , \qquad \text{with} \quad \vp \in G \quad \text{and} \quad \psi \in N
 \end{equation}
 given by  \eqref{Kac1} and \eqref{Kac2}. The two components are uniquely
 determined, because evidently $G \cap N = \{e\}$.

 \endproof
\bigskip

Reversing the order in the above decomposition one simultaneously
obtains two well-defined operations, of $N$ on $G$ and of $G$ on
$N$:
 \begin{equation} \label{actions}
  \psi \circ \vp \, = \,  (\psi \rt \vp) \circ  (\psi \lt \vp) ,
  \qquad \text{for} \quad \vp \in G \quad \text{and} \quad \psi \in N
 \end{equation}
 \bigskip

\begin{proposition}
 The operation  $\rt$ is a left action of $N$ on $G$, and $\lt$ is
a right action of $G$ on $N$. Both actions leave the neutral element
fixed.
\end{proposition}
\smallskip

\proof Let $\psi_1, \psi_2 \in N$ and $\vp \in G$. By
\eqref{actions}, on the one hand
\begin{equation*}
 (\psi_1\circ \psi_2) \circ \vp =  \big( (\psi_1\circ \psi_2)\rt \vp \big) \circ
 \big((\psi_1\circ \psi_2) \lt \vp \big) ,
   \end{equation*}
and on the other hand
\begin{eqnarray*}
 \psi_1\circ (\psi_2 \circ \vp) &=&   \psi_1\circ (\psi_2 \rt \vp) \circ  (\psi_2 \lt \vp)\\
 &=& \big(\psi_1 \rt (\psi_2 \rt \vp)\big) \circ \big(\psi_1 \lt (\psi_2 \rt \vp)\big)
 \circ  (\psi_2 \lt \vp) .
   \end{eqnarray*}
   Equating the respective components in $G$ and $N$ one obtains:
  \begin{eqnarray} \label{action1}
(\psi_1\circ \psi_2)\rt \vp &=&  \psi_1 \rt (\psi_2 \rt \vp) \, \in
G , \qquad \text{resp.} \\ \label{reaction1}
 (\psi_1\circ \psi_2) \lt \vp &=& \big(\psi_1 \lt (\psi_2 \rt \vp)\big)
 \circ  (\psi_2 \lt \vp) \,  \in N .
   \end{eqnarray}
 Similarly,
 \begin{equation*}
 \psi \circ (\vp_1 \circ \vp_2) =  \big(\psi \rt (\vp_1 \circ \vp_2)\big) \circ
 \big(\psi \lt (\vp_1 \circ \vp_2) \big) ,
   \end{equation*}
while
\begin{eqnarray*}
( \psi \circ \vp_1) \circ \vp_2 &=& (\psi \rt \vp_1) \circ  (\psi \lt \vp_1) \circ \vp_2 \\
 &=& (\psi \rt \vp_1) \circ   \big((\psi \lt \vp_1) \rt \vp_2\big) \circ  \big((\psi \lt \vp_1) \lt \vp_2)\big)  ,
   \end{eqnarray*}
whence
  \begin{eqnarray} \label{action2}
\psi \lt (\vp_1 \circ \vp_2) &=&  (\psi \lt \vp_1) \lt \vp_2 \, \in
N , \qquad \text{resp.} \\ \label{reaction2} \psi \rt (\vp_1 \circ
\vp_2) &=&(\psi \rt \vp_1) \circ   \big((\psi \lt \vp_1) \rt
\vp_2\big)
  \,  \in G .
   \end{eqnarray}
Specializing $\vp = e$, resp. $\psi = e$, in the definition
\eqref{actions}, one sees that $e=\Id$ acts trivially, and at the
same time that both actions fix $e$.
\endproof
\bigskip

Via the identification  $\, G \simeq F_{\Pi}\Rb^m$, one can
recognize the action $\rt$ as the usual action of diffeomorphisms on
the frame bundle, {\em cf.}  \eqref{frameact}.
\smallskip

\begin{lemma} \label{rtf}
 The left action $\rt$ of $N$ on $G$ coincides with the restriction of
the natural action of $\Gb$ on $F_{\Pi}\Rb^m$.
\end{lemma}

\proof Let $\phi = \psi \rt \vp \in {\Gb}$, with $\psi \in N$ and
$\vp \in G$. The associated frame, {\it cf.} \eqref{frameid}, is $\,
\big( \phi (0), \phi' (0) \big)$.  By \eqref{actions},
\begin{equation*}
 \psi \big(\vp (0)\big) \, =\,  (\psi \rt \vp) \big((\psi \lt \vp)(0)\big) = \phi (0) ,
   \end{equation*}
since $\, (\psi \lt \vp)(0) = 0$. On differentiating
\eqref{actions} at $0$ one obtains
\begin{equation*}
 \psi' \big(\vp (0)\big) \cdot \vp'(0)\, =\,  (\psi \rt \vp)' \big((\psi \lt \vp)(0)\big) \cdot (\psi \lt \vp)'(0)
 \, =\, (\psi \rt \vp)' (0)  \, =\,  \phi' (0)  ,
   \end{equation*}
since $\, (\psi \lt \vp)'(0) = \Id$. Thus, $\, \big( \phi (0), \phi'
(0) \big) = \tilde{\psi} \big( \vp (0), \vp'(0)\big)$, as in the
definition \eqref{frameact}.
\endproof

\bigskip

\begin{definition}
 The {\em coordinates} of $\,  \psi \in N$ are
the coefficients of the Taylor expansion of $\psi$ at $0 \in \Rb^m$.
The algebra of functions on $N$  generated by these coordinates will
be denoted $\, \Fc (N) $, and its elements will be called {\em
regular functions}.
\end{definition}

Explicitly, $\Fc (N)$ is generated by the functions
\begin{equation*}
\a^i_{j j_1j_2\dots j_r}(\psi)= \p_{j_r}\dots \p_{j_1} \p_j
\psi^i(x)\mid_{x=0} , \, 1\leq i, j, j_1, j_2, \dots,  j_r \leq m ,
 \, \psi \in N;
\end{equation*}
note that  $\, \a^i_j (\psi) = \d^i_j$, because $\psi^\prime (0) =
\Id$, while for $r \geq 1$ the coefficients $\a^i_{j j_1j_2\dots
j_r} (\psi)$ are symmetric in the lower indices but otherwise
arbitrary. Thus,  $\Fc (N)$ can be viewed as the free commutative
algebra over $\Cb$ generated by the indeterminates $\{\a^i_{j
j_1j_2\dots j_r} ; \, 1\leq i, j, j_1, j_2, \dots,  j_r \leq m \}$.
\bigskip

The algebra $\Fc := \Fc (N)$ inherits from the group $N$ a canonical
Hopf algebra structure, in the standard fashion.
\smallskip

\begin{proposition} \label{Fhopf}
With the coproduct $\, \D : \Fc \ra \Fc \ot \Fc $, the antipode \\
$\, S : \Fc \ra \Fc$, and the counit $\ve : \Fc \ra \Cb$ determined
by the requirements
\begin{eqnarray} \label{Fcop}
\D(f)(\psi_1, \psi_2) &=& f(\psi_1 \circ\psi_2) , \qquad \fl \,
\psi_1, \psi_2 \in N , \\ \notag S(f)(\psi)  &=& f(\psi^{-1}) ,
\qquad \fl \, \psi \in N, \quad \fl \, f \in \Fc , \\ \notag \e (f)
&=& f(e) ,
\end{eqnarray}
  $\Fc (N)$ is a Hopf algebra.
\end{proposition}

\proof The fact that these definitions give rise to a Hopf algebra
is completely routine, once they are shown to
 make sense. In turn, checking that
\begin{equation}  \label{FDS}
\D(\a^i_{j j_1j_2\dots j_r}) \in \Fc \ot \Fc \quad \text{and} \quad
S(\a^i_{j j_1j_2\dots j_r}) \in \Fc \ot \Fc ,
\end{equation}
only involves elementary manipulations with the chain rule. For
instance, in the case of $\a^i_{jk}$ the verification goes as
follows. First, for the coproduct,
\begin{align*}
&\D(\a^i_{jk})(\psi_1,\psi_2) \, = \,
\a^i_{jk}(\psi_1\circ \psi_2)= \p_j\p_k(\psi_1\circ\psi_2)^i(x)\mid_{x=0}\\
&=\p_j\big((\p_\mu\psi^i_1)(\psi_2(x)) \p_k\psi^\mu_2(x)\big)\mid_{x=0}\\
 &= (\p_\nu\p_\mu\psi^i_1)(\psi_2(x))  \mid_{x=0} \p_j \psi^\nu_2(x)\
 \mid_{x=0}\p_k\psi^\mu_2(x)\mid_{x=0}+ \\
 &\qquad +(\p_\mu(\psi^i_1
 (\psi_2(x))\mid_{x=0}\p_j\p_k\psi^\mu_2(x))\mid_{x=0}\\
 &= \p_j\p_k\psi^i_1(x)\mid_{x=0}+ \p_j\p_k\psi^i_2(x)\mid_{x=0}\, = \, (\a^i_{j,k}\ot 1+1\ot
 \a^i_{j,k})(\psi_1,\psi_2) ,
\end{align*}
where we have used that $\psi_1(0)=\psi_2(0)=0$ and
$\psi_1'(0)=\psi'_2(0)=\Id$.

To deal with the antipode, one differentiates the identity
$\psi^{-1}(\psi (x))=x$:
\begin{align*}
\d^i_j = \p_j \big((\psi^{-1})^i(\psi (x)) \big) =
\p_\lb(\psi^{-1})^i (\psi(x)) \, \p_j\psi^\lb(x) ,
\end{align*}
which yields under further differentiation
 \begin{align*}
\p_\mu\p_\lb(\psi^{-1})^i (\psi(x)) \, \p_k \psi^\mu (x)\,
 \p_j\psi^\lb(x)+\p_\lb(\psi^{-1})^i(\psi(x))\,
\p_k\p_j\psi^\lb(x)=0;
\end{align*}
evaluation at $x=0$ gives $ \quad \a^i_{jk}(\psi^{-1}) \,  + \,
\a^i_{jk} (\psi) \, = \, 0$.

Taking higher derivatives one proves \eqref{FDS} in a similar
fashion.
   \endproof
   \bigskip

We shall need an alternative description of the algebra $\Fc$, which
will be used to recognize it as being identical to the Hopf
subalgebra of $\Hc(\Pi)$ generated by the $\d^i_{j k \ell_1\dots
\ell_r}$'s.
   \medskip

\begin{lemma} \label{regular}
 The coefficients of the Taylor expansion of $\tpsi$ at $e \in G$,
\begin{equation} \label{fcoord}
 \eta^i_{j k \ell_1 \ldots \ell_r} (\psi): =  \g^i_{j k \ell_1 \ldots
\ell_r} (\psi)(e) , \qquad \psi \in N ,
 \end{equation}
define regular functions on $N$, which generate the algebra $\, \Fc
(N)$.
\end{lemma}

\proof Evaluating the expression \eqref{highg} at $e=(0, \Id) \in
F_{\Pi}\Rb^m$ gives
\begin{equation} \label{highd0}
\g_{jk \,  \ell_1 \ldots \ell_r}^i (\psi) (e) \, = \, \part_{\ell_r}
\ldots \part_{\ell_1} \left(({\psi}^{\prime} (x)^{-1})^i_\nu
 \part_j \part_k {\psi}^\nu (x) \right) \mid_{x=0}.
\end{equation}
The derivatives of ${\psi}^{\prime} (x)^{-1}$ are sums of terms each
of which is a product of derivatives of  ${\psi}^{\prime} (x)$
interspaced with ${\psi}^{\prime} (x)^{-1}$ itself. Since $\psi' (0)
=\Id$, the right hand side of \eqref{highd0} is thus seen to define
a regular function on $N$.

A more careful inspection actually proves the converse as well.
First, by the very definition,
\begin{equation} \label{aeta1}
 \eta^i_{jk} \, =\, \a^i_{jk} \, .
  \end{equation}
Next, one has
\begin{align*}
&\eta^i_{jk\ell}(\psi) \, = \,
\p_\ell\left(({\psi}^{\prime} (x)^{-1})^i_\nu \p_j \p_k {\psi}^\nu (x) \right) \mid_{x=0} \\
&\, =\, \p_\ell\left(({\psi}^{\prime} (x)^{-1})^i_\nu
\right)\mid_{x=0} \p_j \p_k {\psi}^\nu (x) \mid_{x=0} \, + \,
\p_\ell \p_k \p_j {\psi}^i (x) \mid_{x=0} \, ;
\end{align*}
 on differentiating $\, ({\psi}^{\prime} (x)^{-1})^i_\mu \, \p_\nu \psi^\mu (x) = \d^i_\nu$
one sees that
\begin{equation*}
 \p_\ell\left(({\psi}^{\prime} (x)^{-1})^i_\nu \right)\mid_{x=0} \,  +\,
 \p_\ell \p_\nu\psi^i (x)\mid_{x=0} \, = \, 0 ,
 \end{equation*}
 and therefore
\begin{equation}  \label{aeta2}
\eta^i_{jk\ell}(\psi) \, = \, \a^i_{jk\ell}(\psi)\, -\,
\a^i_{\ell\nu}(\psi)\, \a^\nu_{jk}(\psi) .
 \end{equation}
By induction, one shows that
\begin{equation}  \label{aeta3}
\eta^i_{j k \ell_1 \ldots \ell_r}=\a^i_{j k \ell_1 \ldots \ell_r} \,
+ \, P^i_{j k \ell_1 \ldots \ell_r} (\a^\lambda_{\mu \nu}, \dots ,
\a^{\rho}_{\tau \s p_1 \dots p_{r-1}}) ,
\end{equation}
where $P^i_{j k \ell_1 \ldots \ell_r}$ is a polynomial. The
triangular form of the identities \eqref{aeta1}-\eqref{aeta3} allows
to reverse the process and express the $ \a^i_{j k \ell_1 \ldots
\ell_r} $'s in a similar fashion:
\begin{equation}  \label{aeta4}
\a^i_{j k \ell_1 \ldots \ell_r}=\eta^i_{j k \ell_1 \ldots \ell_r} \,
+ \, Q^i_{j k \ell_1 \ldots \ell_r} (\eta^\lambda_{\mu \nu}, \dots ,
\eta^{\rho}_{\tau \s p_1 \dots p_{r-1}}) .
\end{equation}
\endproof
\bigskip

Let $\Hc(\Pi)_{\rm ab}$ denote the (abelian) Hopf subalgebra of
$\Hc(\Pi)$ generated by the operators
 $\{ \d^i_{j k \ell_1 \ldots \ell_r} ; \, 1\leq i, j,k, \ell_1,  \dots,  \ell_r \leq m \}$.
\medskip

\begin{proposition} \label{isoiota}
There is a unique isomorphism of Hopf algebras \\
$ \iota : \Hc(\Pi)^{\rm cop}_{\rm ab} \ra \Fc (N)$ with the property
that
\begin{equation} \label{iota}
\iota (\d^i_{j k \ell_1 \ldots \ell_r}) \, = \, \eta^i_{j k \ell_1
\ldots \ell_r} , \qquad \fl \, 1\leq i, j,k, \ell_1,  \dots,  \ell_r
\leq m \ .
\end{equation}
\end{proposition}
\smallskip

\proof In view of \eqref{aeta1}-\eqref{aeta2}, the generators $
\eta^i_{j k \ell_1 \ldots \ell_r} $ satisfy the analogue of the
Bianchi  identity  \eqref{bianchi}. Indeed,
\begin{equation*}
\eta^i_{jk \ell}-\eta^i_{j\ell k}= \a^i_{jk \ell}-\a^i_{\ell
\rho}\a^\rho_{jk}-\a^i_{j \ell k}+\a^i_{k\rho}\a^\rho_{j\ell}=
\eta^i_{k\rho}\eta^\rho_{j\ell}-\eta^i_{\ell \rho}\eta^\rho_{jk} .
\end{equation*}
From Theorem \ref{hopfGHn} (or rather the proof of Corollary
\ref{quotient}) it then follows that the assignment \eqref{iota}
does give rise to a well-defined algebra homomorphism $\,  \iota :
\Hc(\Pi)_{\rm ab} \ra \Fc (N)$, which by Lemma \ref{regular} is
automatically surjective.

To prove that $\,  \iota : \Hc(\Pi)_{\rm ab} \ra \Fc (N)$ is
injective, it suffices to show that the monomials $\{\eta_K ; K =
\text{increasingly ordered multi-index}\}$, defined in the same way
as the $\d_K$'s of the Poincar\'e-Birkhoff-Witt basis of $\Hc(\Pi)$
({\em cf.} Proposition \ref{free}), are linearly independent. This
can be shown by induction on the height. In the height $0$ case the
statement is obvious, because of \eqref{aeta1}.
 Next, assume
 \begin{equation*}
 \sum_{|J| \leq N-1} \, c_J\,  \eta_J +  \sum_{|K| = N} \, c_K \,  \eta_K = 0 .
  \end{equation*}
Using the identities \eqref{aeta3} and \eqref{aeta3}, one can
replace $ \eta_K$ by $ \a_K +$ {\em lower height}. Since the
$\a^{\bullet}_{\bullet \ldots \bullet}$'s are free generators, it
follows that  $c_K =0$ for each $K$ of height $N$, and thus we are
reduced to
 \begin{equation*}
 \sum_{|J| \leq N-1} \, c_J\,  \eta_J  = 0 ;
  \end{equation*}
the induction hypothesis now implies $\, c_J =0$, for all $J$'s.
\smallskip

It remains to prove that $\,  \iota :  \Hc(\Pi)^{\rm cop}_{\rm ab}
\ra \Fc (N)$ is a coalgebra map, which amounts to checking that
\begin{equation} \label{Diota}
\iota \ot \iota (\D \d^i_{j k \ell_1 \ldots \ell_r}) \, = \, \D^{\rm
op} \eta^i_{j k \ell_1 \ldots \ell_r} .
\end{equation}
Recall, {\em cf.} \eqref{D}, that $\D : \Hc(\Pi) \ra \Hc(\Pi) \ot
\Hc(\Pi)$
 is determined by a Leibniz rule, which for $\, \d^i_{j k \ell_1 \ldots \ell_r}$
 takes the form
  \begin{equation*}
  \d^i_{j k \ell_1 \ldots \ell_r} (U^*_{\phi_1} U^*_{\phi_2})
 \,= \,  \sum c_j^{iAB} \, \d^i_{j A} (U^*_{\phi_1} ) \, \d^i_{j B} (U^*_{\phi_2} ) ,
 \qquad \phi_1, \phi_2 \in {\Gb} ,
   \end{equation*}
which is equivalent to
 \begin{equation}\label{Dg}
  \g^i_{j k \ell_1 \ldots \ell_r} (\phi_2 \circ \phi_1)\, =\,
 \sum c_j^{iAB}  \, \g^i_{j A} (\phi_1) \, \g^i_{j B} (\phi_2)\circ \tilde{\phi}_1 .
  \end{equation}
Restricting \eqref{Dg} to $\psi_1, \psi_2 \in N$ and evaluating at
$e \in G$, one obtains
\begin{equation*}
\D^{\rm op} \eta^i_{j k \ell_1 \ldots \ell_r} (\psi_1 , \psi_2) :=
\, \eta^i_{j k \ell_1 \ldots \ell_r} (\psi_2 \circ \psi_1)\, = \,
 \sum c_j^{iAB}  \, \eta^i_{j A} (\psi_1) \, \eta^i_{j B} (\psi_2) .
\end{equation*}
 \endproof
\bigskip

  The right action $\lt $ of $G$ on $N$ induces an action of $G$ on $\Fc (N)$,
  and hence a left action $\rt $ of $\Uc(\Fg)$ on $ \Fc (N)$, defined by
 \begin{eqnarray}\label{u>f}
(X \rt f)(\psi) = \frac{d}{dt}\mid_{t=0} f (\psi \lt \exp tX) ,
\qquad f \in \Fc , \quad X \in \Fg.
\end{eqnarray}
  On the other hand, there is a natural action of $\Uc(\Fg)$ on
  $ \Hc(\Pi)_{\rm ab} $,
  induced by the adjoint action of $\Fg$ on $\Fh (\Pi)$, extended as
 action by derivations on the polynomials in
 $ \d^i_{j k \ell_1 \ldots \ell_r}$'s.
  In order to relate these two actions, we need a preparatory lemma.
  \medskip

  \begin{lemma} \label{ltact}
   Let $\vp \in G$ and $\phi \in \Gb$. Then for any $ 1\leq i, j,k, \ell_1,  \dots,  \ell_r \leq m$,
   \begin{align} \label{ltact1}
   &\g^i_{j k \ell_1 \ldots \ell_r} (\vp \circ \phi)\, = \, \g^i_{j k \ell_1 \ldots \ell_r} ( \phi) ,\\  \label{ltact2}
   &\g^i_{j k \ell_1 \ldots \ell_r} (\phi \circ \vp)\, = \,
   \g^i_{j k \ell_1 \ldots \ell_r} ( \phi) \circ \tilde{\vp}  .
   \end{align}
  \end{lemma}
\smallskip

 \proof  Both identities can be verified by direct computations,
  using the explicit formula
\eqref{highg} for $\g^i_{j k \ell_1 \ldots \ell_r} $, in conjunction
with the fact that $\vp$ has the simple affine expression $\vp (x) =
{\ab} \cdot x+b$, $\, {\ab} \in G_0(\Pi)$, $b \in \Rb^m$.

 An alternative and more elegant explanations relies on
  the left invariance of the vector fields $X_k$, {\em cf.}
Proposition \ref{leftinv}. The identity \eqref{ltact1} easily
follows from the cocycle property \eqref{gcocy} and the fact that
$\vp$  is affine,
  \begin{align*}
\g^i_{j k} (\vp \circ \phi) \, = \, \g^i_{j k \ell_1 \ldots \ell_r}
( \vp) \circ \tpsi \, +
   \, \g^i_{j k} ( \phi)\, = \, \g^i_{j k} ( \phi) \, ,
 \end{align*}
 because $\, \g^i_{j k } (\vp) \, = \,0$.
To check the second equation one starts with
  \begin{align*}
\g^i_{j k} (\phi \circ \vp) \, = \, \g^i_{j k} ( \phi) \circ
\tilde{\vp} \,+ \, \g^i_{j k } (\vp) \,
 = \, \g^i_{j k} ( \phi) \circ \tilde{\vp} ,
 \end{align*}
and notice that the invariance property $\, U_{\vp} \, X \,
U^*_{\vp} \, = \, X$, for any $X \in \Fg$, implies
 \begin{align*}
X \big(\g^i_{j k} (\phi ) \circ \tilde{\vp} \big) \, = \, X
\big(\g^i_{j k} (\phi)\big) \circ \tilde{\vp} .
 \end{align*}
 \endproof
 \medskip

We are now in a position to formulate the precise relation between
the canonical action of $\, \Uc(\Fg)$ on  $\, \Hc(\Pi)_{\rm ab}$ and
the action $\rt$  on $\, \Fc (N)$.
   \medskip

  \begin{proposition}\label{moda}
  The algebra isomorphism $\, \iota : \Hc(\Pi)_{\rm ab} \ra \Fc (N)$
  identifies the  $\,\Uc(\Fg)$-module  $\Hc(\Pi)_{\rm ab}$
  with the $\,\Uc(\Fg)$-module $\Fc (N)$. In particular $\, \Fc (N)$ is
  $\,\Uc(\Fg)$-module algebra.
  \end{proposition}
\smallskip

\proof  We denote below by $\vp_t$ the $1$-parameter subgroup $\exp
tX$ of $G$ corresponding to $X \in \Fg$, and employ the abbreviated
notation $\eta = \eta^i_{j k \ell_1 \ldots \ell_r}$,  $\, \g =
\g^i_{j k \ell_1 \ldots \ell_r}$. From \eqref{ltact1} it follows
that
\begin{equation*}
 \g (\psi \lt \vp_t) \, = \, \g (\psi \circ \vp_t) ,
\end{equation*}
whence
\begin{eqnarray*}
(X \rt \eta)(\psi) = \frac{d}{dt}\mid_{t=0} \eta (\psi \lt \vp_t) \,
= \,
  \frac{d}{dt}\mid_{t=0} \g (\psi \circ \vp_t) (e) .
\end{eqnarray*}
Now using  \eqref{ltact2}, one can continue as follows:
\begin{eqnarray*}
\frac{d}{dt}\mid_{t=0} \g (\psi \circ \vp_t) (e) \, = \,
 \frac{d}{dt}\mid_{t=0}\g (\psi)( \tilde{\vp_t} (e))  \, = \, X\big(\g (\psi)\big) (e) .
\end{eqnarray*}
By iterating this argument one obtains, for any $u \in \Uc(\Fg)$,
 \begin{equation}  \label{giota}
(u \rt \eta^i_{j k \ell_1 \ldots \ell_r}) (\psi)\, = \,u
\big(\g^i_{j k \ell_1 \ldots \ell_r} (\psi)\big)(e) , \qquad \psi
\in N .
\end{equation}
The right hand side of \eqref{giota}, before evaluation at  $e \in
G$, describes the effect of the action of $u \in \Uc(\Fg)$ on
$\d^i_{j k \ell_1 \ldots \ell_r} \in \Hc(\Pi)_{\rm ab}$. In view of
the defining relation \eqref{iota} for the isomorphism $\iota$, this
achieves the proof.
\endproof
\bigskip

 We proceed to equip $\Uc(\Fg)$ with a right $\Fc (N)$-comodule structure.
 To this end, we assign to each
 element $\, u \in \Uc(\Fg)$ a $\Uc(\Fg)$-valued function on $N$ as follows :
   \begin{equation} \label{comod3}
({\Db} u)  (\psi) \, = \,  \tilde{u} (\psi) (e) , \quad \text{where}
\quad
 \tilde{u} (\psi) \, = \,  U_\psi \, u\, U^{\ast}_\psi .
 \end{equation}
 We claim that  ${\Db} u$ belongs to $ \Uc(\Fg) \ot \Fc (N)$, and therefore
the above assignment defines a linear map $\,\Db:\Uc(\Fg)
\ra\Uc(\Fg) \ot \Fc (N)$.
 Indeed, let $\{Z_I\}$ be the PBW basis  of $\Uc(\Fg)$ defined in \S 5, {\em cf.}
 \eqref{pbw}. We identify $\Uc(\Fg)$ with the algebra of left-invariant
 differential operators on $G$, and regard the $Z_I$'s as a linear basis
 for these operators. In particular, one can uniquely express
 \begin{equation} \label{comod1}
 U_\psi \, Z_I \, U^{\ast}_\psi \, = \, \sum_J  \b_I^J (\psi)\, Z_J, \qquad \psi \in N ,
 \end{equation}
with $\b_I^J  (\psi)$ in
  the algebra of functions on $G$ generated by $\{ \g^i_{j K} (\psi) \}$.
The definition \eqref{comod3} then takes the explicit form
  \begin{equation} \label{comod2}
 \Db Z_I = \sum_J Z_J \ot \zeta_I^J  , \qquad \text{where} \quad
  \zeta_I^J (\psi) = \b_I^J  (\psi) (e) .
 \end{equation}
For example, by \eqref{Ydisp}, \eqref{Xphi} and  \eqref{gcocy}, one
has
\begin{align}\label{Ycoaction}
&\Db Y^i_j= Y^i_j\ot 1 , \\ \label{Xcoaction} &\Db X_k= X_k\ot 1
+Y^j_i\ot\eta^i_{jk} .
\end{align}
 Thus, $\,\Db:\Uc(\Fg) \ra\Uc(\Fg) \ot \Fc (N)$ is well-defined.
 \medskip

 \begin{lemma} \label{cms}
The map $\,\Db:\Uc(\Fg) \ra\Uc(\Fg) \ot \Fc (N)$ endows $\Uc(\Fg)$
with a $\Fc (N)$-comodule structure.
 \end{lemma}

 \proof On the one hand, applying \eqref{comod1} twice one obtains
\begin{align} \notag
&U_{\psi_1} U_{\psi_2}  Z_I  U^{\ast}_{\psi_2} U^{\ast}_{\psi_1} \,
=\,
  \sum_{J}  \b_I^J (\psi_2)\circ \psi_1^{-1}\,  U_{\psi_1} \, Z_J  U^{\ast}_{\psi_1} \\ \label{coac1}
 &\qquad \qquad =
 \,  \sum_{K} \left(\sum_{J} \b_I^J (\psi_2)\circ \psi_1^{-1}\,  \b_J^K (\psi_1)\right) \, Z_K  ,
  \end{align}
 while on the other hand, the same left had side can be expressed as
 \begin{align}  \label{coac2}
U_{\psi_1  \psi_2}  Z_I   U^{\ast}_{\psi_1  \psi_2} \, =\, \sum_K
\b_I^K (\psi_1  \psi_2)\, Z_K ;
\end{align}
therefore
\begin{align}  \label{coac3}
 \b_I^K (\psi_1  \psi_2) \, =\,\sum_{J} \b_J^K (\psi_1)\,
 \b_I^J (\psi_2)\circ \psi_1^{-1} .
\end{align}
By the very definition \eqref{comod2}, the identity \eqref{coac1}
gives
\begin{align*}
({\Db} \ot \Id) ( \Db Z_I ) \, = \, \sum_{K}  Z_K \ot \sum_{J}
\zeta_J^K \ot \zeta_I^J ,
\end{align*}
 while the definition \eqref{Fcop} and  \eqref{coac3} imply
\begin{align*}
\D  \zeta_I^K  \, = \,  \sum_{J} \zeta_J^K \ot \zeta_I^J .
\end{align*}
One concludes that
\begin{align*}
({\Db} \ot \Id) ( \Db Z_I ) \, = \, \sum_{K}  Z_K \ot \Db  \zeta_I^K
\, = \, (\Id \ot \D ) ( \Db Z_I ) .
\end{align*}
 \endproof
\bigskip

\begin{proposition} \label{comcoa}
Equipped with the coaction $\,\Db:\Uc(\Fg) \ra\Uc(\Fg) \ot \Fc (N)$,
$\Uc (\Fg)$ is a  right $\Fc(N)$-comodule coalgebra.
\end{proposition}

\proof It is obvious from the definition that, for any $u\in \Uc
(\Fg)$,
\begin{align} \label{comcoa1}
\ve(u\ns{0})u\ns{1} \, = \, \ve(u)1 .
\end{align}
We just have to check that
\begin{align} \label{comcoa2}
u\ns{0}\ps{1}\ot u\ns{0}\ps{2}\ot u\ns{1}= u\ps{1}\ns{0}\ot
u\ps{2}\ns{0}\ot u\ps{1}\ns{1}u\ps{2}\ns{1} .
\end{align}
In terms of the alternative definition \eqref{comod3}, this amounts
to showing that
\begin{align} \label{comcoa3}
\D \big(({\Db} u) (\psi)\big) \, = \, \Db (\D u) (\psi) , \qquad \fl
\, \psi \in N ,
 \end{align}
where
  \begin{equation*}
({\Db} (\D u))  (\psi) := \,  \widetilde{\D u} (\psi) (e, e) ,
\qquad \text{with} \quad
 \widetilde{\D u} (\psi) \:= \,  (U_\psi \ot U_\psi)\, \D u\, (U^{\ast}_\psi \ot U^{\ast}_\psi) .
 \end{equation*}
 To this end we shall use the fact that, as it follows for instance from
  Proposition \ref{multid}, the decomposition $\D u = u_{(1)} \ot u_{(2)}$ is
  equivalent to the Leibniz rule
  \begin{align*}
u (a b) \, = \, u_{(1)}(a) \, u_{(2)}(b) ,  \qquad \fl \, a, b \in
C^\infty (G) .
 \end{align*}
Thus, since
  \begin{align*}
 & \tilde{u} (\psi) (a b)\, = \, U_\psi u \big(U^{\ast}_\psi (a) U^{\ast}_\psi (b)\big) \, = \,
U_\psi \big(u_{(1)}(U^{\ast}_\psi (a)) u_{(2)}(U^{\ast}_\psi (b))\big) \\
& \quad = \, \big(U_\psi u_{(1)}U^{\ast}_\psi \big) (a) \, = \,
\big(U_\psi  u_{(2)} U^{\ast}_\psi \big)(b) \, = \,
\widetilde{u_{(1)}}(\psi) (a) \, \widetilde{u_{(2)}}(\psi)(b) ;
 \end{align*}
evaluating at $\, e \in G$, one obtains
 \begin{equation*}
({\Db} u) (\psi) (a b) \, = \, ({\Db} u_{(1)}) (\psi)(a) \, ({\Db}
u_{(2)}) (\psi)(b)  ,
 \qquad \fl \, a, b \in C^\infty (G) ,
  \end{equation*}
 which is tantamount to \eqref{comcoa3}.
  \endproof

\begin{lemma}\label{mp4lemma}
For  any  $u,v\in U(\Fg)$ one has
\begin{equation}\label{mcoact}
\Db(uv)=u\ps{1}\ns{0}v\ns{0}\ot u\ps{2}\ns{1}(u\ps{2}\rt v\ns{1})
\end{equation}
\end{lemma}
\begin{proof}
 Without loss of generality, we may assume
 $\, u = Z_I $, $\, v = Z_J$. By the definition of the coaction one has
\begin{align*}
&\tilde{u} (\psi)= U_\psi \, u \, U_\psi^\ast=\sum \b_I^K(\psi) Z_K,
\quad \tilde{v} (\psi)= U_\psi \, v \, U_\psi^\ast=\b_J^L(\psi) Z_L,
\end{align*}
which yields
\begin{align*}
&\widetilde{uv} (\psi) =  \tilde{u} (\psi) (\sum_L\b_J^L(\psi) Z_L)=
\sum_L \tilde{u} (\psi)_{(1)} (\b_J^L(\psi) )\,
\tilde{u} (\psi)_{(2)} Z_L  \\
&=  \sum_L \tilde{u}_{(1)} (\psi) (\b_J^L(\psi) )\, \tilde{u}_{(2)}
(\psi) Z_L  ;
\end{align*}
where the last equality follows from  \eqref{comcoa3}. Denoting
$$ \tilde{u}_{(1)} (\psi) = \sum_M \b_{(1)}^M (\psi) Z_M , \quad
 \tilde{u}_{(2)} (\psi) = \sum_N \b_{(2)}^N(\psi) Z_N ,
 $$
one can continue as follows:
\begin{align*}
&\widetilde{uv} (\psi) = \sum_{L, M, N}\b_{(1)}^M (\psi) Z_M
(\b_J^L(\psi) )\, \b_{(2)}^N(\psi) Z_N  Z_L  .
\end{align*}
Evaluating at $e$, one obtains
\begin{align*}
\Db(uv) \, =  \, \sum_{L, M, N}\z_{(1)}^M (\psi) Z_{M \, e}
 (\b_J^L(\psi)) \, \z_{(2)}^N(\psi) Z_{N \, e}  Z_{L \, e} ;
 \end{align*}
taking into account that $U(\Fg)$ is co-commutative, this is
precisely the right hand side of \eqref{mcoact}.
\end{proof}
\medskip

\begin{lemma}\label{mp2lemma}
For any $u\in \Uc(\Fg)$ and any $f\in \Fc(N)$ one has
\begin{equation}\label{newmp4}
\D(u\rt f)=u\ps{1}\ns{0}\rt f\ps{1}\ot u\ps{1}\ns{1}(u\ps{2}\rt
f\ps{2})
\end{equation}
\end{lemma}
\begin{proof}
By Proposition \ref{isoiota} we may assume $f\in \Fc(N)$ of the form
\begin{align*}
f (\psi) = \tilde{f} (\psi) (e) ,
\end{align*}
with $\tilde{f}$ in the algebra generated by $\{\g^i_{j k \ell_1
\ldots \ell_r} \, ; \, 1\leq i, j,k, \ell_1,  \dots,  \ell_r \leq m
\}$. Then
\begin{align} \label{mp2lemma1}
\D(u\rt f) (\psi_1 , \psi_2) =(u\rt f)(\psi_1 \circ \psi_2)= u
(\tilde{f} (\psi_1 \circ \psi_2))(e) .
\end{align}
Now $\tilde{f}$ corresponds to an element $\tilde{\d} \in
\Hc(\Pi)_{\rm ab}$, via the $\Uc(\Fg)$-equivariant isomorphism
$\iota : \Hc(\Pi)^{\rm cop}_{\rm ab} \ra \Fc(N)$; explicitly,
\begin{align*}
\tilde{\d} (g \, U^*_\psi) \, =\,  \tilde{f} (\psi) \,g \, U^*_\psi
.
\end{align*}
Accordingly,
\begin{align*}
&\tilde{f} (\psi_1 \circ \psi_2) \,U^*_{\psi_2} U^*_{\psi_1} \, = \,
\tilde{\d} (U^*_{\psi_2} U^*_{\psi_1}) \, = \,
\tilde{\d}_{(1)} (U^*_{\psi_2}) \, \tilde{\d}_{(2)} (U^*_{\psi_1}) \, =\\
&\tilde{f}_{(1)} (\psi_2) U^*_{\psi_2} \,\tilde{f}_{(2)} (\psi_1)
U^*_{\psi_1} \, =\, \tilde{f}_{(1)} (\psi_2)  \,
\left(\tilde{f}_{(2)} (\psi_1)\circ \psi_2 \right)\,  U^*_{\psi_2}
U^*_{\psi_1} ,
\end{align*}
whence
\begin{align*}
\tilde{f} (\psi_1 \circ \psi_2)  \, = \, \tilde{f}_{(1)} (\psi_2)
\,  \left(\tilde{f}_{(2)} (\psi_1)\circ \psi_2 \right).
\end{align*}
Thus, we can continue \eqref{mp2lemma1} as follows
\begin{align*}
&\D(u\rt f) (\psi_1 , \psi_2) \, =\,
  u \left( \tilde{f}_{(1)} (\psi_2)  \,  (\tilde{f}_{(2)} (\psi_1)\circ \psi_2) \right)(e)= \\
&u_{(1)} ( \tilde{f}_{(1)} (\psi_2))(e)  \,  u_{(2)}\left(\tilde{f}_{(2)} (\psi_1)\circ \psi_2) \right)(e)= \\
&u_{(1)} ( \tilde{f}_{(1)} (\psi_2))(e)  \,
u_{(2)}\left(\tilde{f}_{(2)} (\psi_1)\circ \psi_2) \right)
(\psi_2^{-1}(e))= \\
&u_{(1)} ( \tilde{f}_{(1)} (\psi_2))(e)  \, \left( U_{\psi_2}
u_{(2)} U^*_{\psi_2} \right) (\tilde{f}_{(2)} (\psi_1))(e) .
 \end{align*}
Since $\iota $ switches the antipode with its opposite, the last
line is equal to
\begin{align*}
 (u_{(1)} \rt f_{(2)}) (\psi_2) \,u\ps{2}\ns{1}(\psi_2) \, (u\ps{2}\ns{0} \rt f_{(1)}) (\psi_1) .
\end{align*}
Remembering that $U(\Fg)$ is co-commutative, one finally obtains
\begin{align*}
&\D(u\rt f) (\psi_1 , \psi_2) \, =\,
 (u_{(2)} \rt f_{(2)}) (\psi_2) \,u\ps{1}\ns{1}(\psi_2) \, (u\ps{1}\ns{0} \rt f_{(1)}) (\psi_1) .
 \end{align*}
\end{proof}
\bigskip

\begin{proposition} \label{match}
The Hopf algebras $\Uc:=U(\Fg)$ and $\Fc:=F(N)$ form a matched pair
of Hopf algebras.
\end{proposition}
\begin{proof}
Proposition \ref{moda} together with  Proposition \ref{comcoa} show
that with the action and coaction d defined in \eqref{u>f} and
\eqref{comod3} $\Fc$ is $\Uc$ module algebra and $\Uc$ is a comodule
coalgebra. In addition we shall show that the action and coaction
satisfy \eqref{mp1}\dots \eqref{mp5}. Since $\Uc$ is cocommutative
and $\Fc$ is  commutative \eqref{mp5} is automatically satisfied.
The conditions \eqref{mp2} and \eqref{mp4} are correspondingly
proved in Lemma \ref{mp2lemma} and Lemma \ref{mp4lemma}. Finally the
conditions \eqref{mp3} and \eqref{mp1} are obviously held.
\end{proof}
Now it is the time for the main result of this section.
\begin{theorem}\label{bicrossed}
The Hopf algebras $\Hc({\Pi})^{\rm cop}$ and $\Fc\acl \Uc$ are
isomorphic.
\end{theorem}
\begin{proof}
 Proposition \ref{free}  provides us with $\d_KZ_I$ as a basis for the Hopf algebra
$\Hc:=\Hc(\Pi)^{cop}$. Let us define
$$\Ic:\Hc\ra \Fc\acl \Uc,$$ by  $\Ic(\d_KZ_I)=\iota(\d_K)\acl Z_I$, where
$\iota$ is defined in Proposition \ref{isoiota}, and linearly extend
it on $\Hc$. First let see why $\Ic$ is well-defined. It suffices to
show that the $\Ic$ preserves the relations between elements of
$\Uc$ and $\Fc$. Let $X\in \Fg$ and $f\in \Fc$, by using Proposition
\ref{isoiota}, we have
\begin{align*}
&\Ic(Xf-fX)= \iota(Xf-fX)\acl 1=\\
& X\rt f\acl 1=(\iota(f)\acl
 X)-(1\acl X)(\iota(f)\acl 1 )=\Ic(X)\Ic(f)-\Ic(f)\Ic(X).
\end{align*}

Now we show  $\Ic$ is injective.

 This can be shown by induction on the height. In the height $0$ case the
statement is obvious because because of \eqref{aeta1} and the fact
that $\a^i_{jk}\ot Z_I$ is part of the basis of $\Fc\acl \Uc$.
 Next, assume
 \begin{equation*}
 \sum_{|J| \leq N-1} \, c_{J,I}\,  \eta_J \ot Z_I+  \sum_{|K| = N} \, c_{K,L} \,  \eta_K\ot Z_L = 0.
  \end{equation*}
Using the identities \eqref{aeta3} and \eqref{aeta3}, one can
replace $ \eta_K$ by $ \a_K +$ {\em lower height}. Since the
$\a^{\bullet}_{\bullet \ldots \bullet}$'s are free generators, it
follows that  $c_{K,L} =0$ for each $K$ of height $N$, and thus we
are reduced to
 \begin{equation*}
 \sum_{|J| \leq N-1} \, c_{J,I}\,  \eta_J\ot Z_I  = 0 ;
  \end{equation*}
the induction hypothesis now implies $\, c_{J,I} =0$, for all
$J,I$'s.

So $\Hc$ and $\Fc\acl \Uc$ are isomorphic as algebras. We now show
they are isomorphic as coalgebras as well. It is enough to show
$\Ic$ commutes with coproducts.
\begin{align*}
&\D_{\Fc\acl \Uc}(\Ic(u))=\D_{\Fc\acl \Uc}(1\acl u)= 1\acl
u\ps{1}\ns{0}\ot u\ps{1}\ns{1}\acl u\ps{2}.
\end{align*}
On the other hand let $u\ps{1}\ot U_\vp u\ps{2}U^\ast_\vp=u\ps{1}\ot
\sum \b^I(\vp)Z_I$. We have
\begin{align*}
&u(fU^\ast_\vp g U^\ast_\psi)=\\
&u(fg\circ \td\vp)U^\ast_\vp U^\ast_\psi=\\
&u\ps{1}(f)u\ps{2}(g\circ \td\vp)U^\ast_\vp U^\ast_\psi=\\
&u\ps{1}(f)U^\ast_\vp U_\vp u\ps{2}U^\ast_\vp(g)U^\ast_\psi=\\
&u\ps{1}(f)U^\ast_\vp\b^I(\vp)Z_I(g)U^\ast_\psi,
\end{align*}
which shows that
$\D_{\Hc^{cop}}(u)=u\ps{1}\iota^{-1}(u\ps{2}\ns{1})\ot
u\ps{2}\ns{0}$. Since $U(\Fg)$ is cocommutative one has
$(\Ic\ot\Ic)\D_\Hc(u)=\D_{\Fc\acl \Uc}(\Ic(u))$.
\end{proof}

\subsection{The non-flat case}

We now take up the case of the contact pseudogroup $\Pi_{\rm cn}$, in which case
the Kac decomposition is given by Lemma \ref{HKac}.

As in the flat case, we define the
{\em coordinates} of an element $\,  \psi \in N_{\rm cn}$ as being
 the coefficients of the Taylor expansion of $\psi$ at $0 \in \Rb^{2n+1}$,
  \begin{equation*}
\a^i_{j j_1j_2\dots j_r}(\psi)= \p_{j_r}\dots \p_{j_1} \p_j
\psi^i(x)\mid_{x=0} , \qquad 0 \leq i, j, j_1, j_2, \dots,  j_r \leq 2n .
\end{equation*}
The algebra they generate will be denoted $\, \Fc (N_{\rm cn}) $. It is the
free commutative algebra generated by the indeterminates $\{\a^i_{j
j_1j_2\dots j_r} ; \, 0 \leq i, j, j_1, j_2, \dots,  j_r \leq 2n , \, r \in \Rb \}$, that are
symmetric in all lower indices.
\medskip

\begin{proposition} \label{Fcnhopf}
$\Fc (N_{\rm cn}) $ is a Hopf algebra, whose
coproduct, antipode and counit are uniquely determined by the requirements
\begin{eqnarray} \label{Fcncop}
\D(f)(\psi_1, \psi_2) &=& f(\psi_1 \circ\psi_2) , \qquad \fl \,
\psi_1, \psi_2 \in N_{\rm cn} , \\ \notag S(f)(\psi)  &=& f(\psi^{-1}) ,
\qquad \fl \, \psi \in N_{\rm cn}, \\ \notag \e (f)
&=& f(e) , \quad \fl \, f \in \Fc  (N_{\rm cn}).
\end{eqnarray}
\end{proposition}

\proof The proof is almost identical to that of Proposition \ref{Fhopf}. There
 are $2n$ new coordinates in this case, namely
$\a^i_0$,  $i =1, \dots, 2n$. for which one checks that the coproduct
is well-defined as follows:
 \begin{align*}
&\D \a^i_0 (\psi_1,\psi_2)\, = \, \a^i_{0}(\psi_1\circ \psi_2)=
 \p_0 (\psi_1\circ\psi_2)^i(0)  =\\
&= \p_0\psi^i_1(\psi_2(0)) \p_0\psi^0_2(0) \, + \,
\sum_{j=1}^{2n} \p_j \psi^i_1(\psi_2(0)) \p_0\psi^j_2(0)= \\
&= \p_0\psi^i_1(0) \p_0\psi^0_2(0) \, + \,
\sum_{j=1}^{2n} \p_j \psi^i_1(0) \p_0\psi^j_2(0) =\\
 &= \p_0\psi^i_1(0) \, + \,
\sum_{j=1}^{2n} \d_j^i \p_0\psi^j_2(0)\,  = \,\p_0\psi^i_1(0)
\, + \, \p_0\psi^i_2(0)\,  = \\
&\qquad =  \, (\a^i_0\ot 1+1\ot
 \a^i_0)(\psi_1,\psi_2). ;
\end{align*}
we have been using above the fact that, for any $\, \psi \in N_{\rm cn}$,
\begin{align*}
 \psi (0) \, = \,0 , \qquad \text{and} \qquad
 \psi^H_\ast (0) \, = \,\Id .
 \end{align*}
Taking $\psi_1 = \psi^{-1}$, $\psi_2 = \psi$, one obtains from the above
 \begin{align*}
 \a^i_0 (\psi^{-1}) \, + \, \a^i_0 (\psi) \, = \,0 , \qquad \text{hence} \qquad
 S  \a^i_0  \, = \,- \a^i_0  .
 \end{align*}
\endproof

\begin{lemma} \label{cnregular}
 The coefficients of the Taylor expansion of $\tpsi$ at $e \in N_{\rm cn}$,
\begin{equation} \label{cnfcoord}
 \eta^i_{j k_1 \ldots k_r} (\psi): =  \G^i_{j k_1 \ldots k_r} (\psi)(e) , \qquad
 \psi \in N_{\rm cn} ,
 \end{equation}
define regular functions on $N_{\rm cn}$, which generate the algebra $\, \Fc
(N_{\rm cn})$.
\end{lemma}

\proof The formula \eqref{a0jpsi} shows that
 \begin{align*}
  \eta^i_0  \, = \, - \a^i_0 , \qquad i=1, \ldots , 2n ,
 \end{align*}
 while
  \begin{align*}
  \eta^i_j   \, =  \,\d^i_j \, =  \, \a^i_j , \qquad i , j =1, \ldots , 2n .
 \end{align*}
To relate their higher derivatives, we observe that, in view of \eqref{X0}, \eqref{Xj}
 \begin{align*}
   X_{k_1} \ldots X_{k_r} \mid_{(0, \Id, 1)} \, = \, E_{k_1} \ldots E_{k_r} \mid_{(0, \Id, 1)} ;
 \end{align*}
on the other hand, it is obvious that the jet at $0$ with respect to the frame
$\{ E_0, \ldots , E_{2n} \}$ is equivalent to the jet at $0$ with respect to the
standard frame $\{ \p_0, \ldots , \p_{2n} \}$. This proves the statement for the
`new' coordinates. For the other coordinates the proof is similar to that of
Lemma \ref{regular}.
\endproof

This lemma allows to recover the analog of Proposition \ref{isoiota} by
identical arguments.

\begin{proposition} \label{cnisoiota}
There is a unique isomorphism of Hopf algebras \\
$ \iota : \Hc(\Pi_{\rm cn})^{\rm cop}_{\rm ab} \ra \Fc (N_{\rm cn})$ with the property
that
\begin{equation*}
\iota (\D^i_{j k_1 \ldots k_r}) \, = \, \eta^i_{j k_1
\ldots k_r} , \qquad  0 \leq i, j, k_1, k_2, \dots,  k_r \leq 2n .
\end{equation*}
\end{proposition}

 Next, one has the tautological counterpart of Lemma \ref{ltact}.

  \begin{lemma} \label{cnltact}
   Let $\vp \in G_{\rm cn}$ and $\phi \in \Gb_{\rm cn}$. Then
   \begin{align*}
   &\G^i_{j  k_1 \ldots k_r} (\vp \circ \phi)\, = \, \G^i_{j  k_1 \ldots k_r} ( \phi) ,\\
   &\G^i_{j  k_1 \ldots k_r} (\phi \circ \vp)\, = \,
   \G^i_{j  k_1 \ldots k_r} ( \phi) \circ \tilde{\vp}  .
   \end{align*}
  \end{lemma}

 \proof As is the case with its sibling result, this is simply a consequence of
 the left invariance the vector fields $\{ X_0, \ldots , X_{2n} \}$.
 \endproof

 In turn, the above lemma allows to recover the analog of Proposition \ref{moda}.

  \begin{proposition}\label{cnmoda}
  The algebra isomorphism $ \iota : \Hc(\Pi_{\rm cn})^{\rm cop}_{\rm ab} \ra \Fc (N_{\rm cn})$
  identifies the  $\,\Uc(\Fg_{\rm cn})$-module  $\Hc(\Pi_{\rm cn})_{\rm ab}$
  with the $\,\Uc(\Fg_{\rm cn})$-module $\Fc (N_{\rm cn})$. In particular $\, \Fc (N_{\rm cn})$ is
  $\,\Uc(\Fg_{\rm cn})$-module algebra.
  \end{proposition}

  Furthermore,
  $\Uc(\Fg_{\rm cn})$ can be endowed with a right $\Fc (N_{\rm cn})$-comodule structure
   $\,\Db:\Uc(\Fg_{\rm cn}) \ra\Uc(\Fg_{\rm cn}) \ot \Fc (N_{\rm cn})$
  in exactly the same way as in the flat case, \cf \eqref{comod3}, Lemma \ref{cms}, and
  is in fact a  right $\Fc(N_{\rm cn})$-comodule coalgebra ({\it comp.} Prop. \ref{comcoa}).
  Likewise, the analog of Proposition \ref{match} holds true, establishing that
   $\Uc(\Fg_{\rm cn})$ and $\Fc(N_{\rm cn})$   form a matched pair
of Hopf algebras. Finally, one concludes in a similar fashion with the
bicrossed product realization theorem for the contact case.

\begin{theorem}\label{cnbicrossed}
The Hopf algebras $\Hc(\Pi_{\rm cn})^{\rm cop}$ and
$\Fc(N_{\rm cn})\acl \Uc(\Fg_{\rm cn})$ are canonically
isomorphic.
\end{theorem}



\section{Hopf cyclic cohomology}

After reviewing some of the most basic notions in Hopf cyclic
cohomology, we focus on the case of the Hopf algebras $\Hc (\Pi)$
constructed in the preceding section and show how their Hopf  cyclic
cohomology can be recovered from a bicocyclic complex manufactured
out of the matched pair. We then illustrate this procedure by
computing the relative periodic Hopf cyclic cohomology of $\Hc_n$
modulo $\Fg \Fl_n$. For $n=1$, we completely calculate the
non-periodized Hopf cyclic cohomology as well.

\subsection{Quick synopsis of Hopf cyclic cohomology}
Let $\Hc$ be a Hopf
algebra,  and let $C$ be a left $\Hc$-module coalgebra,
  such that its comultiplication and counit are $\Hc$-linear, \ie
\begin{align*}
\D(hc)=h\ps{1}c\ps{1}\ot h\ps{2}c\ps{2}, \quad \ve(hc)=\ve(h)\ve(c).
\end{align*}
We recall from \cite{hkrs1}
 that a right module $M$ which is also  a left comodule is called
{\em right-left stable anti-Yetter-Drinfeld module} (SAYD for short)
over the Hopf algebra  $\Hc$ if it satisfies
the following conditions, for any $h\in \Hc$, and $m\in M$ :
\begin{align*}
&m\sns{0}m\sns{-1}=m\\\label{SAYD2}
 &(mh)\sns{-1}\ot (mh)\sns{0}=
S(h\ps{3})m\sns{-1}h\ps{1}\ot m\sns{0}h\ps{2},
\end{align*}
where the coaction of $\Hc$ was denoted by $\Db_M(m)=m\sns{-1}\ot
m\sns{0}$.

\medskip
Having such a datum $(\Hc, C, M)$, one defines (\cf~\cite{hkrs2}) a cocyclic
module $\, \{C^n_\Hc(C,M), \p_i,\s_j,\tau \}_{n\ge 0}$ as follows.
 \begin{equation*}
 \Cc^n:= C^n_\Hc(C,M)=M\ot_\Hc C^{\ot n+1}, \quad n\ge 0,
 \end{equation*}
 with the cocyclic structure given by the operators
\begin{align*}
&\p_i:\Cc^n\ra \Cc^{n+1}, & 0\le i\le n+1\\
& \s_j: \Cc^n\ra \Cc^{n-1}, & 0\le j\le n-1,\\
& \tau:\Cc^n\ra \Cc^n,
\end{align*}
defined explicitly as follows:  \begin{align*} & \p_i(m\ot_\Hc \td c
)=m\ot_\Hc c^0\odots \Delta(c_i)\odots c^n,
\\
 &\p_{n+1}(m\ot_\Hc \td c)=m\sns{0}\ot_\Hc c^0\ps{2}\ot c^1\odots c^n\ot
 m\sns{-1}c^0\ps{1},\\
 &\sigma_i(m\ot_\Hc \td c)=m\ot_\Hc c^0\odots \epsilon(c^{i+1})\odots c^n, \\
 &\tau(m\ot_\Hc \td c)=m\sns{0}\ot_\Hc c^1\odots c^n\ot m\sns{-1}c^0 ;
 \end{align*}
 here we have used the abbreviation $\td c=c^0\odots c^n$.

One checks~\cite{hkrs2} that $ \p_i,\s_j,$ and $\tau$ satisfy the
following identities, which define the structure of a
cocyclic module (\cf~\cite{cng}):
\begin{equation}\label{ds}
\p_j  \p_i = \p_i  \p_{j-1}, \, \, i < j  , \qquad \s_j \s_i = \s_i
\s_{j+1},  \, \,  i \leq j
\end{equation}
\begin{equation}
\s_j  \p_i = \left\{ \begin{matrix} \p_i  \s_{j-1} \hfill &i < j
\hfill \cr 1_n \hfill &\hbox{if} \ i=j \ \hbox{or} \ i = j+1 \cr
\p_{i-1}  \s_j \hfill &i > j+1  ;  \hfill \cr
\end{matrix} \right.
\end{equation}
\begin{eqnarray}
\tau_n  \p_i  = \p_{i-1}  \tau_{n-1} ,
 \quad && 1 \leq i \leq n ,  \quad \tau_n  \p_0 =
\p_n \\ \label{cj} \tau_n  \s_i = \s_{i-1} \tau_{n+1} , \quad &&
1 \leq i \leq n , \quad \tau_n  \s_0 = \s_n  \tau_{n+1}^2 \\
\label{ce} \tau_n^{n+1} &=& 1_n  \, .
\end{eqnarray}

The motivating example for the above notion is the
cocyclic complex associated to a Hopf algebra $\Hc$ endowed with a {\em modular pair
in involution}, (MPI for short), $(\d,\s)$, which we recall from
\cite{cm2}.   $\d$ is an algebra map $\Hc\ra \Cb$, and $\s\in \Hc$ is
a group-like element,  or equivalently a coalgebra map $\Cb\ra \Hc$.
 The pair $(\d, \s)$ is
called MPI if $\d(\s)=1$, and  $\td{ S}^2_\d=Ad \s$; the twisted
antipode $\td S_\d$ is defined by
\begin{equation*}
\td{ S_\d}(h)=(\d\ast S)(h)=\d(h\ps{1})S(h\ps{2}).
\end{equation*}

One views  $\Hc$ as a left $\Hc$-module coalgebra via left
multiplication. On the other hand if one lets $M =~^\s\Cb_\d$ to be
the ground field $\Cb$ endowed with the left $\Hc$-coaction  via
$\s$ and right $\Hc$-action via the character $\d$, then   $(\d,\s)$
is a MPI if and only if $^\s\Cb_\d$ is a SAYD. Thanks to the
multiplication and the antipode of $\Hc$, one identifies $C_\Hc(\Hc,
M)$ with $M\ot \Hc^{\ot n}$ via the map
\begin{align*}
&\Ic: M\ot_\Hc \Hc^{\ot (n+1)}\ra M\ot \Hc^{\ot n},\\
&\Ic(m\ot_\Hc h^0\odots h^n)=mh^0\ps{1}\ot S(h\ps{2})\cdot(h^1\odots
h^n).
\end{align*}
As a result, $\p_i$, $\s_j$, and $\tau$ acquire the simplified form of
 the original definition~\cite{cm2}, namely
\begin{eqnarray*}
\p_0 (h^1 \ot \ldots \ot h^{n-1}) &=& 1 \ot h^1 \ot \ldots \ot
h^{n-1} , \\  \nonumber \p_j (h^1 \ot \ldots \ot h^{n-1}) &=& h^1
\ot \ldots \ot \D h^j \ot \ldots \ot h^{n-1} , \qquad 1 \leq j \leq
n-1
\\  \nonumber
\p_n (h^1 \ot \ldots \ot h^{n-1}) &=& h^1 \ot \ldots \ot h^{n-1} \ot
\s, \\ \nonumber
 \s_i (h^1 \ot \ldots \ot h^{n+1}) &=& h^1 \ot \ldots \ot \ve
(h^{i+1}) \ot \ldots \ot h^{n+1} , \qquad 0 \leq i \leq n \, , \\
\nonumber
   \tau_n (h^1 \ot \ldots \ot h^n) &=& (\D^{p-1}  \wt S (h^1)) \cdot h^2
\ot \ldots \ot h^n \ot \s .
\end{eqnarray*}

For completeness, we record below the bi-complex
$\, (C C^{*, *} (C, \Hc,  M), \, b , \, B )$
that computes the Hopf cyclic cohomology of a coalgebra $C$ with
coefficients in a SAYD module $M$ under the symmetry of a Hopf
algebra $\Hc$:
\begin{equation*} \label{CbB}
CC^{p, q} (C, \Hc; M)  \, = \, \left\{
\begin{matrix}  C^{q-p}_\Hc (C, M ) \, , \quad q \geq p \, , \cr
  0 \, ,  \quad \qquad  \qquad q <  p \, , \end{matrix}   \right.
\end{equation*}
where    $\, b: C^{n}_\Hc (C,M) \ra C^{n+1} _\Hc(C,M)$ is given by
\begin{equation*}  \nonumber
 b =
\sum_{i=0}^{n+1} (-1)^i \p_i \, ;
\end{equation*}
the operator $B:C^{n}_\Hc (C,M) \ra C^{n-1}_\Hc (C,M)$ is
defined by the formula
\begin{equation} \nonumber
  B = A \circ B_{0} \, ,  \quad n \geq 0 \, ,
\end{equation}
where
\begin{equation*}
B_{0}= \s_{n-1}\tau(1-(-1)^n\tau)
\end{equation*}
 and
\begin{equation} \nonumber
A = 1 + \lb  + \cdots +  \lb^n \,  , \qquad \text{with} \qquad \lb =
(-1)^{n-1} \tau_n \, .
\end{equation}

The groups $\, \{H C^n  (\Hc; \delta,\sigma)\}_{n \in \Nb} \,$ are
computed from the first quadrant total complex $\, ( TC^{*} (\Hc;
\delta,\sigma),\, b+B ) \,$,
\begin{equation*}
    TC^{n}(\Hc; \delta,\sigma) \,  = \, \sum_{k=0}^{n} \, C C^{k, n-k} (\Hc; \delta,\sigma) \, ,
\end{equation*}
and the periodic groups $\, \{H P^i  (\Hc; \delta,\sigma)\}_{i \in
\Zb/2} \,$ are computed from the full total complex $\,(T P^{*}(\Hc;
\delta,\sigma),\,  b+B ) \,$,
\begin{equation*}
      T P^{i}(\Hc; \delta,\sigma) \,  = \, \sum_{k \in \Zb} \, C C^{k, i-k} (\Hc; \delta,\sigma) \, .
\end{equation*}
\bigskip

We note that, in defining the Hopf cyclic
cohomology as above, one has the option of viewing
the Hopf algebra $\Hc$ as a left
$\Hc$-module coalgebra or as
 a right $\Hc$-module coalgebra.  It
was the first one which  was
 selected as the definition in \cite{cm2}.
 The other choice would have given to the cyclic operator the expression
\begin{equation*}
\tau_n(h^1\odots h^n)= \s S(h^n\ps{n})\ot h^1 S(h^n\ps{n-1})\odots
h^{n-1}S(h^n\ps{1})\d(S^{-1}(h^n\ps{2})).
\end{equation*}
As it happens, the choice originally selected
 is not the best suited for the
 situations involving a right action. To restore the naturality of the notation,
 it is then convenient to pass from the Hopf algebra $\Hc$
 to the  co-opposite Hopf algebra $\Hc_n\cop$. This transition does not
 affect the Hopf cyclic
   cohomology, because
         for any $\Hc$-module coalgebra $C$  and any SAYD
         module $M$ one has a canonical equivalence
 \begin{equation}\label{equiv}
 (C^\ast_\Hc( C,M), b, B)\simeq  (C^\ast_{\Hc^{\rm cop}}(C^{\rm cop},  M\cop), b, B) ;
 \end{equation}
 $ M\cop:= M$  is SAYD module for $\Hc\cop$,
with the action of $\Hc\cop$  the same as the action of $\Hc$, but with the coaction
$\,  \Db:M\cop \ra \Hc\cop\ot M\cop$ given by
 \begin{equation}
 \Db(m) \, = \, S^{-1}(m\sns{-1})\ot m\sns{0}.
 \end{equation}
 The equivalence \eqref{equiv} is realized by the map
\begin{align}
 &\Tc : C^n_\Hc(C,M)\ra  C^n_{\Hc\cop}(C\cop,M\cop) , \\\notag
 &\Tc(m\ot c^0\odots c^n)=m\sns{0}\ot   m\ns{-1}c^0\ot c^{n} \odots c^1 .
 \end{align}
    \begin{proposition}\label{cop}
  The map $\Tc$  defines an isomorphism of  mixed  complexes.
  \end{proposition}
  \begin{proof}  The map $\Tc$ is well-defined because
  \begin{align*}
   &\Tc(mh\ot c^0\odots c^n)= mh\ps{2}\ot S(h\ps{3})m\sns{-1}h\ps{1}c^0\ot c^n\odots c^1=\\
  & m\sns{0}\ot  h\ps{n+2} S(h\ps{n+3}) m\sns{-1}h\ps{1} c^0\ot h\ps{n+1}c^n\odots h\ps{2}c^1\\
    &m\sns{0}\ot m\sns{-1} h\ps{1}c^0\ot h\ps{n+1} c^n\odots h\ps{2}c^1=\\
    &\Tc(m\ot h\ps{1}c^0\odots h\ps{n+1}c^n).
   \end{align*}
   We denote
  the cyclic structure of $C^n_\Hc(C,M)$, resp. $C_{\Hc\cop}(C\cop,M\cop)$ by $\p_i$, $\s_i$, and
  $\tau$, resp. $d_i$, $s_j$, and $t$.
    We need to show that $\Tc$  commutes with $b$ and $B$. One has in fact
 a stronger  commutation property, namely
  \begin{equation}
  \Tc\p_i=d_{n+1-i}\Tc, \quad 0\le i\le n+1.
  \end{equation}
  Indeed,
    \begin{align*}
   & d_{n+1} \Tc(m\ot c^0\odots c^n)=\\
   &d_{n+1} (m\sns{0}\ot   m\ns{-1}c^0\ot c^{n} \odots c^1)=\\
   &m\sns{0} \ot m\sns{-3} c^0\ps{1}\ot c^n\odots c^1\ot       S^{-1}(m\sns{-1})m\sns{-2}c^0\ps{2}   \\
   &m\sns{0}\ot m\sns{-1} c^0\ps{1}\ot c^n\odots c^1\ot c^0\ps{2}\\
              &\Tc (m\ot c^0\ps{1}\ot c^0\ps{2}\ot c^1\odots c^n)=\\
       &\Tc\p_0(m\ot c^0\odots c^n)
      \end{align*}
   \begin{align*}
  & \Tc\p_{n+1}(m\ot c^0\odots c^n)=\\
   &\Tc(m\sns{0}\ot c^0\ps{2}\ot c^1\odots c^n\ot m\sns{-1}c^0\ps{1})\\
  & m\sns{0}\ot m\sns{-1}c^0\ps{2} \ot m\sns{-2}c^0\ps{1}\ot c^n\odots c^1\\
  &d_0(m\sns{0}\ot m\ns{-1}c^0\ot c^n\odots c^1)=\\
  &d_0\Tc(m\ot c^0\odots c^n) ;
      \end{align*}
     on the other hand, for $1\le i\le n$,
      \begin{align*}
      &\Tc\p_i(m\ot c^0\odots c^n)=\\
      &\Tc(m\ot c^0\odots c^{i-1}\ot  c^i\ps{1}\ot c^i\ps{2}\ot c^{i+1}\odots c^n)=\\
      &m\sns{0}\ot m\sns{-1}c^0 \ot c^n\odots c^{i+1}\ot c^i\ps{2}\ot c^i\ps{1}\ot c^{i-1}\odots c^1=\\
      &d_{n+1-i}( m\sns{0}\ot m\sns{-1}c^0\ot c^n\odots c^1)=\\
      &d_{n+1-i}\Tc(m\ot c^0\dots c^n) .
      \end{align*}
     Thus,  $\Tc b=(-1)^{n+1}b\Tc$.

      Next,  we check that $\Tc\tau=t^{-1}\Tc=t^n\Tc$ as follows:
       \begin{align*}
       &\Tc\tau(m\ot c^0\odots c^n)=\\
       &\Tc(m\sns{0}\ot c^1\odots c^n\ot m\sns{-1}c^0)=\\
       &m\sns{0}\ot m\sns{-1}c^1\ot m\sns{-2}c^0\ot c^n\odots c^2=\\
       &t^{-1}(m\sns{0}\ot m\sns{-1}c^0\ot c^n\odots c^1)=\\
       &t^{-1}\Tc(m\ot c^0\odots c^n) .
       \end{align*}
       It is easy to see that $\Tc\s_i=s_{n-1-i}\Tc$, for $0\le i\le n-1$.
       Using the above identities  and  $ts_0=s_{n-1}t^2$ one obtains
           \begin{align*}
       &\Tc B=\Tc \sum_{j=0}^{n-1}(-1)^{(n-1)j}\tau^j \s_{n-1}\tau(1-(-1)^n\tau)=\\
       &=\sum_{j=0}^{n-1}(-1)^{(n-)j}t^{-j}s_0t^{-1}(1-(-1)^nt^{-1})\Tc=\\
       &=\sum_{j=0}^{n-1}(-1)^{(n-1)j}t^{n-j} t^{-1}s_{n-1}t^2   (1-(-1)^nt^{-1})    \Tc=\\
      &=(-1)^{n+1}\sum_{j=0}^{n-1}(-1)^{(n-1)j}t^{n-1-j} s_{n-1}t  (1- (-1)^nt)    \Tc=\\
       &=(-1)\sum_{k=0}^{n-1}(-1)^{(n-1)k}t^{k}s_{n-1}t(1-(-1)^nt)\Tc= \quad -B\Tc ,
       \end{align*}
       which completes the proof.
           \end{proof}
\bigskip

We next recall from \cite{cm6} the setting for relative Hopf cyclic cohomology.
Let $\Hc$ be an arbitrary Hopf algebra and $\Kc \sbs \Hc$ a Hopf
subalgebra.  Let
\begin{equation}
    \Cc \, = \, \Cc (\Hc, \Kc) \, := \, \Hc \ot_{\Kc} \Cb   \, ,
\end{equation}
where $\Kc$ acts on $\Hc$ by right multiplication and on $\Cb$ by
the counit. It is a left $\Hc$-module in the usual way, via
left multiplication. As such, it  can be identified with the
quotient module $\, \Hc/\Hc \Kc^+ $, $\, \Kc^+ = \Ker \ve | \Kc \,$, via the
isomorphism  induced by
\begin{equation}
h \in \Hc \, \longmapsto \, \dot{h} = h \ot_{\Kc} 1 \in \Hc \ot_{\Kc} \Cb \, .
\end{equation}
Moreover, thanks to the right action of $\Kc$ on $\Hc$,  $\, \Cc =
\Cc (\Hc, \Kc) \,$ is an $\Hc$-module coalgebra. Indeed, its
coalgebra structure is given by the coproduct
 \begin{equation} \label{copr}
     \D_{\Cc} \, (h \ot_{\Kc} 1) \, = \, (h\ps{1} \ot_{\Kc} 1) \ot (h\ps{2} \ot_{\Kc} 1) \, ,
\end{equation}
inherited from that on $\Hc$,  and is compatible with the action of
$\Hc$ on $\Cc$ by left multiplication:
\begin{equation*}
 \D_{\Cc} \, (gh \ot_{\Kc} 1) \, = \,  \D(g) \D_{\Cc} \, (h \ot_{\Kc} 1) \, ;
 \end{equation*}
similarly, there is an inherited counit
\begin{equation} \label{coun2}
     \ve_{\Cc} \, (h \ot_{\Kc} 1) \, = \,  \ve (h) \, , \quad \fl \, c \in C \, ,
\end{equation}
that satisfies
\begin{equation*}
  \ve_{\Cc} \, (gh \ot_{\Kc} 1) \, = \,  \ve (g) \, \ve_{\Cc} \, (h \ot_{\Kc} 1) \,  .
 \end{equation*}

 Thus, $\Cc$ is a
$\Hc$-module coalgebra. The relative Hopf
cyclic cohomology of  $\Hc$ with respect to the  $\Kc$ and with
coefficients in $M$, to be denoted by $HC(\Hc,\Kc;M)$, is by definition
the Hopf cyclic cohomology
of  $\Cc$ with coefficients  in $M$.
Thanks to the antipode of $\Hc$ one simplifies the cyclic complex
as follows (\cf \cite[\S 5]{cm6}):
\begin{equation*}
C^{\ast} (\Hc, \Kc ; M) = \{ C^n  (\Hc, \Kc ; M) : = \, M
\ot_{\Kc} \Cc^{\ot n} \}_{n \geq 0} ,
\end{equation*}
where $\Kc$ acts diagonally on $ \Cc^{\ot n}$,
 \begin{align*}
 &\p_0 (m \ot_{\Kc} c^1 \ot \ldots \ot c^{n-1}) =
 m \ot_{\Kc} \dot{1} \ot c^1 \ot \ldots \ot
  \ldots \ot c^{n-1} ,   \\  \nonumber
&\p_i  (m \ot_{\Kc} c^1 \ot \ldots \ot c^{n-1})  = m \ot_{\Kc} c^1
\ot \ldots \ot
 c^i_{(1)} \ot c^{i}_{(2)} \ot
\ldots \ot c^{n-1} ,   \\  \nonumber &  \fl \quad 1 \leq i \leq n-1
\, ; \\  \nonumber
&\p_n  (m \ot_{\Kc} c^1 \ot \ldots \ot c^{n-1}) =
m_{(0)} \ot_{\Kc}  c^1 \ot \ldots \ot c^{n-1} \ot \dot{m}_{(-1)}\, ;
\\ \nonumber & \s_i  (m \ot_{\Kc} c^1 \ot \ldots \ot c^{n+1})  =
 m \ot_{\Kc}  c^1 \ot \ldots \ot \ve
(c^{i+1}) \ot \ldots \ot c^{n+1} ,  \\ \nonumber &
 \fl \quad 0 \leq
i \leq n \, ; \\    \nonumber
 &  \tau_n  (m \ot_{\Kc}  \dot{h}^1 \ot c^2 \ot \ldots \ot c^n) =
   m_{(0)} h^1_{(1)} \ot_{\Kc}  S( h^1_{(2)}) \cdot (c^2
   \ot \ldots \ot c^n \ot \dot{m}_{(-1)});
\end{align*}
the above operators are well-defined and
endow $C^{\ast} (\Hc, \Kc ; M) $ with a cyclic structure.
\smallskip

Since
\begin{equation*}
k_{(1)}\,  h\,  S(k_{(2)} ) \ot_{\Kc} 1 \, = \, k_{(1)}\, h\ot_{\Kc}
\ve (k_{(2)} )\, 1 \, = \, k \, c \ot_{\Kc} 1\, , \quad k \in \Kc ,
\end{equation*}
the restriction to $\Kc$ of the
 left action of $\Hc$ on $\Cc = \Hc\ot_{\Kc} \Cb $ can also be regarded
 as `adjoint action', induced by conjugation.
\bigskip

When specialized to Lie algebras this definition recovers
 the relative Lie algebra homology, \cf \cite[Thm. 16]{cm6}. Indeed,
let $\Fg$ be a Lie algebra over the field $F$, let $\Fh \sbs \Fg$ be
a reductive subalgebra in $\Fg$, and let $M $ be a $\Fg$-module.  We equip $M $
with the \textit{trivial} $\Fg$-comodule structure
\begin{equation} \label{cotriv}
  \Db_M (m) \, = \,  1   \ot m  \, \in \Hc \ot M \,  ,
 \end{equation}
and note the stability  condition  is then trivially satisfied,
while the AYD one  follows from (\ref{cotriv}) and the
cocommutativity of the universal enveloping algebra $\FA (\Fg)$. The
relative Lie algebra homology and cohomology of the pair $ \Fh \sbs
\Fg$ with coefficients in $M$ is computed from the
Chevalley-Eilenberg  complexes
\begin{equation*}
\{C_* (\Fg, \Fh ; M) , \d \} , \quad C_n (\Fg, \Fh ; M) := M
\ot_{\Fh} \bigwedge ^n (\Fg/\Fh).
\end{equation*}
Here the action of $\Fh$ on $\Fg/\Fh$ is induced by the adjoint
representation and the differentials are given by the formulae
\begin{align} \notag
&\d(m\ot_{\Fh}\dot{X_1} \wdg \ldots \wdg \dot{X}_{n+1} ) \, = \,
\sum_{i=1}^{n+1}  (-1)^{i + 1}  m X_i \ot_{\Fh} \dot{X}_1 \wdg
\ldots \wdg \Check{\dot{X}}_i \ldots  \wdg \dot{X}_{n+1}
\\\label{dh} & + \sum_{i < j} (-1)^{i + j} m \ot_{\Fh}
\dot{\widehat{[X_i , X_j ]}} \wdg \dot{X}_1 \wdg \ldots \wdg
\Check{\dot{X}}_i \ldots \Check{\dot{X}}_j \ldots  \wdg
\dot{X}_{n+1}
\end{align}
where $\dot{X} \in \Fg/\Fh$ stands for the class modulo $\Fh$ of $X
\in \Fg$ and the superscript $\Check{}$ signifies the omission of
the indicated variable.

There are canonical isomorphisms between the periodic relative Hopf
cyclic cohomology  of the pair $\FA(\Fh) \subset \FA(\Fg) $, with
coefficients in any $\Fg$-module $M$, and the relative Lie algebra
homology  with coefficients of the pair $\Fh \subset \Fg $ ( see
\cite{cm6}):
\begin{eqnarray} \label{Lie}
HP^\e (\FA(\Fg), \FA(\Fh) ; M) &\cong& \bigoplus_{n  \equiv \e \mod 2}
H_n (\Fg, \Fh ; M) .
\end{eqnarray}


\subsection{Reduction to diagonal mixed complex} \label{secmixed}

 In this section we develop an apparatus for computing Hopf cyclic cohomology
 of certain cocrossed coproduct coalgebras. These are
 made of two coalgebras endowed with actions by a Hopf algebra.
 Both  coalgebras are  Hopf module coalgebras, and one of them
 is a comodule coalgebra as well, which moreover is an
 Yetter-Drinfeld module. Under these circumstances, one can unwind
  the Hopf cyclic structure of the cocrossed product coalgebra
 and identify it with the diagonal of a cylindrical module.
We then construct a spectral sequence that
 computes the cohomology of the total complex of the cylindrical module,
  which by Eilenberg-Zilber
  theorem is quasi isomorphic to the diagonal of the cylindrical
  module.

 \bigskip

Let  $\Hc$ be a Hopf algebra and let $C$ and $D$ be two left
$\Hc$-module coalgebras. In addition, we assume that $C$ is a left
comodule coalgebra, which makes $C$ a YD module on $\Hc$. We denote
the coaction of $C$ by $\Db_C $ and use the abbreviated notation
$\Db_C(c)=c\ns{-1}\ot c\ns{0}$. We recall that the YD condition stipulates that
\begin{align}\label{yd}
\Db_C(hc)=h\ps{1}c\ns{-1}S(h\ps{3})\ot h\ps{2}c\ns{0}.
\end{align}

Using the coaction of $\Hc$ on $C$ and  its action on $D$, one
constructs a coalgebra structure on $C\ot D$ defined by the coproduct
\begin{align}\label{cocrossedcoproduct}
\D(c\ot d)= c\ps{1}\ot c\ps{2}\ns{-1}d\ps{1}\ot c\ps{2}\ns{0}\ot
d\ps{2}.
\end{align}
We denote this coalgebra by $\nobreak{C\cl D}$.  The import of the
YD condition is revealed by the following result.
\begin{lemma}\label{daction}
Via the diagonal action of $\Hc$,  $\nobreak{C\cl D} $ becomes an
$\Hc$ module coalgebra.
\end{lemma}
\begin{proof} One has
\begin{align*}
&\D(h\ps{1}c\cl h\ps{2}d)=\\
& h\ps{1}c\ps{1}\cl (h\ps{2}c\ps{2})\ns{-1}(h\ps{3}d\ps{1})\ot
(h\ps{2}c\ps{2})\ns{0}\cl
h\ps{4}d\ps{2}=\\
&h\ps{1}c\ps{1}\cl h\ps{2}c\ps{2}\ns{-1}S(h\ps{4})h\ps{5}d\ps{1}\ot
h\ps{3}c\ps{2}\ns{0}\cl
h\ps{6}d\ps{2}=\\
&h\ps{1}c\ps{1}\cl h\ps{2}c\ps{2}\ns{-1}d\ps{1}\ot
h\ps{3}c\ps{2}\ns{0}\cl h\ps{4}d\ps{2}=\\
&h\ps{1}(c\ps{1}\cl c\ps{2}\ns{-1}d\ps{1})\ot
h\ps{2}(c\ps{2}\ns{0}\cl d\ps{2})=\\
&h\ps{1}((\nobreak{c\cl d})\ps{1})\ot  h\ps{2}((\nobreak{c\cl
d})\ps{2}).
\end{align*}
\end{proof}

Now let $M$ be an SAYD over $\Hc$.  We endow  $M\ot C^{\ot q}$  with
the following $\Hc$ action and coaction:
\begin{align*}
 &(m\ot\td c)h= mh\ps{1}\ot S(h\ps{2})\td c, \\
 &\Db(m\ot \td c)= c^0\ns{-1}\dots c^n\ns{-1}
 m\sns{-1}\ot m\sns{0}\ot c^0\ns{0}\odots c^{n}\ns{0},
\end{align*}
\begin{lemma}\label{YD}
Let $C$ be an YD module over $\Hc$. Then via the diagonal action and
coaction $C^{\ot q}$ is also an YD module over $\Hc$.
\end{lemma}
\begin{proof}
We verify that $\td c=c^1\odots c^q\in C^{\ot n}$ and $h\in \Hc$
satisfy \eqref{yd}. Indeed,
\begin{align*}
&\Db(h\td c)=\Db(h\ps{1}c^1\odots h\ps{q}c^q)=\\
&(h\ps{1}c^1)\ns{-1}\dots (h\ps{q}c^q)\ns{-1}\ot
(h\ps{1}c^1)\ns{-1}\odots (h\ps{q}c^q)\ns{-1}=\\
&h\ps{1}c^1\ns{-1}S(h\ps{3})h\ps{4}c^2\ns{-1}S(h\ps{6})\dots
h\ps{3q-2}c^q\ns{-1}S(h\ps{3q}) \ot\\
& h\ps{2}c^1\ns{0}\ot h\ps{5}c^2\ns{0}\odots h\ps{3q-1}c^q\ns{0}=\\
&h\ps{1}\td c\ns{-1}S(h\ps{3})\ot h\ps{2}\td c\ns{0}.
\end{align*}
\end{proof}

\begin{proposition}\label{AYD}
Equipped with  the above action and coaction, $M\ot C^{\ot q}$ is
 an AYD module.
\end{proposition}
\begin{proof}
Let $\td c\in C^{\ot q}$, and $h\in \Hc$. By Lemma \ref{YD},
we can write
\begin{align*}
&\Db((m\ot \td c)h)=\Db(m h\ps{1}\ot S(h\ps{2})\td c)= \\
&(S(h\ps{2})\td c)\ns{-1}(mh)\ns{-1}\ot (mh\ps{1})\ns{0}\ot
(S(h\ps{1})\td c)\ns{0}=\\
&S(h\ps{6})\td c\ns{-1}S^2(h\ps{4})S(h\ps{3})m\ns{-1}h\ps{1}\ot
m\ns{0}h\ps{2}\ot S(h\ps{5})\td c\ns{0}=\\
&S(h\ps{3})\td c\ns{-1}m\ns{-1}h\ps{1}\ot (m\ns{0}\ot \td
c\ns{0})h\ps{2}.
\end{align*}
\end{proof}

We define the following bigraded module, inspired by
\cite{ft,gj,ak,mr}, in order to obtain a cylindrical module for cocrossed
product coalgebras. Set
$$\FX^{p,q}:=   M\ot_\Hc D^{\ot p+1}\ot C^{ \ot q+1} ,$$
 and endow $\FX$ with  the operators
  \begin{align}\label{hd}
&\hd_i:\FX^{(p,q)}\rightarrow \FX^{(p+1,q)}, && 0\le i\le p+1\\
&\hs_j: \FX^{(p,q)}\rightarrow \FX^{(p-1,q)}, &&0\le j\le
p-1\\\label{hta}
 &\hta:  \FX^{(p,q)}\rightarrow \FX^{(p,q)},
\end{align}
defined by
\begin{align*}
&\hd_{i}(m\ot \td{d} \ot \td{c})= m\ot d^0\odots
\D(d^i)\odots d^p\ot \td c,\\
&\hd_{p+1}(m\ot \td{d} \ot \td{c})=m\sns{0}\ot d^0\ps{2}\odots
d^p\ot\td c\ns{-1}m\sns{-1}d^0\ps{1} \ot \td c\ns{0},\\
 &\hs_j(m\ot \td{d} \ot \td{c})= m\ot  d^1\odots
\ve(d^{j})\odots d^p\ot  \td{c},\\
 &\hta(m\ot \td{d} \ot \td{c})= m\sns{0}\ot d^1\odots d^p\ot \td c\ns{-1}m\sns{-1}d^0\ot  \td c\ns{0} ;
\end{align*}
the vertical structure is just the cocyclic structure of
$C(\Hc,\Kc;\Kc^{\ot p+1}\ot M)$, with
\begin{align} \label{cvd}
&\vd_i= :\FX^{(p,q)}\rightarrow \FX^{(p,q+1)}, &&0\le i\le q+1\\
\label{cvs} &\vs_j: \FX^{(p,q)}\rightarrow \FX^{(p,q-1)}, && 0\le
j\le q-1\\ \label{cvt} &\vta: \FX^{(p,q)}\rightarrow \FX^{(p,q)},
\end{align}
defined by
\begin{align*}
 &\vd_{i}(m\ot \td{d}\ot\td{c})=m\ot \td d\ot c^0\odots \D(c^i)\odots c^q,\\
&\vd_{q+1}(m\ot \td{d}\ot  \td{c})= \\
&=m\sns{0}\ot S^{-1}(c^0\ps{1}\ns{-1}) \td d\ot c^0\ps{2}\ot  c^1\odots c^q\ot  m\sns{-1}c^0\ps{1}\ns{0},\\
 &\vs_j(m\ot \td{d}\ot  \td{c})= m\ot \td d \ot m\ot  c^0\odots \ve(c^{j})\odots c^q,\\
 &\vta(m\ot \td{d}\ot \td{c})= m\sns{0}\ot S^{-1}(c\ns{-1})\cdot\td d\ot c^1\ot  c^2\odots c^q\ot  m\sns{-1}c^0\ns{0}.
\end{align*}
\smallskip

\begin{lemma}\label{cyclicoperators}
The  horizontal and vertical  operators defined in \eqref{hd},
\dots, \eqref{cvt} are well-defined and the $\tau$-operators are
invertible.
\end{lemma}
\begin{proof}
In view of Lemma \ref{AYD},  $M\ot \Cc^{\ot q}$ is AYD module and,
since the $q$th row of the above bigraded complex is the Hopf cyclic
complex of $\Kc$ with coefficients in $M\ot \Cc^{\ot q+1}$,  all
horizontal operators are well-defined \cite{hkrs2}. By contrast, the
columns are not Hopf cyclic modules of coalgebras in  general,
except in some special cases such as the case  in Subsection
\ref{uchern}.

Let us check that the vertical $\tau$-operator is well-defined,
which implies that all the others are well-defined. One has
\begin{align*}
&\vta(mh\ot \td d\ot  c^0\odots c^q)= \\
& (mh)\sns{0}\ot S^{-1}(c^0\ns{-1})\td d\ot c^1\odots c^n\ot
(mh)\sns{-1}c^0\ns{0}=\\
&m\sns{0}h\ps{2}\ot S^{-1}(c^0\ns{-1})\td d\ot c^1\odots c^n\ot
S(h\ps{3})m\sns{-1}h\ps{1}c^0\ns{0}=\\
&m\sns{0}\ot h\ps{2}S^{-1}(c^0\ns{-1})\td d\ot h\ps{3}c^1\odots
h\ps{n+2}c^n\ot m\sns{-1}h\ps{1}c^0\ns{0}=\\
&m\sns{0}\ot S^{-1}((h\ps{2}c^0\ns{-1}S(h\ps{4})) h\ps{1}\td d\ot
h\ps{5}c^1\odots h\ps{n+4}c^n\ot \\
&\ot m\sns{-1}(h\ps{3}c^0)\ns{0}=\\
&m\sns{0}\ot S^{-1}((h\ps{2}c^0)\ns{-1}) h\ps{1}\td d\ot
 h\ps{3}c^1\odots h\ps{n+2}c^n\ot m\sns{-1}(h\ps{2})\ns{0}=\\
& \vta(m\ot h\ps{1}\td d\ot h\ps{2}c^0\odots h\ps{n+2}c^n).
\end{align*}
To show that the $\tau$-operators are invertible, we write explicitly
 their inverses:
\begin{align*}
&\vta^{-1}(m\ot \td d\ot \td c)=\\
&m\sns{0}\ot m\sns{-3}c^q\ns{-1}S(m\sns{-1})\td
 d\ot m\sns{-2}c^q\ns{0} c^0\odots c^{q-1}.
\end{align*}
\begin{align*}
\hta^{-1}(m\ot \td d\ot \td c)= m\sns{0}\ot S(\td
c\ns{-1}m\sns{-1})d^p\ot d^0\odots d^{q-1}\ot \td c\ns{0}.
\end{align*}
Using the invertibility of the antipode of $\Hc$ and also the YD
module property of $C$, one checks that
$\vta\circ\vta^{-1}=\vta^{-1}\circ\vta=\Id$ and
$\hta\circ\hta^{-1}=\hta^{-1}\circ\hta=\Id$.
\end{proof}

\begin{proposition}\label{cylind}
Let  $C$  be a (co)module coalgebra and  $D$ a module coalgebra over
$\Hc$. Assume that $C$ is  an YD module over $\Hc$, and $M$ is an AYD
over $\Hc$. Then the bigraded module $\FX^{p,q}$ is a cylindrical module. If in addition
  $M\ot C^{\ot q}$ is stable, then $\FX^{p,q}$ is bicocyclic.
\end{proposition}

\begin{proof}
Lemma \ref{AYD} shows that  $M\ot \Cc^{\ot q}$ is AYD module, and
since the $q$th row of the above bigraded complex is the Hopf cyclic
complex of $C^\ast_\Hc(\Kc; M\ot \Cc^{\ot q+1})$, it defines a
paracocyclic module \cite{hkrs2}.

The columns are not necessarily Hopf cyclic modules of coalgebras
though, except in some special cases. However, one can show that the
columns are paracocyclic modules. The verification of the fact that
$\vta$, $\vd_i$ and $\vs_j$ satisfy \eqref{ds}\dots \eqref{cj} is
straightforward.
 The only nontrivial relations are
those that involve $\vta$ and $\vd_{p+1}$, the others being the same
as for cocyclic module associates to coalgebras. One
needs thus to prove that
\begin{align*}
&\vta\vd_i=\vd_{i-1}\vta,&  1\le i\le q+1,\\
&\vta\vd_0=\vd_{q+1},\\
&\vta \vs_j=\vs_{j-1}\vta, & 1\le j\le q-1,\\
&\vta\vs_0=\vs_{q-1}\vta^2
\end{align*}
 To verify these identities, first let $1\le i\le q$; one has
\begin{align*}
&\vta\vd_i(m\ot \td d\ot  c^0\odots c^q)=\\
&\vta(m\ot \td d\ot c^0\odots c^i\ps{1}\ot c^i\ps{2}\odots c^q)=\\
&m\sns{0}\ot S^{-1}(c^0\ns{-1})\cdot \td{d}\ot c^1\odots
c^i\ps{1}\ot
c^i\ps{2}\odots c^q\ot m\sns{-1}c^0\ns{0}=\\
&\vd_{i-1}(m\sns{0}\ot S^{-1}(c^0\ns{-1})\cdot \td{d}\ot c^1\odots
c^q\ot m\sns{-1}c^0\ns{0})=\\
&\vd_{i-1}\vta(m\ot \td d\ot  c^0\odots c^q).
\end{align*}
Next let $i=q+1$; by using the fact that $C$ is $\Hc$ module
coalgebra, one obtains
\begin{align*}
&\vta\vd_{q+1}(m\ot \td d\ot  c^0\odots c^q)=\\
&\vta(m\sns{0}\ot S^{-1}(c^0\ps{1}\ns{-1})\td d\ot
 c^0\ps{2}\odots c^q\ot m\sns{-1}c^0\ps{1}\ns{0})=\\
&m\sns{0}\ot S^{-1}(c^0\ps{2}\ns{-1})S^{-1}(c^0\ps{1}\ns{-1})\cdot
\td{d}\ot c^1\odots\\
& c^q\ot m\sns{-2}c^0\ps{1}\ns{0}\ot m\sns{-1}c^0\ps{2}\ns{0}=\\
&m\sns{0}\ot S^{-1}(c^0\ps{1}\ns{-1}c^0\ps{2}\ns{-1})\cdot
\td{d}\ot c^1\odots\\
& c^q\ot m\sns{-2}c^0\ps{1}\ns{0}\ot m\sns{-1}c^0\ps{2}\ns{0}=\\
&m\sns{0}\ot S^{-1}(c^0\ns{-1})\cdot \td{d}\ot c^1\odots\\
& c^q\ot m\sns{-2}c^0\ns{0}\ps{1}\ot m\sns{-1}c^0\ns{0}\ps{2}=\\
&\vd_{q}(m\sns{0}\ot S^{-1}(c^0\ns{-1})\cdot \td{d}\ot c^1\odots
 c^q\ot m\sns{-2}c^0)=\\
 &\vd_q\vta(m\ot \td d\ot  c^0\odots c^q).
\end{align*}
Finally let $i=0$; one has
\begin{align*}
&\vta\vd_0(m\ot \td d\ot  c^0\odots c^q)=\\
&\vta(m\ot \td d\ot  c^0\ps{1}\ot c^0\ps{2}\ot c^1\odots c^q)= \\
&m\sns{0}\ot S^{-1}(c^0\ps{1}\ns{-1})\td{d}\ot c^0\ps{2}\ot
c^1\odots c^q\ot m\sns{-1}c^0\ps{1}\ns{0}=\\
&\vd_{q+1}(m\ot \td d\ot  c^0\odots c^q).
\end{align*}
The other identities are checked in a similar fashion.

We next show that the vertical operators commute with
 horizontal operators. Using Lemma \ref{cyclicoperators}, and since in any parcocyclic
module with  invertible $\tau$-operator  one has
\begin{align}\label{conj}&\p_j= \tau^{-j}\p_0\tau^j,\quad 1\le j\le
n, &\s_i= \tau^{-i}\s_{n-1}\tau^i,\quad 1\le i\le n-1,
\end{align}
it suffices to verify the identities
\begin{align*}
&\vta\hta=\hta\vta , \, \quad \vta\hd_0=\hd_0\vta , \, \quad
\hta\vd_0=\vd_0\hta ,  \\
&\vta\hs_{p-1}=\hs_{p-1}\vta ,
\quad \text{and} \quad \hta\vs_{q-1}=\vs_{q-1}\hta .
\end{align*}
Let us check  the first; one has
\begin{align*}
&\vta\hta(m\ot k^0\odots k^p\ot c^0\odots c^q)=\\
&\vta(m\sns{0}\ot k^1\odots k^p\ot \td c\ns{-1}m\sns{-1}k^0\ot \td
c\ns{0})=\\
&m\sns{0}\ot S^{-1}(c^0\ns{-1})\cdot(k^1\odots k^p)\ot
c^1\ns{-1}\dots  c^q\ns{-1}m\sns{-2}k^0\ot\\
& c^1\ns{0}\odots c^q\ns{0}\ot m\sns{-1}c^0\sns{0}=\\
&m\sns{0}\ot  S^{-1}(c^0\ns{-2})\cdot(k^1\odots k^p)\ot\\
&\ot c^1\ns{-1}\dots c^q\ns{-1} m\sns{-4}c^{0}\ns{-1}S(m\sns{-2})
 m\sns{-1}S^{-1}(c^0\ns{-2})k^0\ot \\
& \ot c^1\ns{0}\odots c^q\ns{0}\ot  m\sns{-3}c^0\ns{0}=\\
 &m\sns{0}\ot  S^{-1}(c^0\ns{-1})\cdot(k^1\odots k^p)\ot\\
&\ot c^1\ns{-1}\dots c^q\ns{-1}(m\sns{-2}c^0\ns{0})\ns{-1}
  m\sns{-1}S^{-1}(c^0\ns{-2})k^0\ot \\
& \ot c^1\ns{0}\odots c^q\ns{0}\ot  (m\sns{-2}c^0)\ns{0}=\\
&\hta(m\sns{0}\ot S^{-1}(c^0\ns{-1})\cdot\td
 k\ot c^1\odots c^q\ot m\sns{-1}c^0)=\\
 &\hta\vta(m\ot k^0\odots k^p\ot c^0\odots c^q)
\end{align*}
Since  $C$ is a $\Hc$-module coalgebra one can write
\begin{align*}
&\vta\hd_0(m\ot d^0\odots d^p\ot c^0\odots c^q)=\\
&\vta(m\ot d^0\ps{1}\ot d^0\ps{2}\ot d^1\odots d^p\ot c^0\odots
c^q)=\\
&(m\sns{0}\ot S^{-1}(c^0\ns{-1})(d^0\ps{1}\ot d^0\ps{2}\ot
d^\ot\dots\\& \dots \ot d^p)\ot
 c^1\odots c^q\ot m\sns{-1}c^0\ns{0})=\\
 &\hd_0\vta(m\ot d^0\odots d^p\ot c^0\odots c^q).
\end{align*}
 To show the $\hta\vd_0=\vd_0\hta$ one uses only the module coalgebra
 property of $C$; thus,
\begin{align*}
&\hta\vd_0(m\ot \td d\ot \td c)=\hta(m\ot \td d\ot  c^0\ps{1}\ot
c^0\ps{2}\ot c^1\odots
c^q)=\\
&m\sns{0}\ot d^1\odots d^p\ot
c^0\ps{1}\ns{-1}c^0\ps{2}\ns{-1}c^1\ns{-1}\dots\\
&\dots c^q\ns{-1}m\ns{-1}d^0\ot c^0\ps{1}\ns{0}\ot
c^0\ps{2}\ns{0}\ot c^1\ns{0}\odots c^q\ns{0}=\\
&m\sns{0}\ot d^1\odots d^p\ot
c^0\ns{-1}\ns{-1}c^1\ns{-1}\dots\\
&\dots c^q\ns{-1}m\ns{-1}d^0\ot c^0\ns{0}\ps{1}\ot
c^0\ns{0}\ps{2}\ot c^1\ns{0}\odots c^q\ns{0}= \\
&\vd_0\hta(m\ot \td d\ot \td c).
\end{align*}
The remaining relations, $\hta\vs_{q-1}=\vs_{q-1}\hta$
and $\vta\hs_{p-1}= \hs_{p-1}\vta $, are obviously true.

Finally, by using the stability of $M$, we verify the cylindrical
condition $\nobreak{\, \hta^{p+1}\vta^{q+1}=\Id}$ as follows:
\begin{align*}
&\hta^{p+1}\vta^{q+1}(m\ot d^0\odots d^p\ot c^0\odots c^q)=\\
&\hta^{p+1}(m\sns{0}\ot S^{-1}(c^0\ns{-1}\dots
c^q\ns{-1})\cdot(d^0\odots d^p)\ot\\
&\ot m\sns{-q-1}c^0\ns{0}\odots
 m\sns{-1}c^q\ns{0})=\\
 &m\sns{0}\ot m\sns{-p-q-2}d^0\odots m\sns{-q-2}d^p\ot m\sns{-q-1}c^0\odots
 m\sns{-1}c^q=\\
 &m\ot d^0\odots d^p\ot c^0\odots c^q.
\end{align*}

\end{proof}
The diagonal of any cylindrical  module is a  cocyclic module
\cite{gj} whose cyclic structures is given by,
\begin{align}
&\p_i:= \vd_i\hd_i: \FX^{n,n}\ra \FX^{n+1,n+1},\\
&\s_i:= \vs_i\hd_i: \FX^{n,n}\ra \FX^{n-1,n-1},\\
&\tau:= \vta\hta: \FX^{n,n}\ra \FX^{n,n}.
\end{align}

In order to relate
 the Hopf cyclic cohomology of $\nobreak{C\cl D}$ to the diagonal complex of
the $\FX$, we define a map similar to the one used in \cite{mr}, and
use the fact that $C$ is a YD module over $\Hc$ to show that is
well-defined. Specifically,
\begin{align}\label{Psi}
&\Psi: C^n_\Hc(\nobreak{C\cl D}; M)\ra \FX^{n,n},\\\notag
 &\Psi(m\ot
c^0\cl d^0\odots c^n\cl d^n)=\\\notag
 &m\ot  c^0\ns{-n-1}d^0\odots
c^0\ns{-1}\dots c^n\ns{-1}d^n\ot c^0\ns{0}\odots c^n\ns{0}.
\end{align}

\begin{proposition}\label{diagonal}
The  map defined in \eqref{Psi}  establishes a cyclic isomorphism
between the complex
 $\nobreak{C^\ast_\Hc(C\cl D; M)}$ and the diagonal of $\FX^{\ast,\ast}$.
\end{proposition}
\begin{proof}
First we show that the above map is well-defined.
 The fact that $C$ is
YD module helps in two ways: firstly,  Lemma \ref{daction} shows
that $\Hc$ is acting diagonally on $\nobreak{C\cl D}$, and  secondly
the twisting $\top: C\ot D\ra D\ot C$ where$\top(c\ot
d)=c\ns{-1}d\ot c\ns{0}$, is $\Hc$-linear; indeed,
\begin{align*}
&\top (h\ps{1}c\ot\ps{2}d)= (h\ps{1}c)\ns{-1}(h\ps{2}d)\ot
(h\ps{1})\ns{0}=\\
& h\ps{1}c\ns{-1}S(h\ps{3}) h\ps{4}d\ot h\ps{2}c\ns{0}=\\
&h\ps{1}c\ns{-1}d\ot h\ps{2}c\ns{0}.
\end{align*}
This ensures that $\Psi$ is well-defined,  because it is obtained
out of $\top$ by iteration.

In order to prove that $\Psi$ is a cyclic map, it
suffices to check that $\Psi$ commutes with the
$\tau$-operatorss, the first
coface,  and the last codegeneracy, because the rest of the operators
are  made of these \eqref{conj}. One verifies that
$\Psi$ commutes with cyclic operators as follows. On the one hand,
\begin{align*}
&\Psi\tau_{\nobreak{C\cl D}}(m\ot c^0\cl d^0\odots c^n\cl d^n)=\\
&\Psi(m\sns{0}\ot c^1\cl d^1\odots c^n\cl d^n\ot m\ns{-1}c^0\cl
m\ns{-2}d^0)=\\
&m\sns{0}\ot  c^1\ns{-n-1}d^1\odots  c^1\ns{-2}\dots c^n\ns{-2}
d^n\ot\\
& \ot  c^1\ns{-1}\dots c^n\ns{-1} (m\sns{-2}c^0)\ns{-1}m\sns{-1}d^0
\ot \\
&\ot c^1\ns{0}\odots c^n\ns{0}\ot (m\sns{-2}c^0)\ns{-1}=\\
&m\sns{0}\ot  c^1\ns{-n-1}d^1\odots  c^1\ns{-2}\dots c^n\ns{-2}
d^n\ot \\
&  \ot c^1\ns{-1}\dots c^n\ns{-1} m\sns{-2}c^0\ns{-1} d^0 \ot
c^1\ns{0}\odots c^n\ns{0}\ot m\sns{-1}c^0\ns{0}.
\end{align*}
On the other hand,
\begin{align*}
&\hta\vta\Psi(m\ot c^0\cl d^0\odots c^n\cl d^n)=\\
&\hta\vta(m\ot  c^0\ns{-n-1}d^0\odots   c^0\ns{-1}\dots
c^n\ns{-1}d^n\ot c^0\ns{0}\odots c^n\ns{0})=\\
&\hta(m\sns{0}\ot
  c^0\ns{-n-1}c^1\ns{-n-1}d^1\odots c^0\ns{-2}\dots
c^n\ns{-2}d^n\ot \\
& c^0\ns{-1}\dots c^n\ns{-1} m\sns{-1}c^0\ns{-n-2}d^0 \ot
c^0\ns{0}\odots c^n\ns{0})=\\
& m\sns{0}\ot
 S^{-1}(c^0\ns{-1})\cdot (c^0\ns{-n-1}c^1\ns{-n-1}d^1\odots c^0\ns{-2}\dots
c^n\ns{-2}d^n\ot \\
& c^0\ns{-1}\dots c^n\ns{-1} m\sns{-1}c^0\ns{-n-2}d^0) \ot
c^1\ns{0}\odots c^n\ns{0} \ot m\sns{-1}c^0\ns{0}=\\
&m\sns{0}\ot  c^1\ns{-n-1}d^1\odots  c^1\ns{-2}\dots c^n\ns{-2}
d^n\ot \\
&  \ot c^1\ns{-1}\dots c^n\ns{-1} m\sns{-2}c^0\ns{-1} d^0 \ot
c^1\ns{0}\odots c^n\ns{0}\ot m\sns{-1}c^0\ns{0}.
\end{align*}

The next to check is the equality
${\p_0}\Psi=\Psi{\p_0}_{\nobreak{C\cl D}}$. Using the fact that $C$
is a module coalgebra, one has
\begin{align*}
&\p_0\Psi(m\ot m\ot c^0\cl d^0\odots c^n\cl d^n)=\\
&\p_0(m\ot  c^0\ns{-n-1}d^0\odots   c^0\ns{-1}\dots c^n\ns{-1}d^n\ot
c^0\ns{0}\odots c^n\ns{0})=\\
&m\ot  c^0\ns{-n-2}d^0\ps{1}\ot  c^0\ns{-n-1}d^0\ps{2}\odots
c^0\ns{-1}\dots c^n\ns{-1}d^n\ot\\
&c^0\ns{0}\ps{1}\ot c^0\ns{0}\ps{2}\ot c^1\ns{0}\odots c^n\ns{0}=\\
&m\ot  c^0\ps{1}\ns{-n-2}c^0\ps{2}\ns{-n-2}d^0\ps{1}\ot
c^0\ps{1}\ns{-n-1}c^0\ps{2}\ns{-n-1}d^0\ps{2}\odots\\
&c^0\ps{1}\ns{-1}c^0\ps{2}\ns{-1}c^1\ns{-1}\dots c^n\ns{-1}d^n\ot
c^0\ps{1}\ns{0}\ot c^0\ps{2}\ns{0}\ot c^1\ns{0}\ot\dots\\& \dots\ot
c^n\ns{0}.
\end{align*}
On the other hand,
\begin{align*}
&\Psi{\p_0}_{\nobreak{C\cl D}}(m\ot c^0\cl d^0\odots c^n\cl d^n)= \\
&\Psi(m\ot c^0\ps{1}\cl c^0\ps{2}\ns{-1}d^0\ps{1}\ot
c^0\ps{2}\ns{0}\cl d^0\ps{2}\ot  c^1\cl d^1\ot\dots\\
&\dots\ot c^n\cl d^n )=\\
&m\ot  c^0\ps{1}\ns{-n-2}c^0\ps{2}\ns{-1}d^0\ps{1}\ot
c^0\ps{1}\ns{-n-1}c^0\ps{2}\ns{0}\ns{-n-1}d^0\ps{2} \ot\\
& c^0\ps{1}\ns{-n}c^0\ps{2}\ns{0}\ns{-n}c^1\ns{-n}d^1 \odots
c^0\ps{1}\ns{-1}c^0\ps{2}\ns{0}\ns{-1}\dots c^n\ns{-1}d^n\ot\\
&c^0\ps{1}\ns{0}\ot c^0\ps{2}\ns{0}\ns{0}\ot c^1\ns{0}\odots
c^n\ns{0}=\\
&m\ot  c^0\ps{1}\ns{-n-2}c^0\ps{2}\ns{-n-2}d^0\ps{1}\ot
c^0\ps{1}\ns{-n-1}c^0\ps{2}\ns{-n-1}d^0\ps{2}\odots\\
&c^0\ps{1}\ns{-1}c^0\ps{2}\ns{-1}c^1\ns{-1}\dots c^n\ns{-1}d^n\ot
c^0\ps{1}\ns{0}\ot c^0\ps{2}\ns{0}\ot c^1\ns{0}\ot\dots\\
&\dots\ot c^n\ns{0}.
\end{align*}

The equality $\s_{n-1}\Psi=\Psi{\s_{n-1}}_{\nobreak{C\cl D}}$ is
obvious.

 In order to show that $\Psi$ is an isomorphism,  one constructs  again
an iterated map
 made  out of factors
 $\perp: D\ot C\ra C\ot D$, defined by $\perp(d\ot c)=c\ns{0}\ot
 S^{-1}(c\ns{-1})d$. It is easy to see that $\top $ and $\perp$ are
 inverse to one another. Explicitly the inverse of $\Psi $ is defined by
 \begin{align*}
 &\Psi^{-1}(m\ot d^0\odots d^n\ot c^0\odots c^n)=\\
 &=m\ot c^0\ns{0}\cl S^{-1}(c^0\ns{-1})d^0\ot\dots\\
&\dots\ot c^n\ns{0}\cl S^{-1}(c^0\ns{-n-1}c^1\ns{-n}\dots
 c^n\ns{-1})d^n
 \end{align*}
\end{proof}

Applying now the cyclic version of the Eilenberg-Zilber  theorem \cite{gj} one
obtains the sought-for quasi-isomorphism of mixed complexes.

  \begin{proposition} \label{mix+}
  The mixed complexes
    $(C^*_\Hc(D\cl C,M), \,b,\, B)$  and \\
    $(\Tot(\FX), \,b_T,\, B_T))$ are quasi-isomorphic.
  \end{proposition}

\bigskip

The rest of this section is devoted to showing that the Hopf algebras
that make the object of this
paper satisfy  the conditions of the above proposition. To put
this in the proper setting, we let $\Hc$ be a Hopf algebra, and
we consider a pair of $\Hc$-module coalgebras $C,D$ such
that $C$ is $\Hc$-module coalgebra and via its action  and coaction
it is an YD module over $\Hc$. Let
$\Lc\subset \Kc$ be Hopf subalgebras of $\Hc$. One defines the coalgebra
$\Cc:= \Hc\ot_\Kc\Cb$,  where $\Kc$ acts on $\Hc$ by multiplication
and on $\Cb$ via counit (cf \cite[\S 5]{cm6}). In the same fashion one defines
the coalgebra  $\KL :=\Kc\ot_\Lc\Cb$. If $h\in \Hc$ and $c=\dot h$
be its class in $\Cc$, then
\begin{align*}
\D(c)= c\ps{1}\ot c\ps{2}:= \dot h\ps{1}\ot \dot h\ps{2}, \quad
\ve(c):=\ve(h)
\end{align*}
This coalgebra has  a natural coaction from  $\Hc$.
\begin{align*}
\Db(c)=c\ns{-1}\ot c\ns{0}:= h\ps{1}S(h\ps{3})\ot \dot h\ps{2}\,.
\end{align*}

\begin{lemma}
The above action and coaction are well-defined and  make  $\Cc$ a
(co)module coalgebra, respectively. In addition both make $\Cc$ an
YD module over $\Hc$.
\end{lemma}
\begin{proof} First let us check that the action and coaction are well
defined. With  $h,g\in \Hc$,  and $k\in \Kc$, one has
\begin{align*}
h( gk\ot 1)= hgk\ot_\Kc 1= hg\ot_\Kc\epsilon(k)= h(g\ot
\epsilon(k)) ,
\end{align*}
which verifies the claim for the action. For the coaction, we write
\begin{align*}
&\Db(hk\ot 1)= h\ps{1}k\ps{1}S(k\ps{3})S(h\ps{3})\ot
(h\ps{2}k\ps{2}\ot_\Kc 1)=\\
 &h\ps{1}k\ps{1}S(k\ps{3})S(h\ps{3})\ot
(h\ps{2}\ot_\Kc \epsilon(k\ps{2}))= \\
&= h\ps{1}k\ps{1} S(k\ps{2})
S(h\ps{3})\ot (h\ps{2} \ot_\Kc 1)=\\
& h\ps{1} S(h\ps{3}))\ot (h\ps{2} \ot_\Kc \epsilon(k))=\Db(h\ot \epsilon(k)).
\end{align*}
It is obvious that these are indeed action, resp. coaction, and thus define
a module coalgebra structure, resp. a
module coalgebra structure over $\Cc$.  One checks that the
action and coaction satisfy \eqref{yd} as follows:
\begin{align*}
&\Db(h( \dot g)=\Db( \wid{hg})=
h\ps{1}g\ps{1}S(g\ps{3})S(h\ps{3})\ot
\wid{h\ps{2}g\ps{2}}=\\
& h\ps{1}(g\ps{1}S(g\ps{3}))S(h\ps{3})\ot h\ps{2}\dot{g\ps{2}}.
\end{align*}
\end{proof}

  Let us assume that  $\Hc$  acts on $\Kc$, as well and
   via this action $\Kc$ is  a module
coalgebra. Similarly to $\Cc$ in relation to $\Hc$, $\KL$ inherits in a
natural way an action from $\Hc$. This makes $\KL$ a module
coalgebra over $\Hc$. Precisely,  if $h\in \Hc$ and $k\in \Kc$,
denoting $k\ot_\Lc 1$ by $\dot k$,  one defines
\begin{align*}
h\cdot \dot k= hk\ot_\Lc 1= \wid{hk}.
\end{align*}

 Now $\Cc$ and  $\KL$ together with their action and
coaction satisfy all conditions of Proposition \ref{cylind}. As a
result one can form the
 cylindrical module  $\FX(\Hc, \Cc,\KL; M)$.
Due to the special properties discussed above, we may expect
some simplification.
Indeed, let
\begin{align*}
\FY^{p,q}=M\ot_\Kc \KL^{\ot p+1}\ot \Cc^{\ot q}
\end{align*}
We define the  following map from $\FX$ to $\FY$:
\begin{align}\label{Phi_1}
&\Phi_1:\FX^{p,q}\ra \FY^{p,q},\\\notag &\Phi_1(m\ot_\Hc \td {\dot
k}\ot \dot{h^0}\odots \dot{h^n})= \\\notag &mh^0\ps{2}\ot_\Kc
S^{-1}(h^0\ps{1})\cdot \td{ \dot k}\ot
S(h^0\ps{3})\cdot(\dot{h^1}\odots \dot{h^q}).
\end{align}
\begin{lemma}
The map $\Phi_1$ defined in \eqref{Phi_1} is a well-defined
isomorphisms of vector spaces.
\end{lemma}
\begin{proof}
To check that it is well-defined, we clarify the ambiguities
in the definition of $\Phi_1$ as follows. Let $k\in \Kc$, $\td{\dot
k}\in \KL^{\ot (p+1)}$, and $g, h^0\dots h^q\in \Hc$; then
\begin{align*}
&\Phi_1(m\ot \td {\dot k}\ot (h^0k\ot_\Kc 1)\ot \dot{h^1} \odots \dot{h^q} ))=\\
&mh^0\ps{2}k\ps{2}\ot_\Kc S^{-1}(h^1\ps{1}k\ps{1})\cdot\td{\dot
k}\ot
S(h^0\ps{3}k\ps{3})\cdot (\dot{h^1}\odots \dot{h^q})=\\
&mh^0\ps{2}\ot_\Kc k\ps{2}S^{-1}(k\ps{1})S^{-1}(h^1\ps{1})\cdot\td{
\dot k}\ot k\ps{3}S(k\ps{4})S(h^0\ps{3})\cdot (\dot{h^1}\odots
\dot{h^q})=\\
& \Phi_1(m\ot \td {\dot k}\ot (h^0\ot_\Kc \epsilon(k))\ot \dot{h^1}
\odots \dot{h^n} )).
\end{align*}
Also,
\begin{align*}
&\Phi_1(m\ot g\ps{1}\cdot \td {\dot k}\ot g\ps{2}\cdot(\dot{h^0}\ot
\dot{h^1}
\odots \dot{h^q} ))=\\
&m\ot_\Kc g\ps{1}\cdot \td {\dot k}\ot g\ps{2}\cdot(\dot{h^0}\ot
\dot{h^1}
\odots \dot{h^q} )=\\
&m g\ps{3}h^0\ps{2}\ot_\Kc
S^{-1}(h^0\ps{1})S^{-1}(g\ps{2})g\ps{1}\cdot\td {\dot k}\ot\\
&\ot S(h^0\ps{3})S(g\ps{4})g\ps{5}\cdot (\dot {h^1}\odots
\dot{h^q})=\\
&m g h^0\ps{2}\ot_\Kc S^{-1}(h^0\ps{1})\cdot\td {\td k}\ot
S(h^0\ps{3})\cdot
 (\dot {h^1}\odots \dot{h^q})=\\
&\Phi_1(mg\ot \td {\dot k}\ot \dot{h^0}\ot \dot{h^1} \odots
\dot{h^q}).
\end{align*}
One easily checks that the following map defines an inverse for $\Phi_1$:
 \begin{align*}
&\Phi_1^{-1}: \FY^{p,q}\ra \FX^{p,q},\\\notag
 &\Phi_1^{-1}(m\ot_\Kc \td {\dot k}\ot \td {\dot{h}})= m\ot_\Hc \td {\dot k}\ot \dot{1}\ot \td {\dot{h}}.
\end{align*}
\end{proof}
\smallskip

We next push forward the cylindrical structure of $\FX$ to get the
following cylindrical structure on $\FY$:
\begin{align*}
&\hd_{i}(m\ot \td{\dot k} \ot \td{{\dot{h}}})= m\ot {\dot k}^0\odots
\D({\dot k}^i)\odots {\dot k}^p\ot \td {\dot{h}},\\
&\hd_{p+1}(m\ot \td{\dot k} \ot \td{{\dot{h}}})=m\sns{0}\ot {\dot
k}^0\ps{2}\odots
{\dot k}^p\ot \td {\dot{h}}\ns{-1}m\sns{-1}\dot{k}^0\ps{1} \ot \td {\dot{h}}\ns{0},\\
 &\hs_j(m\ot \td{\dot k} \ot \td{{\dot{h}}})= m\ot  {\dot k}^1\odots
\ve({\dot k}^{j})\odots {\dot k}^p\ot  \td{{\dot{h}}},\\
 &\hta(m\ot \td{{\dot k}} \ot \td{{\dot{h}}})= m\sns{0}\ot {\dot k}^1
 \odots {\dot k}^p\ot \td {\dot{h}}\ns{-1}m\sns{-1}{\dot k}^0\ot  \td {\dot{h}}\ns{0}
\end{align*}
\begin{align*}
&\vd_{0}(m\ot \td{{\dot k}}\ot\td{{\dot{h}}})=m\ot \td {\dot k}\ot \dot 1\ot {\dot{h}}^1\odots  {\dot{h}}^q,\\
 &\vd_{i}(m\ot \td{{\dot k}}\ot\td{{\dot{h}}})=m\ot \td {\dot k}\ot {\dot{h}}^0\odots \D({\dot{h}}^i)\odots {\dot{h}}^q,\\
&\vd_{q+1}(m\ot \td{{\dot k}}\ot  \td{{\dot{h}}})= m\sns{0}\ot \td {\dot k}\ot  {\dot{h}}^1\odots {\dot{h}}^q\ot m\sns{-1},\\
 &\vs_j(m\ot \td{{\dot k}}\ot  \td{{\dot{h}}})= m\ot \td {\dot k} \ot m\ot  {\dot{h}}^1\odots \ve({\dot{h}}^{j+1})\odots {\dot{h}}^q,\\
 &\vta(m\ot \td{{\dot k}}\ot \td{{\dot{h}}})=\\
 & m\sns{0} h^1\ps{2}\ot S^{-1}(h^1\ps{1})
 \cdot\td {\dot k}\ot S(h^1\ps{3})\cdot({\dot{h}}^2\odots {\dot{h}}^q\ot  m\sns{-1}).
\end{align*}

To obtain a further simplification, we
assume that  the action of $\Kc\subset\Hc$ on $\Kc$  coincides with the
multiplication by $\Kc$.  We define a map from
$\FY^{p,q}$ to
\begin{align*} \FZ^{p,q}:=\FZ^{p,q}(\Hc,
\Kc, \Lc;M):= M\ot_\Lc \KL^{\ot p}\ot \Cc^{\ot q},
\end{align*}
as follows:
\begin{align}\label{Phi_2}
&\Phi_2: \FY^{p,q}\ra \FZ^{p,q},\\\notag & \Phi_2(m\ot_\Kc {\dot
k}^0\odots {\dot k}^p\ot \td {\dot{h}})= mk^0\ps{1}\ot_\Lc
S(k^0\ps{2})\cdot({\dot k}^1\odots {\dot k}^p\ot \td {\dot{h}}).
\end{align}

\begin{lemma}
 The map $\Phi_2$ defined in \eqref{Phi_2} is a well-defined isomorphism of vector spaces.
\end{lemma}
\begin{proof} One has
\begin{align*}
& \Phi_2(m\ot k\cdot( {\dot k}^0\odots
{\dot{k}}^p\ot \td {\dot{h}}))=\\
& \Phi_2(m\ot k \ps{1}\cdot {\dot{k}}^0\ot  k\ps{2}\cdot( {\dot{k}}^1\odots {\dot{k}}^p\ot \td {\dot{h}}))=\\
& m k\ps{1}{\dot{k}}^0\ps{1}\ot S(k\ps{2}{\dot{k}}^0\ps{2})\cdot
(k\ps{3}\cdot({\dot{k}}^2\odots {\dot{k}}^p\ot \td {\dot{h}}))=\\
& m k\ps{1}\epsilon(k\ps{2}{\dot{k}}^0\ps{1})\ot S({\dot{k}}^0\ps{2})\cdot ({\dot{k}}^2\odots {\dot{k}}^p\ot \td {\dot{h}})=\\
&\Phi_2(mk\ot {\dot{k}}^0\odots {\dot{k}}^p\ot \td {\dot{h}}) .
\end{align*}
The inverse of $\Phi_2$ is given by
\begin{align*}
&\Phi_2^{-1}: \FZ^{p,q}\ra \FY^{p,q},\\
&\Phi_2^{-1}(m\ot {\dot{ k}}^1\odots {\dot{k}}^p\ot
{\dot{h}}^1\odots {\dot{h}}^q)= \\
&m\ot_\Kc \dot 1\ot {\dot{k}}^1\odots {\dot{k}}^p\ot
{\dot{h}}^1\odots {\dot{h}}^q.
\end{align*}
\end{proof}

We now push forward the cylindrical structure of   of $\FY$ on $\FZ$
to get the following operators on $\FZ^{\ast,\ast}$ :
\begin{align*}
&\hd_{0}(m\ot \td{{\dot{k}}} \ot \td{{\dot{h}}})= m\ot  \dot 1\ot
{\dot{k}}^1\odots
{\dot{k}}^p\ot \td {\dot{h}},\\
&\hd_{i}(m\ot \td{{\dot{k}}} \ot \td{{\dot{h}}})= m\ot
{\dot{k}}^1\odots
\D({\dot{k}}^{i})\odots {\dot{k}}^p\ot \td {\dot{h}},\\
&\hd_{p+1}(m\ot \td{{\dot{k}}} \ot \td{{\dot{h}}})=m\sns{0}\ot
{\dot{k}}^1\odots
{\dot{k}}^p\ot \td {\dot h}\ns{-1}\wid{m\sns{-1}} \ot \td {\dot{h}}\ns{0},\\
 &\hs_j(m\ot \td{{\dot{k}}} \ot \td{{\dot{h}}})= m\ot  {\dot{k}}^1\odots
\ve({\dot{k}}^{j+1})\odots {\dot{k}}^p\ot  \td{{\dot{h}}},\\
 &\hta(m\ot \td{{\dot{k}}} \ot \td{{\dot{h}}})= m\sns{0}k^1\ps{1}\ot S(k^1\ps{2})\cdot ({\dot{k}}^2\odots
  {\dot{k}}^p\ot \td {\dot{h}}\ns{-1}\wid{m\sns{-1}}\ot  \td
 {\dot{h}}\ns{0}).
\end{align*}
\begin{align*}
&\vd_{0}(m\ot \td{\dot k}\ot\td{{\dot{h}}})=m\ot \td {\dot k}\ot\dot 1\ot {\dot{h}}^1\odots  {\dot{h}}^q,\\
 &\vd_{i}(m\ot \td{\dot k}\ot\td{{\dot{h}}})=m\ot \td {\dot k}\ot
  {\dot{h}}^0\odots \D({\dot{h}}^i)\odots {\dot{h}}^q,\\
&\vd_{q+1}(m\ot \td{\dot k}\ot  \td{{\dot{h}}})= m\sns{0}\ot \td {\dot k}\ot  {\dot{h}}^1\odots {\dot{h}}^q\ot\wid{ m\sns{-1}} ,\\
 &\vs_j(m\ot \td{\dot k}\ot  \td{{\dot{h}}})= m\ot \td {\dot k} \ot  {\dot{h}}^1\odots \ve({\dot{h}}^{j+1})\odots {\dot{h}}^q,\\
 &\vta(m\ot \td{\dot k}\ot \td{{\dot{h}}})=m\sns{0} h^1\ps{4}S^{-1}
(h^1\ps{3}\cdot 1_\Kc)\ot \\
&S(S^{-1}(h^1\ps{2})\cdot 1_\Kc)\cdot\left(
S^{-1}(h^1\ps{1})\cdot\td{\dot k}\ot
S(h^1\ps{5})\cdot({\dot{h}}^2\odots {\dot{h}}^q\ot \wid{ m\sns{-1}})
\right).
\end{align*}

\begin{lemma}
The above operators defined on $\FZ$ are well-defined and yield a
cylindrical module.
 \end{lemma}
\begin{proof}
The second part of the lemma holds by the very definition, in view of the fact that
$\FX$ is cylindrical module. We check the first claim for $\vta$ and $\hta$,
for the other operators being obviously true. Since $M$ is AYD
and $\Cc$ is $YD$, we have
\begin{align*}
&\hta(m\ot l\cdot(\td {\dot k} \ot \td{\dot{h}}))=\\
&m\sns{0}l\ps{1}k^1\ps{1}\ot_\Lc S(l\ps{2}k^1\ps{2})\cdot
(l\ps{3}\dot{k}^2\odots l\ps{p+1}\dot{k}^p\ot\\
& \overbrace{(l\ps{p+2}\td{\dot{h}})\ns{-1}m\sns{-1}}\ot
(l\ps{p+2}\td{\dot{h}})\ns{0})=\\
&m\sns{0}l\ps{1}k^1\ps{1}\ot_\Lc S(l\ps{2}k^1\ps{2})\cdot
(l\ps{3}\dot{k}^2\odots l\ps{p+1}\dot{k}^p\ot\\
& \overbrace{l\ps{p+2}\td{\dot{h}}\ns{-1} S(l\ps{p+4})m\sns{-1}}\ot
l\ps{p+3}\td{\dot{h}}\ns{0})=\\
&m\sns{0}l\ps{1}k^1\ps{1}\ot_\Lc S(k^1\ps{2})\cdot (\dot{k}^2\odots
\dot{k}^p\ot  \overbrace{\td{\dot{h}}\ns{-1} S(l\ps{2})m\sns{-1}}\ot
\td{\dot{h}}\ns{0})=\\
&m\sns{0}l\ps{2}k^1\ps{1}\ot_\Lc S(k^1\ps{2})\cdot (\dot{k}^2\odots
\dot{k}^p\ot  \overbrace{\td{\dot{h}}\ns{-1}
S(l\ps{3})m\sns{-1}}l\ps{1}\ot
\td{\dot{h}}\ns{0})=\\
&\hta(ml\ot \td{\dot k}\ot \td{\dot h}).
\end{align*}

For the vertical cyclic operator one only uses the AYD property of
$M$; thus,
\begin{align*}
&\vta(m\ot l\cdot(\td{{\dot{k}}}\ot \td{{\dot{h}}}))=m\sns{0}
l\ps{5}h^1\ps{4}S^{-1}
(l\ps{4}h^1\ps{3}\cdot 1_\Kc)\ot_\Lc \\
&S(S^{-1}(l\ps{3}h^1\ps{2})\cdot 1_\Kc)\cdot\left(
S^{-1}(l\ps{2}h^1\ps{1})\cdot l\ps{1}\td{\dot{k}}\ot
S(l\ps{6}h^1\ps{5})\cdot(l\ps{7}{\dot{h}}^2\ot\right. \\
&\left. \dots\ot l\ps{5+q}{\dot{h}}^q\ot \dot{m\sns{-1}})\right)=
m\sns{0} l\ps{1}h^1\ps{4}S^{-1}
(h^1\ps{3}\cdot 1_\Kc)\ot_\Lc \\
&S(S^{-1}(h^1\ps{2})\cdot 1_\Kc)\cdot\left( S^{-1}(h^1\ps{1})\cdot
\td{\dot{k}}\ot S(h^1\ps{5})\cdot({\dot{h}}^2\ot\dots \right.\\
&\left.\dots\ot{\dot{h}}^q\ot S(l\ps{2})\wid{m\sns{-1}})\right)=\\
&m\sns{0} l\ps{2}h^1\ps{4}S^{-1} (h^1\ps{3}\cdot 1_\Kc)\ot_\Lc
S(S^{-1}(h^1\ps{2})\cdot 1_\Kc)\cdot\left( S^{-1}(h^1\ps{1})\cdot
\td{\dot{k}}\ot\right.\\
&\left.\ot S(h^1\ps{5})\cdot({\dot{h}}^2\odots {\dot{h}}^q\ot
\wid{S(l\ps{3})m\sns{-1}l\ps{1}})\right)=\vta(ml\ot
\td{{\dot{k}}}\ot \td{{\dot{h}}}).
\end{align*}
\end{proof}

\begin{proposition} The following map defines
an isomorphism of cocyclic modules,
\begin{equation}\label{teta}
\Theta:= \Phi_2\circ\Phi_1\circ\Psi:C^\ast_\Hc(\Cc\acl\Kc_\Lc; M)\ra
\FZ_\Lc^{\ast,\ast},
 \end{equation}
where $\Psi$, $\Phi_1$, and $\Phi_2$ are defined in
\eqref{Psi},\eqref{Phi_1} and \eqref{Phi_2} respectively.

\end{proposition}
\begin{proof}
The maps $\Ph_2$, $\Phi_1$, and $\Psi$ are isomorphisms.
\end{proof}

Define now the map
\begin{align} \notag
&\Phi:\Hc\ot_\Lc 1\ra \Cc\cl \KL, \qquad \text{by the formula}\\ \label{phimap}
&\Phi(h\ot_\Lc 1)=
\wid{h\ps{1}}\cl \wid{h\ps{2}\cdot 1_\Kc}= (h\ps{1}\ot_\Kc 1)\cl
(h\ps{2}\cdot 1_\Kc\ot_\Lc 1).
\end{align}
\begin{proposition}\label{phi}
The  map $\Phi$ is a map of  $\Hc$-module coalgebras.
\end{proposition}
\begin{proof}
Using that if $h\in \Hc$ and $l\in \Lc$  then  $(hl)\cdot 1_K= (h\cdot
1_\Kc)l$, one checks that  $\Phi$
is well-defined, as follows:
\begin{align*}
&\Phi(hl\ot 1)= (h\ps{1}l\ps{1}\ot_\Kc 1)\cl (h\ps{2}l\ps{2}\cdot
1_\Kc\ot_\Lc 1)=\\
&(h\ps{1}\ot_\Kc \epsilon(l\ps{1}))\cl ((h\ps{2}\cdot 1_\Kc)
l\ps{2}\ot_\Lc 1)=\\
&(h\ps{1}\ot_\Kc \epsilon(l\ps{1}))\cl ((h\ps{2}\cdot 1_\Kc) \ot_\Lc
\epsilon(l\ps{2}))= \Phi(h\ot \epsilon(l)).
\end{align*}

Next, using the coalgebra
structure of $\Cc\cl\KL$, one can write
\begin{align*}
&\D(\Phi(h))= \D(\wid{ h\ps{1}}\cl  \wid{h\ps{2}\cdot 1_\Kc})= \\
&\wid{h\ps{1}}\acl\wid{ h\ps{2}S(h\ps{4})\cdot (h\ps{5}\cdot 1_\Kc)}\ot \wid{ h\ps{3}}\cl \wid{h\ps{6}\cdot 1_\Kc}=\\
&\wid{ h\ps{1}}\acl \wid{(h\ps{2}S(h\ps{4})h\ps{5})\cdot 1_\Kc}\ot \wid{ h\ps{3}}\cl \wid{h\ps{6}\cdot 1_\Kc}=\\
&\wid{ h\ps{1}}\acl \wid{h\ps{2}\cdot1_\Kc}\ot \wid{ h\ps{3}}\cl\wid{ h\ps{4}\cdot 1_\Kc}=\\
&\Ph(h\ps{1})\ot \Phi(h\ps{2}).
\end{align*}
Finally, the $\Hc$-linearity is trivial because $\Hc$ acts on $\Cc\cl \KL$
diagonally.
\end{proof}

Summing up, we conclude with the following result, which applies to a many cases of
interest, in particular to those which make the main object of
this paper.

\begin{theorem}\label{quasi}
Assuming that $\Phi$ is an isomorphism, the Hopf cyclic complex of
the Hopf algebra $\Hc$ relative to the Hopf subalgebra $ \Lc$ with
coefficients in SAYD $M$ is quasi-isomorphic with the total complex
of the mixed complex of  $\FZ(\Hc,\Kc, \Lc;M)$.
\end{theorem}

\bigskip
 As a matter of fact, the cylindrical modules  $\FZ$ and $\FY$ are often
 bicocyclic modules,
   for example if  the   SAYD  $M$ has the  property that
   $$m\sns{-1}\ot m\sns{0}\in
\Kc\ot M , \qquad \fl \, m\in M;
$$
in this case $M$ will be called $\Kc$-SAYD. In this paper we compute
Hopf cyclic cohomology with coefficients in $\Cb_\d$, which is
obviously $\Kc$-SAYD for any Hopf subalgebra $\Kc\subset \Hc$.

 \begin{proposition}
If $M$ is $\Kc$-SAYD module then $M\ot C^{\ot q}$ is SAYD, and hence
$\FZ$ and $\FY$ are  bicocyclic modules.
\end{proposition}
\begin{proof}
Only the stability condition
remains to be proved. We check it only for $q=1$, but
 the same proof works for all $q\ge 1$.
 By the very definition of the coaction,
\begin{align*}
&\D( m\ot h)=h\ps{1}S(h\ps{3}) m\sns{-1}\ot m\sns{0}\ot h\ps{2}
\end{align*}
We need to verify the identity
\begin{align*}
(m\ot h)\sns{0}(m\ot h)\sns{-1} \, = \, m\ot h .
\end{align*}
The left hand side can be expressed as follows:
\begin{align*}
 &(m\ot h)\sns{0}(m\ot h)\sns{-1}=\\
 &m\sns{0}  h\ps{1}S(h\ps{5})m\sns{-2}\ot
 S(m\sns{-1})S^2(h\ps{4})S(h\ps{2})h\ps{3}=\\
 &m\sns{0}  h\ps{1}S(h\ps{3})m\sns{-2}\ot
 S(m\sns{-1})S^2(h\ps{2}).
\end{align*}
 We now recall that the cyclic operator
 $\tau_1: M\ot \Hc\ra M\ot \Hc$ is given by  $\tau_1(m\ot
h)=m\sns{0}h\ps{1}\ot S(h\ps{2})m\sns{-1}$. Since
\begin{align*}
&\tau_1^2(m\ot h)= \tau_1(m\sns{-1}h\ps{1}\ot S(h\ps{1})m\sns{-1})=\\
&m\sns{0}h\ps{2}S(h\ps{5})m\sns{-3}\ot
S(m\sns{-2})S^2(h\ps{4})S(h\ps{3})m\sns{-1}h\ps{1},
\end{align*}
the desired equality simply follows from the facts that
$\tau_1^2=\Id$ and $M$ is $\Kc$-SAYD.
\end{proof}
\subsection{Bicocyclic complex for primitive Hopf algebras}

We now proceed to show that
the Hopf algebras $\Hc (\Pi)$ associated to
a primitive pseudogroups do satisfy all the requirements of the preceding
subsection, and so Theorem \ref{quasi} applies to allow the computation
of their Hopf cyclic cohomology by means of a bicocyclic
complex.

First we note that,
in view of Proposition \ref{cop}, we can replace
 $\Hc (\Pi)$ by $\Hc(\Pi)^{\rm cop}$.
By Theorems \ref{bicrossed} and \ref{cnbicrossed}, the latter
can be identified to the
bicrossed product $\Fc\acl \Uc$, where $\, \Uc:=U(\Fg)$,
$\, \Fc:=F(N)$, and we shall do so from now on without further
warning.

Next, we need both $\Fc$ and $\Uc$ to be $\Hc$-module coalgebras.
Recalling that $\Fc$ is a left $\Uc$-module, we define an action
of $\Hc$ on  $\Fc$ by the formula
\begin{align} \label{newact}
(f\acl u)g= f \, (u\rt g) , \qquad f , g \in \Fc , \, u \in \Uc .
\end{align}

\begin{lemma}  The above formula defines an action, which
makes $\Fc$ an $\Hc$-module coalgebra.
\end{lemma}

\proof
We check directly the action axiom
\begin{align*}
&(f^1\acl u^1)\big((f^2\acl u^1)g\big)= (f^1\acl u^1)(f^2u^2\rt g)=\\
&f^1 u^1\rt(f^2u^1\rt g )= f^1(u^1\ps{1}\rt f^2) (u^1\ps{2}u^2\rt g)=\\
&\left[f^1u^1\ps{1}\rt f^2\acl u^1\ps{2}u^2\right]g= \left[(f^1\acl
u^1)((f^2\acl u^1)\right]g ,
\end{align*}
and the property of being a Hopf action:
\begin{align*}
&\D((f\acl u)g)= \D(f u\rt g)= f\ps{1}u\ps{1}\ns{0}g\ps{1}\ot
f\ps{2}u\ps{1}\ns{1}(u\ps{2}\rt g\ps{2})=\\
&(f\ps{1}\acl u\ps{1}\ns{0})g\ps{1}\ot (f\ps{2}u\ps{1}\ns{1}\acl
u\ps{2})g\ps{2}.
\end{align*}
\endproof

 To realize $\Uc$ as an $\Hc$-module and comodule  coalgebra, we identify it
 with $\Cc:= \Hc\ot_ \Fc\Cb$, via the map $\, \natural:\Hc\ot_\Fc \Cb\ra \Uc$
 defined by
  \begin{align}
 \natural(f\acl u\ot_\Fc 1)=\epsilon(f)u .
\end{align}

 \begin{lemma} \label{natural}
 The map $\, \natural:\Cc \ra \Uc$ is an isomorphism of
 coalgebras.
\end{lemma}

\begin{proof}
The map is well-defined, because
\begin{align*}
&\natural((f\acl u)(g\acl 1)\ot_\Fc 1)= \natural(fu\ps{1}\rt g\acl
u\ps{2}\ot_\Fc 1)=\\
&\epsilon(fu\ps{1}\rt g)u\ps{2}=\epsilon(f)\epsilon(g)u=
\natural(f\acl u\ot_\Fc \epsilon(g)).
\end{align*}
Let us check that $\natural^{-1}(u)=(1\acl u)\ot_\Fc 1$ is its
inverse. It is easy to see that $\natural\circ
\natural^{-1}=\Id_\Uc$. On the other hand,
\begin{align*}
&\natural^{-1}\circ\natural (f\acl u \ot_\Fc
1)=\natural^{-1}(\epsilon(f)u)=\\
& (1\acl u)\ot_\Fc \epsilon(f)= (1\acl
 u\ps{2})\ot_\Fc \epsilon(S^{-1}(u\ps{1})\epsilon(f)= \\
 &\left[(1\acl  u\ps{2})(S^{-1}(u\ps{1})\rt f\acl 1)\right]\ot_\Fc 1= f\acl
 u\ot_\Fc 1.
\end{align*}
Finally one checks the comultiplicativity as follows:
\begin{align*}
&\natural\ot\natural (\D_\Cc(f\acl u\ot_\Fc 1))=\\
& \natural(f\ps{1}\acl u\ps{1}\ot_\Fc 1)\ot
\natural(u\ps{2}\ns{-1}f\ps{2}\acl u\ps{2}\ns{0}\ot_\Fc 1)=\\
&\epsilon(f\ps{2})u\ps{1}\ot
\epsilon(u\ps{2}\ns{-1}f\ps{2})u\ps{2}\ns{0}=\\
&\epsilon(f\ps{1})u\ps{1}\ot
\epsilon(f\ps{2})u\ps{2}=\D_\Uc(\natural(f\acl u\ot_\Fc 1)).
\end{align*}
\end{proof}

By Proposition \ref{phi}, the following is a map of $\Hc$-module
coalgebras.
\begin{align*}
&\Phi:\Hc\ra \Cc\cl \Fc, \\
&\Phi(h)= \dot h\ps{1}\cl h\ps{2}1.
\end{align*}
To be able to use the Theorem  \ref{quasi}, we need to show that
$\Phi$ is bijective. One has
\begin{align*}
&\Phi( f\acl u )= \\
&\wid{(f\ps{1}\acl u\ns{0})}\ot f\ps{2}u\ns{1}=\\
&\wid{ (1\acl (u\ns{0}\ps{1}))(S(u\ns{0}\ps{2})\rt f\ps{1}\acl 1)
)}\ot f\ps{2}u\ns{1}=\\
&\wid{1\acl u\ns{0}}\ot fu\ns{1}\equiv u\ns{0}\ot fu\ns{1} ,
\end{align*}
from which it follows that $\Phi^{-1}(u\ot f)=fS^{-1}(u\ns{1})\acl u\ns{0}$.

We are now in a position to apply Theorem \ref{quasi}
to get  a bicocyclic module for
computing the Hopf cyclic cohomology of $\Hc$. We would like,
before that, to understand  the action and
 coaction of $\Hc$ on $\Cc$. One has
\begin{align*}
& (f\acl u)\cdot \wid{(g\acl v)}= \wid {fu\ps{1}\rt g\acl
u\ps{2}v}\equiv \ve(f)\ve(g)uv,
\end{align*}
 which coincides with the natural action of $\Hc$ on $\Uc$.

Next,  we take up the coaction of $\Cc$.
Recall the  coaction of $\Hc$ on $\Cc$,
\begin{align*}
&\Db:\Cc\ra \Hc\ot  \Cc, \\\notag
 &\Db(\wid{1\acl u})= (1\acl u)\ps{1}S(1\acl u)\ps{3} \ot \wid{(1\acl
 u)\ps{2}}\, .
\end{align*}
\begin{proposition}
The above coaction coincides on the original
coaction of $\Fc$ on $\Uc$.
\end{proposition}
\begin{proof} Using the fact that $\Uc$ is cocommutative, one has
\begin{align*}
&   (1\acl u)\ps{1} S((1\acl u)\ps{3})\ot \wid{ (1\acl u)\ps{2}}=\\
&  (1\acl u\ps{1}\ns{0})S(u\ps{1}\ns{2}u\ps{2}\ns{1}\acl u\ps{3})\ot\wid{ u\ps{1}\ns{1}\acl u\ps{2}\ns{0}}=\\
& (1\acl u\ps{1}\ns{0})(1\acl S(u\ps{3}\ns{0})(S(u\ps{1}\ns{2}u\ps{2}\ns{1}u\ps{3}\ns{1})\acl 1)\ot \\
&\ot\wid{ u\ps{1}\ns{1}\acl u\ps{2}\ns{0}}=\\
&(1\acl u\ps{1}\ns{0}S(u\ps{3}\ns{0}))(S(u\ps{1}\ns{2}u\ps{2}\ns{1}u\ps{3}\ns{1})\acl 1)\ot \wid{u\ps{1}\ns{1}\acl u\ps{2}\ns{0}}=\\
\end{align*}
which means that after identifying $\Cc$ with $\Uc$  one has
\begin{align*}
 & (1\acl u\ps{1}\ns{0}S(u\ps{3}\ns{0}))(S(u\ps{1}\ns{1}u\ps{2}\ns{1}u\ps{3}\ns{1})\acl 1)\ot u\ps{2}\ns{0}
 =\\
 & (1\acl u\ns{0}\ps{1}S(u\ns{0}\ps{3}))(S(u\ns{1})\acl 1)\ot u\ns{0}\ps{2}=\\
 & 1\acl S(u\ns{1})\ot u\ns{0}\equiv S(u\ns{1})\ot u\ns{0}.
\end{align*}
\end{proof}

With the above identifications of actions and coactions,  one has
the following bicyclic module $\FZ(\Hc, \Fc; \Cb_\d)$, where $\d$ is
the modular character:

\begin{align}\label{UF}
\begin{xy} \xymatrix{  \vdots\ar@<.3 ex>[d]^{B_\Uc} & \vdots\ar@<.3 ex>[d]^{B_\Uc}
 &\vdots \ar@<.3 ex>[d]^{B_\Uc} & &\\
\Cb_\d \ot \Uc^{\ot 2} \ar@<.3 ex>[r]^{b_{\Fc}}\ar@<.3
ex>[u]^{b_\Uc} \ar@<.3 ex>[d]^{B_\Uc}&
  \Cb_\d\ot \Fc\ot \Uc^{\ot 2}  \ar@<.3 ex>[r]^{b_{\Fc}}\ar@<.3 ex>[l]^{B_{\Fc}}\ar@<.3 ex>[u]^{b_\Uc}
   \ar@<.3 ex>[d]^{B_\Uc}&\Cb_\d\ot\Fc^{\ot 2} \ot \Uc^{\ot 2}
   \ar@<.3 ex>[r]^{b_{\Fc}}\ar@<.3 ex>[l]^{B_{\Fc}}\ar@<.3 ex>[u]^{b_\Uc}
   \ar@<.3 ex>[d]^{B_\Uc}&\ar@<.3 ex>[l]^{B_{\Fc}} \hdots&\\
\Cb_\d \ot \Uc \ar@<.3 ex>[r]^{b_{\Fc}}\ar@<.3 ex>[u]^{b_\Uc}
 \ar@<.3 ex>[d]^{B_\Uc}&  \Cb_\d\ot\Fc\ot \Uc \ar@<.3 ex>[r]^{b_{\Fc}}
 \ar@<.3 ex>[l]^{B_{\Fc}}\ar@<.3 ex>[u]^{b_\Uc} \ar@<.3 ex>[d]^{B_\Uc}
 &\Cb_\d\ot \Fc^{\ot 2} \ot \Uc \ar@<.3 ex>[r]^{b_{\Fc}}\ar@<.3 ex>[l]^{B_{\Fc}}\ar@<.3 ex>[u]^{b_\Uc}
  \ar@<.3 ex>[d]^{B_\Uc}&\ar@<.3 ex>[l]^{B_{\Fc}} \hdots&\\
\Cb_\d \ot \Cb \ar@<.3 ex>[r]^{b_{\Fc}}\ar@<.3 ex>[u]^{b_\Uc}&
\Cb_\d\ot\Fc\ot \Cb \ar@<.3 ex>[r]^{b_{\Fc}}\ar[l]^{B_{\Fc}}\ar@<.3
ex>[u]^{b_\Uc}&\Cb_\d\ot\Fc^{\ot 2} \ot \Cb \ar@<.3
ex>[r]^{b_{\Fc}}\ar@<.3 ex>[l]^{B_{\Fc}}\ar@<1 ex >[u]^{b_\Uc}
&\ar@<.3 ex>[l]^{B_{\Fc}} \hdots&,  }
\end{xy}
\end{align}

At this stage, we introduce the following Chevalley-Eilenberg-type
bicomplex:

\begin{align}\label{UF+}
\begin{xy} \xymatrix{  \vdots\ar[d]^{\p_\Fg} & \vdots\ar[d]^{\p_\Fg}
 &\vdots \ar[d]^{\p_\Fg} & &\\
\Cb_\d \ot \wg^2\Fg \ar[r]^{\b_\Fc} \ar[d]^{\p_\Fg}&
 \Cb_\d\ot\Fc\ot \wg^2\Fg  \ar[r]^{\b_\Fc} \ar[d]^{\p_\Fg}
 &\Cb_\d\ot\Fc^{\ot 2} \ot \wg^2\Fg \ar[r]^{~~~~~~~~~\b_\Fc}   \ar[d]^{\p_\Fg}& \hdots&\\
\Cb_\d \ot \Fg \ar[r]^{\b_\Fc} \ar[d]^{\p_\Fg}&  \Cb_\d\ot\Fc\ot \Fg
\ar[r]^{\b_\Fc} \ar[d]^{\p_\Fg}&\Cb_\d\ot\Fc^{\ot 2}
 \ot \Fg  \ar[d]^{\p_\Fg} \ar[r]^{~~~~~~~~~\b_\Fc}& \hdots&\\
\Cb_\d \ot \Cb \ar[r]^{\b_\Fc}&  \Cb_\d\ot\Fc\ot \Cb
\ar[r]^{\b_\Fc}&\Cb_\d\ot\Fc^{\ot 2} \ot \Cb \ar[r]^{~~~~~~~\b_\Fc}
& \hdots&,  }
\end{xy}
\end{align}

Here $\p_\Fg$ is the Lie algebra homology boundary  of $\Fg$ with
coefficients in right module $\Cb_\d\ot \Fc^{\ot p}$; explicitly
\begin{align*}
&\p_\Fg(\one \ot \td f\ot X^0\wdots X^{q-1} )= \\
&\sum_i(-1)^i (\one \ot \td f)\lt X^i\ot X^0\wdots
\widehat{X^i}\wdots X^{q-1}+ \\
&(-1)^{i+j}\sum_{i<j} \one \ot \td f\ot [X_i,X_j]\wg X^0\wdots
\widehat{X^i}\wdots \widehat{X^j}\wdots X^{q-1},
\end{align*}
where
\begin{equation}\label{g-action}
(\one \ot \td f)\lt X= \delta(X)\ot \td f+ \one \ot S(1\acl X)\td f.
\end{equation}
The horizontal operator is given by
\begin{align*}
&\b_\Fc(\one \ot \td f\ot X^1\wdots X^q)= \\
&\sum_{i=0}^{q} \one \ot (\Id^{\ot i }\ot\D\ot \Id^{\ot (q-i)})(\td
f)\ot
X^1\wdots X^q +\\
&(-1)^{q+1} \one \ot \td f\ot S(X^1\ns{1}\dots X^q\ns{1})\ot
X^1\ns{0}\wdots X^q\ns{0}
\end{align*}

In other words,  $\b_\Fc$ is just the coalgebra cohomology
coboundary of the coalgebra  $\Fc$ with coefficients in $\wg^q\Fg$
induced from the coaction of $\Fc$ on $\Uc^{\ot q}$.
\begin{align}\label{wgcoaction}
&\Db_{\wg\Fg}:\Cb_\d\ot \wg^q\Fg\ra \Fc\ot
\Cb_\d\ot\wg^q\Fg,\\\notag &\Db(\one \ot X^1\wdots
X^q)=S(X^1\ns{1}\dots X^q\ns{1})\ot \one\ot X^1\ns{0}\wdots
X^q\ns{0}
\end{align}

 One notes that, since $\Fc$
is commutative, the coaction $\Db_{\wg \Fg}$ is well-defined.

\begin{proposition}\label{simple}
The bicomplexes \eqref{UF} and $\eqref{UF+}$ have quasi-isomorphic
total complexes.
\end{proposition}
\begin{proof}
We apply  antisymmetrization map
\begin{align*} &\td \a:\Cb_\d\ot\Fc^{\ot p}\ot
\wg^q\Fg\ra \Cb_\d\ot \Fc^{\ot p}\ot \Uc^{\ot q}\\
&\td\a= \Id \ot \a ,
\end{align*}
where $\a$ is the usual antisymmetrization map
\begin{equation}
\a(X^1\wdots X^p)= \frac{1}{p!} \sum_{\s\in S_p}(-1)^\s  X^{\s(1)}\odots X^{
\s(p)}.
\end{equation}
Its left inverse is $\td \mu:= \Id \ot \mu$ where $\mu$ is the
natural left inverse to $\a$  (see \eg \cite[page 436]{ka}). As in
the proof of Proposition 7 of \cite{cm2}, one has the following
commutative diagram:
\begin{equation*}
\begin{xy} \xymatrix{
\Cb_\d\ot\Fc^{\ot p}\ot \wg^q\Fg\ar[r]^{\td \a}\ar@<1 ex>[d]^{0}&
\Cb_\d\ot
\Fc^{\ot p}\ot \Uc^{\ot q}\ar@<1 ex>[d]^{b_\Uc}\\
\Cb_\d\ot\Fc^{\ot p}\ot \wg^{q+1}\Fg\ar[r]^{\td \a}\ar@<1
ex>[u]^{\p_\Fg}& \Cb_\d\ot \Fc^{\ot p}\ot \Uc^{\ot q+1}\ar@<1
ex>[u]^{B_\Uc}. }
\end{xy}
\end{equation*}
Since $\td\a$ does not affect $\Fc^{\ot p}$, it is easy to see that
the following diagram also commutes:
\begin{equation*}
\begin{xy} \xymatrix{
\Cb_\d\ot\Fc^{\ot p}\ot \wg^q\Fg\ar[r]^{\td \a}\ar@<1
ex>[d]^{\b_\Fc}& \Cb_\d\ot
\Fc^{\ot p}\ot \Uc^{\ot q}\ar@<1 ex>[d]^{b_\Fc}\\
\Cb_\d\ot\Fc^{\ot p+1}\ot \wg^{q}\Fg\ar[r]^{\td \a}\ar@<1
ex>[u]^{\Bc_\Fc}& \Cb_\d\ot \Fc^{\ot p+1}\ot \Uc^{\ot q}\ar@<1
ex>[u]^{B_\Fc}. }
\end{xy}
\end{equation*}

Here $\Bc_\Fc$ is the Connes boundary operator for the cocyclic
module $C^n_\Fc(\Fc, \Cb_\d\ot\wg^q\Fg)$, where $\Fc$ coacts on
$\Cb_\d\ot\wg^q\Fg$ via \eqref{wgcoaction} and acts trivially.

 Since $\Fc$ is commutative and its action on $\Cb_\d\ot \wg^q\Fg$ is trivial,
  Theorem 3.22 of \cite{kr1} implies  that  $\Bc_\Fc\equiv 0$
  in Hochschild cohomology. On the other hand Theorem 7 \cite{cm2}  (which is
  applied here for coefficients in a general module),
   proves that the columns  of \eqref{UF} with coboundary $B_\Uc$ and
   the columns of \eqref{UF+}
  with coboundary $\p_\Fg$ are quasi-isomorphic. This finishes the
  proof.
\end{proof}


\subsection{Applications}

\subsubsection{Hopf cyclic Chern classes}\label{uchern}
 In this section we compute the relative Hopf cyclic cohomology of $\Hc_n$ modulo
the subalgebra  $\Lc:= U(\Fg \Fl_n)$ with coefficients in $\Cb_\d$.

To this end, we first form the Hopf subalgebra  $\Kc:=
\Lc \acl \Fc\subset \Hc$. Here   $\Lc$ acts on $\Fc$ via its action
inherited from $\Lc\subset \Uc$ on $\Fc$, and  $\Fc$ coacts on $\Lc$
trivially. The second coalgebra is   $\Cc:=\Hc\ot_\Kc \Cb$.
 Letting $\Fh:= \Fg \Fl_n$, and
 $\Sc:=S(\Fg/\Fh)$ be the symmetric algebra of the vector space
  $\Fg/\Fh$,
we identify   $\Cc$ with the coalgebra $\Sc$, as follows.
Since $\Fg\cong V\al\Fh$, where $V= \Rb^n$,
 we can regard $\Uc$ as being $ \Uc(V)\al
  \Uc(\Fh)$. The identification of  $\Cc$ with $\Sc$ is achieved by the map
\begin{align*}
(X\al Y)\acl f\ot_\Kc 1 \mapsto \epsilon(Y)\epsilon(f)X ,
\qquad X\in \Uc(V), \, Y\in U(\Fh) , \, f\in \Fc .
\end{align*}
As in Lemma \ref{natural}, one checks that this identification is
  a coalgebra isomorphism. Finally,  since $\Fg/\Fh\cong V$ as vector spaces,
  $\Uc(V)\cong S(\Fg/\Fh)$ as
  coalgebras.

  Similarly, one identifies, $\KL:=\Kc\ot_\Lc \Cb$ with $\Fc$
  as coalgebras, via
  \begin{align*}
Y\acl f\ot_\Lc 1 \mapsto\epsilon(Y)f,
\qquad Y\in U(\Fh) , \, f\in \Fc .
\end{align*}

  Let us recall the mixed complex  $(C^n, b,B)$, where  $C^n:=C^n(\Hc,\Lc; \Cb_\d)= \Cb_\d\ot_\Lc\Cc^{\ot n}$, which computes relative
  Hopf cyclic cohomology of $\Hc$ modulo $\Lc$ with coefficients  in $\Cb_\d$.
   We also recall that  the isomorphism  $\Theta$ defined in \eqref{teta}
     identifies  this  complex with the diagonal complex $\FZ^{\ast,\ast}(\Hc, \Kc, \Lc;\Cb_\d)=
       \Cb_\d\ot_\Lc \Fc^{\ot p}\ot \Sc^{\ot q}$. Now by
       using Theorem \ref{quasi} we get the following bicomplex
         whose total complex is quasi isomorphic to $C^n$ via the
       Eilenberg-Zilber Theorem.

\begin{align}\label{VF}
\begin{xy} \xymatrix{  \vdots\ar@<.3 ex>[d]^{B_\Sc} & \vdots\ar@<.3 ex>[d]^{B_\Sc}   &\vdots \ar@<.3 ex>[d]^{B_\Sc} & &\\
\Cb_\d \ot_\Lc \Sc^{\ot 2} \ar@<.3 ex>[r]^{b_{\Fc}}\ar@<.3
ex>[u]^{b_\Sc} \ar@<.3 ex>[d]^{B_\Sc}& \Cb_\d\ot_\Lc \Fc\ot\Sc^{\ot
2} \ar@<.3 ex>[r]^{b_{\Fc}}\ar@<.3 ex>[l]^{B_{\Fc}}\ar@<.3
ex>[u]^{b_\Sc} \ar@<.3 ex>[d]^{B_\Sc}&\Cb_\d\ot_\Lc \Fc^{\ot 2}
 \ot\Sc^{\ot 2} \ar@<.3 ex>[r]^{~~~~~~~~b_{\Fc}}\ar@<.3 ex>[l]^{B_{\Fc}}\ar@<.3 ex>[u]^{b_\Sc}
 \ar@<.3 ex>[d]^{B_\Sc}&\ar@<.3 ex>[l]^{~~~~~~~~~B_{\Fc}} \hdots&\\
\Cb_\d \ot_\Lc \Sc \ar@<.3 ex>[r]^{b_{\Fc}}\ar@<.3 ex>[u]^{b_\Sc}
\ar@<.3 ex>[d]^{B_\Sc}&
  \Cb_\d\ot_\Lc \Fc\ot\Sc \ar@<.3 ex>[r]^{b_{\Fc}}\ar@<.3 ex>[l]^{B_{\Fc}}\ar@<.3 ex>[u]^{b_\Sc}
   \ar@<.3 ex>[d]^{B_\Sc}&\Cb_\d\ot_\Lc \Fc^{\ot 2}\ot \Sc
    \ar@<.3 ex>[r]^{~~~~~~~~~~~b_{\Fc}}\ar@<.3 ex>[l]^{B_{\Fc}}\ar@<.3 ex>[u]^{b_\Sc}
    \ar@<.3 ex>[d]^{B_\Sc}&\ar@<.3 ex>[l]^{~~~~~~~~~B_{\Fc}} \hdots&\\
\Cb_\d \ot_\Lc\Cb \ar@<.3 ex>[r]^{b_{\Fc}}\ar@<.3 ex>[u]^{b_\Sc}&
\Cb_\d\ot_\Lc\Fc\ot\Cb
 \ar@<.3 ex>[r]^{b_{\Fc}}\ar[l]^{B_{\Fc}}\ar@<.3 ex>[u]^{b_\Sc}&\Cb_\d\ot_\Lc \Fc^{\ot 2}
 \ot \Cb \ar@<.3 ex>[r]^{~~~~~~~~~~b_{\Fc}}\ar@<.3 ex>[l]^{B_{\Fc}}\ar@<1 ex >[u]^{b_\Sc}  &\ar@<.3 ex>[l]^{~~~~~~~~~~B_{\Fc}} \hdots&,  }
\end{xy}
\end{align}

Similar to \eqref{UF+},  we introduce the following bicomplex.

\begin{align}\label{VF+}
\begin{xy} \xymatrix{  \vdots\ar[d]^{\p_{(\Fg,\Fh)}} & \vdots\ar[d]^{\p_{(\Fg,\Fh)}}
 &\vdots \ar[d]^{\p_{(\Fg,\Fh)}} & &\\
\Cb_\d \ot \wg^2V \ar[r]^{\b_\Fc} \ar[d]^{\p_{(\Fg,\Fh)}}&
 \Cb_\d\ot\Fc\ot \wg^2V  \ar[r]^{\b_\Fc} \ar[d]^{\p_{(\Fg,\Fh)}}
 &\Cb_\d\ot\Fc^{\ot 2} \ot \wg^2V \ar[r]^{~~~~~~~~~\b_\Fc}   \ar[d]^{\p_{(\Fg,\Fh)}}& \hdots&\\
\Cb_\d \ot V \ar[r]^{\b_\Fc} \ar[d]^{\p_{(\Fg,\Fh)}}&
\Cb_\d\ot\Fc\ot V \ar[r]^{\b_\Fc}
\ar[d]^{\p_{(\Fg,\Fh)}}&\Cb_\d\ot\Fc^{\ot 2}
 \ot V  \ar[d]^{\p_{(\Fg,\Fh)}} \ar[r]^{~~~~~~~~~\b_\Fc}& \hdots&\\
\Cb_\d \ot \Cb \ar[r]^{\b_\Fc}&  \Cb_\d\ot\Fc\ot \Cb
\ar[r]^{\b_\Fc}&\Cb_\d\ot\Fc^{\ot 2} \ot \Cb \ar[r]^{~~~~~~~\b_\Fc}
& \hdots&,  }
\end{xy}
\end{align}

Here $V$ stands for the vector space $\Fg/\Fh$ and  $\p_{(\Fg,
\Fh)}$ is the relative Lie algebra homology boundary of $\Fg$
relative to $\Fh$ with coefficients in $\Cb_\d\ot \Fc^{\ot p}$ with
the action defined in \eqref{g-action}. The coboundary $\td\b_\Fc$
is induced by $\b_\Fc$ the coalgebra cohomology coboundary of $\Fc$
with trivial coefficients in $\Cb_\d\ot \wg V$. One notes that the
coaction of $\Fc$ on $\Sc$ is induced from the coaction of $\Fc$ on
$\Uc$ via the natural projection $\pi:\Uc\ra \Uc\ot_\Lc \Cb\cong
\Sc$.
\begin{proposition}
The total complexes of the  bicomplexes  \eqref{VF} and \eqref{VF+}
are quasi isomorphic.
\end{proposition}
\begin{proof}
The proof is similar to that of  Proposition \ref{simple} but more
delicate. One  replaces the anitsymmetrization map by its relative
version and instead of Proposition 7 \cite{cm2} one uses Theorem 15
\cite{cm6}. On the other hand since $\Lc\subset \Hc_n$ is
cocommutative, every map in the bicomplex \eqref{VF} is
$\Lc$-linear, including the homotopy maps used in Eilenberg-Zilber
theorem~\cite{kr3}, and in~\cite[pp. 438-442]{ka}, as well as in the
spectral sequence used in \cite[Theorem 3.22]{kr1}.
\end{proof}

To find the cohomology of the above bicomplex we first compute the
subcomplex consisting of $\Fg \Fl_n$-invariants.

\begin{definition}
 We say a map $\Fc^{\ot m}\ra \Fc^{\ot p}$ is an ordered
projection of order $(i_1,\dots, i_q)$, including $\emptyset$ order,
if it of the form
\begin{equation} f^1\odots f^n\mapsto f^1\dots f^{i_1-1}\ot f^{i_1}\dots
f^{i_2-1}\odots f^{i_p}\dots
 f^m.
 \end{equation}
 We denote the set of all such projections by $\Pi^m_q$.
\end{definition}
Let  $S_m$ denote the symmetric group  of order $m$. Let also  fix
$\s \in S_{n-q}$, and $\pi\in \Pi^{n-q}_{p}$, where  $0\le p,q\le
n$. We define
\begin{align} \label{invariants}
&\theta(\s,\pi)= \\\notag
&\sum  (-1)^\mu \one  \ot \pi(\eta^{j_1}_{
\mu(1),j_{\s(1)}} \odots \eta^{j_{n-q}}_{\mu(n-q),j_{\s(n-q)}})\ot
X_{\mu(n-q+1)}\wdots
 X_{\mu(n)}
\end{align}
where the summation is over  all $\, \mu \in S_n$ and all $ \, 1\le
j_1,j_2,\dots j_q\le n$.
\medskip

\begin{lemma}\label{glninv}
The elements $\, \theta(\s,\pi)\in \Cb_\d\ot \Fc^{\ot q}\ot\wg^p V$,
with $\, \s\in S_{n-p}$ and $\, \pi\in \Pi^{n-p}_q$, are
$\Fh$-invariant and span the space $\Cb_\d\ot_\Fh \Fc^{\ot q}\ot
\wg^p V$.
\end{lemma}

\proof We identify $\Cb_\d$ with $\wg^n V^\ast $ as $\Fh$-modules, by
sending $\one \in \Cb_\d$ into the volume element
 $\xi^1\wg \ldots  \wg \xi^n \in \wg^n V^\ast $, where
$\{\xi^1, \ldots , \xi^n \}$
 is the dual basis to  $\{X_1, \ldots , X_n \}$. One obtains thus an isomorphism
 \begin{equation}\label{invariant}
\Cb_\d\ot_\Fh \Fc^{\ot p} \ot  \wg^q V\cong \big(\wg^n V^\ast\ot
\Fc^{\ot p} \ot \wg^q V \big)^\Fh.
\end{equation}
We now observe that the relations \eqref{u>f} and \eqref{Yhighd} show
that the action of $\Fh$ on the jet coordinates $\eta^i_{j j_1,\dots,j_k}$s is
tensorial. In particular, the central element of $\Fh$,
  $Z=\sum_{i=1}^n Y^i_i$, acts as a grading operator, assigning degree
  $1$ to each $X \in V$, degree $-1$ to each $\xi \in V^\ast$ and
  degree $k$ to $\eta^i_{j j_1,\dots,j_k}$. As an immediate consequence,
  the space of invariants \eqref{invariant} is seen to
  be  generated by monomials with the same number of upper and lower indices,
  more precisely
   \begin{equation}\label{invariant2}
\big(\wg^n V^\ast\ot \Fc^{\ot p}  \ot  \wg^q V \big)^\Fh =
\big(\wg^n V^\ast \ot \Fc^{\ot p}[n-q]\ot  \wg^q V \big)^\Fh ,
\end{equation}
 where $\Fc^{\ot p}[n-q] $ designates the homogeneous component of
$\Fc^{\ot q}$ of degree $n-q$.

  Furthermore, using \eqref{aeta3}, \eqref{aeta4}
  we can think of
  $\Fc$ as generated  by $\a^i_{jj_1,\dots,j_k}$s instead of the $\eta^i_{j_1,\dots,j_k}$s.
The advantage is that the $\a^i_{ j_1,\dots,j_k}$s are symmetric in
  all lower indices and freely generate the algebra $\Fc$. Fixing a lexicographic
  ordering on the set of multi-indices $ ^i_{j j_1,\dots,j_k}$,
  one obtains a PBW-basis of $\Fc$, viewed as the
  Lie algebra generated by the $\a^i_{jj_1,\dots,j_k}$s.  This allows us to extend
  the assignment
  \begin{align*}
  \a^i_{j j_1,\dots,j_k} \, \mapsto \, \xi^i\ot X_j \ot X_{j_1} \ot \ldots \ot X_{j_k} \in V^\ast
  \ot V^{\ot k+1 },
  \end{align*}
 first to the products which define the PBW-basis and next
  to an  $\Fh$-equivariant
  embedding of $\Fc^{\ot q}[n-q] $ into a finite direct sum of
  tensor products of the form $(V^\ast)^{\ot r} \ot V^{\ot s} $.
   Sending the wedge product of vectors into the antisymmetrized
   tensor product, we
    thus obtain an embedding
    \begin{equation}\label{invariant3}
\big(\wg^n V^\ast \ot \Fc^{\ot p}[n-q]\ot  \wg^q V \big)^\Fh
\hookrightarrow \big((V^\ast)^{\ot n} \ot \sum_{r, s} (V^\ast)^{\ot
r} \ot V^{\ot s} \ot V^{\ot q} \big)^\Fh.
\end{equation}
  In the right hand side we are now dealing with the classical theory of
  tensor invariants
  for $\GL(n, \Rb)$, \cf~\cite{weyl}. All such invariants are linear combinations
  tensors of the form
     \begin{align*}
& T_\s = \sum
 \xi^{j_1} \ot \ldots \ot \xi^{j_n} \ot  X_{j_{\s(1)}} \ot \ldots \ot X_{j_{\s(p)}} \ot \ldots \\
 & \quad \dots
 \ot \xi^{j_{r}} \ot \ldots \ot \xi^{j_m} \ot \ldots \ot
 X_{j_{\s(s)}} \ot \ldots \ot X_{j_{\s(m)}},
 \end{align*}
where  $\,  \s \in S_m $ and the sum is over  all  $ \, 1\le
j_1,j_2,\dots j_m \le n$.

Because the antisymmetrization is a projection operator, such linear
combinations belong to the subspace in the left hand side of
\eqref{invariant3} only when they are totally antisymmetric in the
first $n$-covectors. Recalling the identification $\, \one \cong
\xi^1\wg \ldots  \wg \xi^n $, and the fact that $\a^i_{j k} =
\eta^i_{j k}$, while the higher $\a^i_{ j j_1,\dots,j_k}$s are
 symmetric in all lower indices, one easily sees that the projected invariants
 belong to the linear span of those of the form \eqref{invariants}.
\endproof
\medskip

For each partition $\lb =  (\lb_1 \geq  \ldots \geq  \lb_k)$
 of the set $ \{ 1, \dots , p \}$, where $\, 1 \leq p \leq n$,
 we let $\lb \in S_p$ also denote a permutation whose
 cycles have lengths $\lb_1 \geq  \ldots \geq  \lb_k$,  \ie
 representing
 the corresponding conjugacy class $[\lb] \in [S_p] $.
 We then define
\begin{equation*}
C_{p, \lb} := \sum  (-1)^\mu \one \ot \eta^{j_1}_{ \mu(1),j_{{\lb
(1)}}} \wdots \eta^{j_{p}}_{\mu(p),j_{{\lb (p)}}}\ot
X_{\mu(p+1)}\wdots
 X_{\mu(n)} ,
\end{equation*}
where the summation is over  all $\, \mu \in S_n$ and all $ \, 1\le j_1,j_2,\dots , j_p\le n$.

\begin{theorem} \label{rclasses}
 The cochains $\{ C_{p, \lb} \, ; \, 1 \leq p \leq n , \quad [\lb] \in [S_p] \}$
 are cocycles and their classes form a basis of
the  group  $HP^\e (\Hc_n , \Uc(\Fg \Fl_n); \Cb_{\d})$, where $\e
\equiv n$ {\em mod} $2$. The complementary parity group
 $HP^{1-\e} (\Hc_n , \Uc(\Fg \Fl_n); \Cb_{\d}) = 0$.
\end{theorem}

\begin{proof}
Let $\Fd$ be the commutative Lie algebra generated by all
$\eta^i_{j,k}$, and  let $F=U(\Fd)$, be the polynomial algebra of
$\eta^i_{j,k}$. Obviously $F$ is a Hopf subalgebra of $\Fc$ and
stable under the action of  $\Fh$. By using Lemma \ref{invariants},
we have:
\begin{equation}\label{u(d)}
\Cb_\d\ot_\Fh  F^{\ot p}\ot \wg^q V\cong\Cb_\d\ot_\Fh \Fc^{\ot p}\ot
\wg^q V\,.
\end{equation}
Hence   we reduce the problem to compute the cohomology of $(
\Cb_d\ot_\Fh F^{\ot \ast}\ot  \wg^q V, b_F)$, where $\b_F$ is
induced by  the Hochschild coboundary of the coalgebra $F$ with
trivial coefficients. One uses the fact that $ (\Cb_\d\ot F^{\ot
\ast}\ot  \wg^q V, b_F)$ and $ (\Cb_\d V\ot \wg^\ast \Fd \ot\wg^q,
0)$ are homotopy equivalent and the homotopy is $\Fh$-linear
\cite{ka}  to see that the $q$th cohomology of the complex under
consideration is $\Cb_\d\ot_\Fh \wg^p \Fd\ot\wg^q V$. Similar to the
proof of Lemma \ref{invariants} we replace $\Cb_\d$ with
$\wg^nV^\ast$. We also replace $\a^i_{j,k}$ with $\td\x^i\ot X'_j\ot
X_k'' $, where $\td\x^i$, $X_j'$ and $X_k''$ are basis for $V^\ast$,
$V$ and $V$ respectively.  Using the above identification one has
$\Cb_\d\ot_h  \wg^p\Fd\ot \wg^q V \cong \big(\wg^n V^\ast\ot
(\wg^pV^\ast\ot \wg^pV\ot\wg^pV) \ot \wg^q V\big)^\Fh$. Via the same
argument as in the proof of Lemma \ref{invariants}, one concludes
that the invariant space is generated by elements of the form
\begin{equation} \label{chern}
C_{p, \s} := \sum  (-1)^\mu \one \ot \eta^{j_1}_{ \mu(1),j_{\s(1)}}
\wdots \eta^{j_{p}}_{\mu(p),j_{\s(p)}}\ot X_{\mu(p+1)}\wdots
 X_{\mu(n)}
\end{equation}
where $\s\in S_p$,  and  the summation is taken as in
\eqref{invariants}. Now let $\s=\s_1\dots\s_k$, where $\s_1=(a,
\s(a), \dots, \s^{\a-1}(a))$, \dots, $\s_k=(z, \s(z), \dots,
\s^{\z-1}(z)).$\,\, Since distinct cyclic permutations commute among
each other we may assume $\a\ge \b\dots\ge \g\ge\z$. We defined the
following two permutations.
 \[\tau:=(1,2,\dots, \a)(\a+1, \dots, \b),\dots,
(\a+\dots+\g+1, \dots, \a+\dots +\z),\]
\[\theta(i)=\left\{
              \begin{array}{ll}
               \s^{i-1}(a), &1\le i\le  \a \\
                \s^{i-1-\a}(b), & 1+\a\le i\le \a+\b \\
                \vdots  & \vdots  \\
                \s^{i-1-\a-\dots-\g}(z), &\a+\dots+\g+1\le i\le
                p\,\\
                     i, &p+1\le i\le n
              \end{array}
            \right.
\]
 We claim
$C_{p,\s}$ and $C_{p,\tau}$ coincide up to a sign. Indeed,
\begin{align*}
&C_{p,\s}= \pm\sum_{\mu}(-1)^\mu\eta^{j_a}_{\mu(a), j_\s(a)}\wg
\eta^{j_{\s(a)}}_{\mu(\s(a)),j_{\s^2(a)}}\wg\dots\wg
\eta^{j_{\s^{\a-1}(a)}}_{\mu(\s^{\a-1}(a)),j_a}\wg\\
&\eta^{j_b}_{\mu(b), j_\s(b)}\wg
\eta^{j_{\s(b)}}_{\mu(\s(b)),j_{\s^2(b)}}\wg\dots\wg
\eta^{j_{\s^{\b-1}(b)}}_{\mu(\s^{\b-1}(b)),j_b}\wdots \\
&\eta^{j_z}_{\mu(z), j_\s(z)}\wg
\eta^{j_{\s(z)}}_{\mu(\s(z)),j_{\s^2(z)}}\wg\dots\wg
\eta^{j_{\s^{\z-1}(z)}}_{\mu(\s^{\z-1}(z)),j_z}\ot X(\mu)=\\
&\pm\sum_{\mu} (-1)^\mu\eta^{j_1}_{\mu(a), j_2}\wg
\eta^{j_2}_{\mu(\s(a)),j_3}\wg\dots\wg
\eta^{j_\a}_{\mu(\s^\a(a)),j_1}\wg\\
&\eta^{j_{\a+1}}_{\mu(b), j_{\a+2}}\wg
\eta^{j_{\a+2}}_{\mu(\s(b)),j_{\a+3}}\wg\dots\wg
\eta^{j_{\a+\b}}_{\mu(\s^\b(b)),j_{\a+1}}\wdots \\
&\eta^{j_{\a+\dots+\g+1}}_{\mu(z), j_{\a+\dots+\g+2}}\wg \eta^{j_{
\a+\dots+\g +2}}_{\mu(\s(z)),j_{\a+\dots \g+3}}\wg\dots\wg
\eta^{j_{\a+\dots +\z}}_{\mu(\s^\z(z)),j_{\a+\dots \g+1}}\ot X(\mu)\ot  =\\
& \pm\sum_{\mu}(-1)^\mu \eta^{j_1}_{ \mu\theta(1),j_{\tau(1)}}
\wdots \eta^{j_{q}}_{\mu\theta(p),j_{\tau(q)}}\ot X(\mu\theta) = \pm
C_{p,\tau}.
\end{align*}
Here $X(\mu)$ stands for $ X_{\mu(p+1)}\wdots X_{\mu(n)}$.

The fact that $\{ C_{p, \lb} \, ; \, 1 \leq p \leq n , \quad [\lb] \in
[S_p] \}$ are linearly independent, and therefore form a basis for
$\Cb_\d\ot_\Fh \Fc^{\ot p}\ot \wg^qV$, follows from the observation
that, if $\s \in S_p$, the term
\[\sum  (-1)^\mu \one\ot \eta^{1}_{ \mu(1),{\s(1)}} \wdots
\eta^{{p}}_{\mu(p),\s(p)}\ot X_{\mu(p+1)}\wdots
 X_{\mu(n)},\] appears  in $C_{p,\tau}$ if and only if
$[\s]=[\tau]$.

These are all
the periodic cyclic classes, because all of them sit on the $n$th
skew diagonal of \eqref{VF+} and there are no other invariants.
Therefore $B\equiv \p_{\Fg,\Fh}=0$ on the invariant space.
\end{proof}

In particular
\begin{equation*}
\dim HP^\ast (\Hc_n , \Uc(\Fg \Fl_n); \Cb_{\d}) = \Fp (0) + \Fp(1) + \ldots
+ \Fp(n) ,
 \end{equation*}
 where $\Fp$ denotes the partition function,
 which is the same as the dimension
 of the truncated Chern ring
 $\Pc_n[c_1, \dots c_n]$. Moreover, as noted in the introduction,
 the assignment
  \begin{equation*}
C_{p, \lb} \mapsto c_{p, \lb} :=  c_{\lb_1} \cdots c_{ \lb_k} ,
\qquad \lb_1 +  \ldots +  \lb_k = p ;
\end{equation*}
 defines a linear isomorphism between
 $HP^\ast (\Hc_n , \Uc(\Fg \Fl_n); \Cb_{\d})$ and
 $\Pc_n[c_1, \dots c_n]$.

 \subsubsection{Non-periodized Hopf cyclic cohomology of $\Hc_1$}
In this section we apply the bicomplex \eqref{UF+} to compute
all cyclic cohomology of $\Hc_1$ with coefficients in the MPI
$(1, \d)$. Actually, invoking the equivalence \eqref{equiv},
we shall compute   $HC(\Hc_1^{\cop},\Cb_\d)$ instead.

Let us recall the presentation of the  Hopf algebra $\Hc_1$.
 As an algebra,
$\Hc_1$ is generated by $X$, $Y$, $\delta_k$,  $k\in \mathbb{N}$,
subject to  the relations
\begin{equation*}
[Y, X]=X,\qquad [Y,\delta_k]=k\delta_k,\qquad
[X,\delta_k]=\delta_{k+1}, \qquad [\delta_j, \delta_k]=0.
\end{equation*}
Its coalgebra structure is uniquely determined by
\begin{align*}
&\Delta(Y)=Y\ot 1+1\ot Y,\\
&\Delta(X)=X\ot 1+1\ot
X+\delta_1\ot Y,\\
 & \Delta(\delta_1)=\delta_1\ot
1+1\ot \delta_1, \\
&\epsilon(X)=\epsilon(Y)=\epsilon(\delta_k)=0 .
\end{align*}

By Theorem \eqref{bicrossed}, one has $\Hc_1^{\cop}\cong
\Fc\acl\Uc$. We next apply Theorem \ref{quasi}, which gives a
quasi-isomorphism between $C^\ast(\Fc\acl\Uc, \Cb_\d)$ and the total
complex of the bicomplex \eqref{UF+}.  One notes that  $Y$ grades
$\Hc_1$ by $\nobreak{[Y\,,\, X]\,=\,X}$, $[Y
\,,\,\d_k]\,=\,k\,\d_k$. Accordingly, $Y$ grades $\Uc$, and also $\Fc$ by
$$
Y\rt \eta^1_{1,\dots,1}\,=\,k\,\eta^1_{1,\dots 1} \, ,
$$
 where $k+1$ is the number of lower indices. Hence $Y$ grades the bicomplex
\eqref{UF+}. One notes that every identification, isomorphism, and
homotopy which has been used to pass from the bicomplex $C^\ast(\Fc\acl
\Uc, \Cb_\d)$ to \eqref{UF+} respects this grading (\cf also
\cite{mr}).  As a result, one can relativize the computations to
each homogeneous component. The degree of $\td f\in
\Fc^{\ot q}$ will be denoted by $\nm{\td f}$.

We next recall Gon\v carova's  results~\cite{gon} concerning the
Lie algebra cohomology of $\Fn$. Taking as basis of
$\Fn$ the vector fields
\begin{equation*}
e_i=x^{i+1}\frac{d}{dx}\in \Fn, \quad i\ge 1 ,
\end{equation*}
and denoting the dual basis by $\{e^i\in \Fn^\ast \mid i\ge 1 \}$, one identifies
$ \wg^k \Fn^*$ with the totally antisymmetric
 polynomials in variables $z_1 \dots z_k$, via the map
 \begin{equation*}
 e^{i_1}\wdots e^{i_k} \mapsto \sum_{\mu\in S_k}(-1)^\mu
 z_1^{i_{\mu(1)}} \dots z_k^{i_{\mu(k)}}.
 \end{equation*}
 According to~\cite{gon},  for each dimension $k\ge 1$,
 the Lie algebra cohomology group $H^k (\Fn , \Cb)$ is
 two-dimensional and is generated by the classes
 \begin{equation*}
 \om_k:= z_1\dots z_k\Pi_k^3, \quad  \om_k':= z_1^2\dots z_k^2\Pi_k^3,
 \end{equation*}
 where $\Pi_k:=  \underset{1\le i<j\le k}{\Pi}(z_j-z_i)$.
 We shall denote by $\xi_k \in H^k(\Fc,\Cb)$, resp.  $\xi_k' \in H^k(\Fc,\Cb)$
the image  in Hochschild cohomology with trivial coefficients of
 $\om_k$, resp. $\om_k'$,
 under the  van Est isomorphism of \cite[Proposition 7(b)]{cm2}.
 \smallskip

 \begin{lemma}\label{homotopy}
 The  map $\g: \Fc^{\ot i}\ra \Fc^{\ot (i+1)}$ defined by
  \begin{align*}
 \g(\td f)=\td f\ot \eta_1, \qquad \text{where} \quad \eta_1:=\eta^1_{11} ,
 \end{align*}
  commutes with the Hochschild coboundary
 and for $i\ge 1$ is null-homotopic.
 \end{lemma}
 \begin{proof}
We check that $\g$ is map of complexes:
 \begin{align*}
 &\g( b\td f)=b(\td f)\ot \eta_1=\\
&b(\td f)\ot \eta_1+(-1)^{i} (\td
f\ot 1\ot \eta_1 -\td f\ot 1\ot \eta_1 +\td f\ot \eta_1+1-\td f\ot \eta_1\ot 1)\\
&= b(\g(\td f)).
 \end{align*}
Using  the bicomplex \eqref{UF+}, one can
show that the above map is zero at the level of Hochschild cohomology.
 Indeed, recalling that $\Hc$ acts
 on $\Fc$ as in \eqref{newact}, and
remembering that we are dealing with the normalized bicomplex, one has
 \begin{align*}
 \p_\Fg
&\b_\Fc(\one\ot \td f \ot X )=\p_\Fg(\one\ot b(\td f) \ot X
+(-1)^{i+1}
\one\ot \td f\ot \eta_1 \ot Y=\\
& -(-1)^{i+1}\one  \ot X b(\td f)+\one \ot \td f\ot \eta_1-
\one\ot Y(\td f\ot \eta_1)=\\
&-(-1)^{i+1}\one \ot X b(\td f)-  \one\ot |\td f |\td f\ot \eta_1.
\end{align*}
On the other hand,
\begin{align*}
\b_\Fc\p_\Fg(\one\ot \td f  \ot X)=-\b_\Fc(\one \ot X \td f) .
\end{align*}
This implies that if $b(\td f)=0$  and $|\td f|\ge 1$,  then
\begin{equation*}
\td f\ot \eta_1=\frac{1}{|\td f|} b(X \td f).
\end{equation*}
 \end{proof}

To compute
  the Hochschild cohomology of $\Hc_1$, we appeal again to
 the double complex \eqref{UF+}.  The $E^{q,p}_1$ term
 is   $\Cb_\d\ot\wg^p \Fg\ot \Fc^{\ot q} $ and the boundary is $\b_\Fc$,
 as in the diagram below:

\begin{center}
$\begin{xy} \xymatrix{ \Cb_\d\ot\wg^2\Fg \ar[r]^{\b_\Fc}& \Cb_\d\ot
\Fc\ot\wg ^2\Fg  \ar[r]^{\b_\Fc}& \Cb_\d\ot \Fc^{\ot2}\ot\wg ^2\Fg\ar[r]^{~~~~~~~~~~~~~\b_\Fc} &\dots  \\
\Cb_\d\ot\Fg \ar[r]^{\b_\Fc}&  \Cb_\d\ot \Fc \ot\Fg\ar[r]^{\b_\Fc}& \Cb_\d\ot \Fc^{\ot 2}\ot\Fg \ar[r]^{~~~~~~~\b_\Fc}&\dots \\
\Cb_\d\ot\Cb \ar[r]^{\b_\Fc} &\Cb_\d\ot \Fc \ar[r]^{\b_\Fc} &
\Cb_\d\ot\Fc^{\ot 2} \ar[r]^{\b_\Fc}&\dots}
\end{xy}$
\end{center}

\begin{lemma}
\begin{align*}
 &E_2^{0,2}= \Cb[\one\ot 1\ot  X\wg Y], \\
 & E_2^{q,2}= \Cb[ \one \ot \xi_q \ot  X\wg Y ]\oplus \Cb[ \one  \ot
 \xi_q' \ot X\wg Y ],& \quad q\ge 1\\
&E_2^{0,0}= \one\ot \Cb, \quad E_2^{q,0}= \one \ot \xi_q\ot
\Cb\oplus \one \ot
 \xi_q'\ot \Cb , &\quad q\ge 1\, .\\
\end{align*}
\end{lemma}
\begin{proof}
The result for the $0$th row is obvious. For the second  row one
observes that
\begin{align*}
&\Db(\one\ot X\wg Y)=\\
&\one \ot 1 \ot  X\ot Y+ \one \ot \eta_1 \ot  Y\ot Y-\one \ot
\eta_1\ot Y\ot
Y -\one \ot 1\ot Y\ot X=\\
& \one\ot  1\ot X\wg Y.
\end{align*}
\end{proof}
\begin{lemma}
\begin{align*}
 &E_2^{0,1}= \Cb[\one\ot  1\ot Y], \\
 & E_2^{1,1}= \Cb[\one\ot \xi_1'\ot Y ]\oplus \Cb[\one\ot\xi_1 \ot X - \one\ot X\xi_1\ot Y ]\oplus \\
 &\oplus \Cb[ \one\ot \xi_1'\ot X - \frac{1}{2}\one\ot X\xi_1'\ot  Y
 ]\\
 & E_2^{p,1}= \Cb[\one\ot \xi_p\ot  Y ]\oplus \Cb[\one\ot \xi_p'\ot Y]\oplus
  \Cb[\one\ot \xi_p \ot X - \frac{1}{|\xi_p|}\one\ot X \xi_p \ot  Y
 ]\oplus\\
 &\oplus \Cb[ \one \ot \xi_p'\ot X -  \frac{1}{|\xi_p'|}\one \ot X \xi_p'\ot Y
 ], \quad \quad p\ge 2\, .\\
\end{align*}
\end{lemma}
\begin{proof}
We filter $E_1^{p,1}$ by setting
\begin{equation*} E_1^{p,1}=\one\ot \Fc^{\ot p}\ot  X
\oplus \one \ot \Fc^{\ot p}\ot Y\supset \one \ot \Fc^{\ot p}\ot
Y\supset 0.
\end{equation*}
The spectral sequence $\Ec$ that computes $E_2^{p,1}$ with respect
to the above filtration is : $\Ec_1^{-1,q}= \one\ot \Fc^{\ot q-1}
\ot X$, $\, \Ec_1^{0,q}= \one\ot \Fc^{\ot q}\ot Y$, and the rest
$=0$,
$$\begin{xy} \xymatrix{  \vdots & \vdots  &\vdots \\
 \one\ot \Fc\ot X \ar[u]^{\Id\ot b} & \one\ot \Fc^{\ot 2}\ot Y\ar[u]^{\Id\ot b}& 0\dots \\
\one\ot \Cb \ot  X \ar[u]^{\Id\ot b}& \one\ot \Fc \ot Y \ar[u]^{\Id\ot b}&0\dots \\
 0\ar[u]^0 & \one\ot \Cb \ot  Y \ar[u]^{\Id\ot b} &0\dots }
\end{xy}$$

The $\Ec_2$ term is then described by

$$\begin{xy} \xymatrix{  \vdots & \vdots  &\vdots \\
 \one\ot H^1(\Fc,\Cb)\ot X \ar[r]^{\td\g} & \one\ot H^2(\Fc, \Cb)\ot Y\ar[r]& 0\dots \\
 \one\ot  H^0(\Fc,\Cb)\ot X \ar[r]^{\td\g}& \one\ot H^1(\Fc, \Cb)\ot Y \ar[r]&0\dots \\
 0\ar[r] & \one\ot H^0(\Fc, \Cb)\ot Y \ar[r]&0\dots .}
\end{xy}$$
Here $\td\g(\one\ot\td f\ot X)= \one\ot \td f\ot\eta_1\ot Y$.  One
applies Lemma \ref{homotopy} to deduce that all maps in the above
diagram are 0, except $\begin{xy} \xymatrix{ \one\ot H^0(\Fc,\Cb)
\ot X \ar[r]^{\td\g}& \one\ot H^1(\Fc, \Cb)\ot Y}
\end{xy}$, which sends the class of  $\one\ot 1\ot X$ to the class of $\one\ot \eta_1\ot Y
=\one\ot \x_1 \ot Y$. The spectral sequence $\Ec$ obviously
collapses at this level.
\end{proof}

In turn, the spectral sequence $E$ collapses at $E_2$, and we obtain
 the following result.

\begin{proposition}
The Hochschild cohomology groups of of $\Hc_1$ are given by
\begin{align*}
&H^0(\Hc_1, \Cb)= \Cb[\one\ot 1\ot 1]\\[.2cm]
 &H^1(\Hc_1, \Cb)= \Cb[\one\ot \xi_1\ot 1]\oplus \Cb[\one\ot
\xi_1'\ot 1]\oplus
\Cb[\one\ot 1\ot Y]\\[.2cm]
 &H^2(\Hc_1, \Cb)= \Cb[\one\ot \xi_2 \ot 1]\oplus \Cb[\one\ot
\xi_2'\ot 1]\oplus \Cb[\one\ot \x_1'\ot Y]\oplus \\
&\oplus \Cb[\one\ot 1 \ot X\wg Y]\oplus \Cb[\one\ot \x_{1}\ot X-
\frac{1}{|\x_{1}|}\one\ot X \xi_{1}\ot Y] \oplus \\
&\oplus \Cb[\one\ot \x_{1}'\ot X- \frac{1}{|\x_{1}'|}\one\ot
X \xi_{p}'\ot Y]\\[.2cm]
 &H^p(\Hc_1, \Cb)= \Cb[\one\ot
\xi_p \ot 1]\oplus \Cb[\one\ot \xi_p'\ot 1]\oplus \Cb[\one\ot
\x_{p-1}\ot Y] \oplus\\
& \oplus \Cb[\one\ot \x_{p-1}'\ot Y]\oplus   \Cb[\one\ot
\x_{p-1\ot X}- \frac{1}{|\x_{p-1}|}\one\ot X \xi_{p-1}\ot Y]\oplus\\
&\oplus  \Cb[\one\ot \x_{p-1}\ot X- \frac{1}{|\x'_{p-1}|}\one\ot X
\xi'_{p-1}\ot
Y]\oplus \\
& \Cb[\one\ot \x_{p-2}\ot X\wg Y] \oplus \Cb[\one\ot \x'_{p-2}\ot
X\wg Y], \qquad\qquad p\ge 3.
\end{align*}
\end{proposition}

In the following proposition  we use the standard notation
 $S:HC^m(\Hc_1, \Cb_\d)\ra HC^{m+2}(\Hc_1,
\Cb_\d)$ for Connes' periodicity operator~\cite{cng}.
\begin{proposition}
The following classes form a basis of $HC^p (\Hc_1, \Cb_\d)$:

for $p=0,$ \quad  \quad  $\t_0 = \one \ot 1\ot 1$;

for $p=1,$ \quad \quad  $GV_1 \equiv \tau_1= \one\ot\xi_1\ot 1$,  \quad \quad
$\tau_1' = \one\ot \xi_1'\ot 1$;

for $p=2q,$ \quad
\begin{align*}
&\tau_p = \one\ot \x_p\ot 1, \quad \quad \tau_p'= \one\ot \x_p'\ot 1,  \\
&\s_p=\one\ot \x_{p-1} \ot X- \frac{1}{|\x_{p-1}|}\one\ot X
\xi_{p-1}\ot Y, \\
&\s_p'  = \one\ot \x_{p-1}'\ot X- \frac{1}{|\x_{p-1}|}\one\ot X\x'_{p-1}\ot
Y, \\
&TF_p = S^{q-1}(\one\ot 1\ot X\wg Y);
\end{align*}

for $p=2q+1,$
\begin{align*}
\tau_p = &\one\ot \x_p\ot 1,\quad \quad
\tau_p' = \one\ot \x_p'\ot 1,  \\
&\s_p=\one\ot \x_{p-1} \ot X- \frac{1}{|\x_{p-1}|}\one\ot X
\xi_{p-1}\ot Y, \\
&\s_p' = \one\ot \x_{p-1}'\ot X- \frac{1}{|\x_{p-1}|}\one\ot X\x'_{p-1}\ot
Y, \\
&GV_p = S^{q}(\one\ot \xi_1\ot1).
\end{align*}
\end{proposition}
\begin{proof}
Evidently, $HC^0(\Hc_1, \Cb_\d) =  H^0(\Hc_1, \Cb)$. From the
bicomplex \eqref{UF+}, one sees that the vertical boundary $\vB$  is
given by the Lie algebra homology boundary of $\Fg$ with
coefficients in $\Cb_\d\ot \Fc^{\ot \ast}$, where $\Fg $ acts by
\eqref{g-action}.  We thus need to compute the homology of the
derived complex $\{ H^\bullet (\Hc_1, \Cb_\d), {\p_\Fg} \}$.  Since
$\, {\p_\Fg}:H^1(\Hc_1,\Cb)\ra H^0(\Hc_1, \Cb) $ is given by
\begin{align*}
  {\p_\Fg}(\one\ot \x_1\ot 1)= {\p_\Fg}(\one\ot \x_1'\ot 1)=0,  \quad \quad {\p_\Fg}(\one\ot 1\ot Y)=\one\ot 1;
\end{align*}
it follows that $HC^1(\Hc_1, \Cb_\d)$ is $2$-dimensional, with the
{\em Godbillon-Vey  class} $\, \one\ot \x_1 \ot 1$ and {\em
Schwarzian class} $\, \one\ot \x_1'\ot 1$ as its basis.

Next, $\,  {\p_\Fg}:H^2(\Hc_1,\Cb)\ra H^1(\Hc_1, \Cb)$ is given by
\begin{align*}
   &{\p_\Fg}(\one\ot \x_2\ot 1)={\p_\Fg}(\one\ot \x_2'\ot 1)=0, \\
 &{\p_\Fg}(\one\ot
 1\ot Y)=1,\\
 & {\p_\Fg}(\one\ot \x_1'\ot Y)= \one\ot( 1-|\x_1'|)\x_1',\\
 & {\p_\Fg}(\one\ot 1\ot X\wg Y)=\one\ot
 1 \ot  X-\one\ot 1\ot X=0, \\
 &{\p_\Fg}(\one\ot \x_{1}\ot X- \frac{1}{|\x_{1}|}\one \ot X \xi_{1}\ot Y)=0\\
 & {\p_\Fg}(\one\ot \x_{1}'\ot X- \frac{1}{|\x_{1}|}\one\ot X
 \xi_{1}'\ot Y)=0.
 \end{align*}
Hence   $HC^2(\Hc_1, \Cb_\d)$ is $5$-dimensional and is generated by
 the  {\em  transverse fundamental class} $\, \one\ot 1\ot X\wg Y$, together
 with the classes
 \begin{align*}
 &\one\ot
 \x_2\ot  1,\;\;\quad \one\ot \x_2'\ot 1,\\
 & \one\ot \x_{1}\ot X- \frac{1}{|\x_{1}|}\one
 \ot X\xi_{1}\ot Y, \quad \;\;\one\ot \x_{1}'\ot X- \frac{1}{|\x_{1}|}\one
 \ot X\xi_{1}'\ot Y.
 \end{align*}
 Finally,  for $p\ge 3$, $\, {\p_\Fg}: H^p(\Hc_1,\Cb)\ra H^{p-1}(\Hc_1, \Cb)$ one has
 \begin{align*}
 &{\p_\Fg}(\one\ot \x_p\ot 1)=\p_\Fg(\one\ot \x_p'\ot 1)=0, \\
 &{\p_\Fg}(\one\ot \x_p\ot Y)= (1-|\x_p|)\one\ot \x_p, \\
 &{\p_\Fg}(\one\ot \x_p'\ot Y)=
 (1-|\x_p'|)\one\ot \x_p',\\
 &{\p_\Fg}(\one\ot \x_{p-2}\ot X\wg Y)=|\x_{p-2}|\one\ot \x_{p-2}\ot X -\one\ot X \x_{p-2}\ot Y, \\
&{\p_\Fg}(\one\ot \x_{p-2}'\ot X\wg Y)=|\x_{p-2}'|\one\ot \x_{p-2}' \ot X-\one\ot X \x_{p-2}'\ot Y, \\
 &{\p_\Fg}(\one\ot \x_{p-1}\ot X- \frac{1}{|\x_{p-1}|}\one\ot X \xi_{p-1}\ot Y)=0\\
 &{\p_\Fg}(\one\ot \x_{p-1}'\ot X- \frac{1}{|\x_{p-1}|}\one\ot X
 \xi_{p-1}'\ot Y)=0 ,
 \end{align*}
 whence the claimed result.
 \end{proof}

\end{document}